\documentclass[11pt]{article}

\RequirePackage{amsthm,amsmath,amsfonts,amssymb}
\usepackage[margin=1in]{geometry}
\usepackage{bbm}
\usepackage{latexsym}
\usepackage{graphicx}
\usepackage{color,wrapfig}
\usepackage{multirow}
\usepackage{bm,url,comment}
\usepackage{authblk}
\usepackage{mathtools,mathrsfs}
\usepackage{algorithmic}
\usepackage{bm}
\usepackage{array}
\usepackage{algorithm}
\usepackage{caption}
\usepackage{subcaption}
\usepackage{hyperref}

\theoremstyle{plain}
\newtheorem{theorem}{Theorem}[section]
\newtheorem*{theorem*}{Theorem}
\newtheorem{lemma}[theorem]{Lemma}
\newtheorem{assumption}[theorem]{Assumption}
\newtheorem{corollary}[theorem]{Corollary}
\newtheorem{proposition}[theorem]{Proposition}

\newtheorem{definition}[theorem]{Definition}
\theoremstyle{remark}

\newtheorem{remark}[theorem]{Remark}

\newcommand{\ub}{\mathbf{u}}

\newcommand{\db}{\mathbf{d}}

\newcommand{\vb}{\mathbf{v}}

\newcommand{\T}{\mathrm{T}}

\newcommand{\norm}[1]{\left\lVert#1\right\rVert}

\newcommand{\wt}{\widetilde}

\theoremstyle{plain} %plain, definition, remark

\newtheorem*{lemma*}{Lemma}

\newtheorem*{corollary*}{Corollary}

\newtheorem*{proposition*}{Proposition}

\newtheorem*{definition*}{Definition}
\theoremstyle{remark}

\newtheorem*{example*}{Example}

\newtheorem*{remark*}{Remark}
\newtheorem*{remarks*}{Remarks}

\DeclareMathOperator*{\argmin}{arg\,min}
\title{Data-Driven Matrix Recovery via Optimal Shrinkage and Spatially Resolved Singular Vector Denoising under High-Dimensional Separable Noise}

\author{Pei-Chun Su\thanks{Department of Applied and Computational Mathematics, Yale University, New Haven, USA. E-mail: pei-chun.su@yale.edu.}}

\date{}

\begin{document}
\maketitle

 \begin{abstract}
This paper develops a spatially resolved perturbation theory for singular vectors under high-dimensional separable noise and applies it to data-driven matrix recovery. In the asymptotic regime where the matrix dimensions are proportional and significantly larger than the signal rank, we derive exact leading-order variance formulas for the singular vector perturbation projected onto any spatial patch. The variance decomposes into a spatially non-uniform component governed by the local noise covariance and a spatially uniform component governed by the global noise level. These formulas provide the foundation for the \emph{extended optimal shrinkage and wavelet shrinkage} (e$\mathcal{OWS}$) algorithm, which recovers low-rank matrices satisfying a mixed H\"older condition. The pipeline begins with optimal shrinkage of singular values, then constructs coupled multiscale partition trees on the row and column spaces from the denoised estimate, generating a tensor Haar-Walsh wavelet basis. Spatially adaptive wavelet shrinkage is applied using data-driven, coefficient-level thresholds derived from the perturbation theory. We establish convergence rates that strictly improve upon both optimal shrinkage and wavelet shrinkage applied in isolation. Numerical simulations demonstrate reliable matrix recovery and accurate reconstruction of the underlying singular subspaces, including an application to fetal ECG extraction.
\end{abstract}
\section{Introduction}

The primary focus of this study is the denoising of a $p \times n$ data matrix $\widetilde{S}$, comprising $n \in \mathbb{N}$ noisy samples of dimension $p \in \mathbb{N}$. The observed data matrix is modeled additively as:
\begin{equation}\label{eq_model}
\widetilde{S} = S + Z \in \mathbb{R}^{p \times n}
\end{equation}
Here, $Z$ represents a noise-only random matrix governed by a specific dependence structure detailed subsequently. The matrix $S$ is the underlying signal matrix, which is assumed to be both low-rank and subject to a mixed H\"older continuity condition. Such additive splitting frameworks are pervasive throughout mathematical analysis and possess broad utility in applied domains, particularly within signal and image processing methodologies.

In the following, we briefly review two shrinkage-based approaches for recovering the underlying matrix $S$. Optimal shrinkage is designed to exploit low-rank structures, leveraging spectral decay to recover the dominant singular components. Conversely, we review a wavelet shrinkage scheme specifically tailored for the mixed H\"older condition, utilizing localized basis shrinkage to capture and preserve multi-scale features. While wavelet shrinkage is a versatile tool applicable to various regularity classes, the specific implementation here leverages the geometry of the data to achieve near-optimal recovery for the mixed H\"older case.

\subsection{Spiked Model and Optimal Shrinkage of Singular Values}\label{sec_os}
\subsubsection{Spiked Model}
Given the additive model in \eqref{eq_model}, assume $S$ is strictly low-rank, admitting a singular value decomposition (SVD) of the form $\sum_{i=1}^r d_i \mathbf{u}_i \mathbf{v}_i^\top$, where the rank $r \geq 1$ is assumed to be significantly smaller than both $p$ and $n$. The vectors $\mathbf{u}_i \in \mathbb{R}^p$ and $\mathbf{v}_i \in \mathbb{R}^n$ represent the left and right singular vectors, respectively, which prescribe the signal structure, while the associated singular values $d_i > 0$ quantify the signal strength and may scale with $n$. For brevity, we will refer to the pair $\mathbf{u}_i$ and $\mathbf{v}_i$ as the $i$-th signal vectors and $d_i$ as the $i$-th signal strength. The asymptotic \emph{high-dimensional} regime is defined such that $p = p(n)$ grows proportionally with $n$, satisfying $p/n \to \beta \in (0,\infty)$ as $n \to \infty$. This model represents a generalization of the traditional high-dimensional spiked covariance model \cite{bai2008central, bai2012sample, baik2006eigenvalues}.\\

\subsubsection{Optimal Shrinkage of Singular Values}
When $S$ is inherently low-rank, a weighting approach applied to the SVD of $\widetilde{S}$---widely known as singular value shrinkage---has been actively studied for the recovery of $S$. To the best of our knowledge, this approach was first introduced in \cite{eckart1936approximation, mirsky1960symmetric, golub1965calculating}. The central idea involves selecting an appropriate, often nonlinear, shrinkage function $\varphi: [0, \infty) \rightarrow [0,\infty)$ to construct the estimator:
\begin{equation}\label{eq_shrink}
\widehat{S}_{\varphi} = \sum_{i=1}^{p \wedge n} \varphi(\widetilde{\lambda}_i) \widetilde{\mathbf{u}}_i \widetilde{\mathbf{v}}_i^\top
\end{equation}
where $\widetilde{\lambda}_1 \geq \widetilde{\lambda}_2 \geq \cdots \geq \widetilde{\lambda}_{p \wedge n} \geq 0$ denote the singular values of $\widetilde{S}$ (i.e., the square roots of the common eigenvalues of $\widetilde{S}\widetilde{S}^\top$ and $\widetilde{S}^\top\widetilde{S}$), and $\{\widetilde{\mathbf{u}}_i\}_{i=1}^p$ and $\{\widetilde{\mathbf{v}}_i\}_{i=1}^n$ denote the empirical left and right singular vectors of $\widetilde{S}$, respectively. The function $\varphi$ is termed a shrinker. By employing a loss function $L_n:\mathbb{R}^{p\times n}\times \mathbb{R}^{p\times n}\to \mathbb{R}_+$ to quantify the discrepancy between $\widehat{S}_{\varphi}$ and $S$, the associated \emph{optimal shrinker}, provided it exists, is defined as:
$$ \varphi^* := \arg\min_{\varphi} \lim_{n\to \infty} L_n(\widehat{S}_{\varphi}, S) $$
Commonly utilized loss functions include the Frobenius norm and the operator norm of the difference $\widehat{S}_{\varphi} - S$. This methodological approach is referred to as \textit{Optimal Shrinkage} ($\mathcal{OS}$) in previous literature \cite{donoho2018optimal, GD}. Under the spiked model framework, numerous studies \cite{nadakuditi2014optshrink, GD, shabalin2013reconstruction, leeb2021optimal, leeb2021matrix, leeb2021rapid, donoho2023screenot,ding2024eigenvector, zhang2024matrix,lin2024eigenvector, su2025data} have investigated this matrix denoising problem and derived the associated $\mathcal{OS}$ under various assumptions regarding the noise matrix $Z$.

In \cite{su2025data}, we analyzed a specific scenario characterized by the separable covariance structure:
\begin{equation}\label{colored and dependent noise model}
Z = \mathcal{A}^{1/2}\mathcal{X} \mathcal{B}^{1/2}
\end{equation}
where $\mathcal{X}$ is a random matrix with independent entries satisfying specific moment conditions, and $\mathcal{A}$ and $\mathcal{B}$ are deterministic, positive-definite matrices of dimensions $p \times p$ and $n \times n$, respectively, which dictate the color and row-column dependence structure of the noise. Building upon the foundational results in \cite{benaych2011fluctuations, limitpaper}, we derived \emph{eOptShrink}, a fully data-driven algorithm designed to estimate the effective rank and the spectral density distribution of $Z$. This formulation yields a precise estimation of the asymptotic bias embedded in the singular values and singular vectors, leading to superior optimal shrinkers when compared to existing state-of-the-art approaches \cite{GD, donoho2023screenot, nadakuditi2014optshrink}.

The application of $\mathcal{OS}$ and eOptShrink spans diverse fields, including diffusion magnetic resonance imaging denoising \cite{veraart2016denoising}, fetal electrocardiogram (fECG) extraction from trans-abdominal maternal ECG (ta-mECG) \cite{su2019recovery, su2025data}, ECG T-wave quality evaluation \cite{su2020robust}, seismic data enhancement \cite{anvari2018enhancing}, otoacoustic emission signal denoising \cite{liu2021denoising}, and the removal of both stimulation artifacts from intracranial electroencephalograms (EEG) \cite{alagapan2019diffusion} and cardiogenic artifacts from standard EEG \cite{chiu2022get}.
Further details regarding $\mathcal{OS}$ and eOptShrink are provided in Section \ref{sec_eoptshrink}.
\subsection{Mixed H\"older Matrix Recovery via Wavelet Shrinkage}
\subsubsection{Wavelet Analysis}
Let us briefly review some foundational concepts in wavelet theory, focusing specifically on the framework of a multiresolution analysis within the space $L^2(\mathbb{R})$. This framework begins with an initial scaling function, $\phi(x)$, from which a family of functions is generated via dyadic dilation and integer translation:
\begin{equation}
\phi_{j,k}(x) = 2^{-j/2}\phi(2^{-j}x - k)
\end{equation}
Let the subspace $V_j$ denote the closed linear span of the set $\{ \phi_{j,k} \}$ across all integers $k$. Provided the function $\phi$ meets specific regularity criteria, these subspaces form a nested sequence such that $V_j \subsetneq V_{j-1}$. Equivalently, the base function $\phi(x)$ can be expressed as a linear combination of the scaled and shifted functions $\phi(2x - k)$. Furthermore, the union of all such subspaces spans the entirety of $L^2(\mathbb{R})$. A collection of subspaces $V_j$ satisfying these structural properties constitutes a ``multiresolution analysis'' of $L^2(\mathbb{R})$. Conceptually, each individual subspace $V_j$ isolates the behavior of a function at a specific dyadic scale or resolution level.

The core of wavelet analysis entails examining $W_j$, defined as the orthogonal complement of the subspace $V_j$ relative to $V_{j-1}$. Given this multiresolution structure, it is possible to construct a distinct function $\psi(x)$ such that its integer translations form an orthonormal basis for $W_0$. Consequently, the scaled and translated family of functions $\psi_{j,k}(x) = 2^{-j/2}\psi(2^{-j}x - k)$ will completely span the subspace $W_j$. In this context, $\phi$ is traditionally referred to as the scaling function (or ``father wavelet''), while $\psi$ is designated as the ``mother wavelet.''

The Haar wavelet system serves as the foundational illustration of these concepts \cite{gavish2010multiscale,coifman2011harmonic,ankenman2014geometry, ankenman2018mixed}. For the Haar system, the scaling function $\phi$ is defined as the indicator function over the unit interval $[0, 1]$. The corresponding mother wavelet is defined as the difference of indicator functions: $\chi_{[0, 1/2]} - \chi_{[1/2, 1]}$. The associated subspace $V_j$ is thus constructed as the span of indicator functions over the dyadic intervals $[2^{-j}k, 2^{-j}(k+1)]$ for any integer $k$. \\

\subsubsection{H\"older and Mixed H\"older Space}
The variation of a function $f$ defined on a metric space $(X, d)$ is frequently quantified using the H\"older semi-norm. For a given smoothness parameter $\alpha > 0$, this is defined by the supremum:
\begin{equation}
\sup_{x \neq y} \frac{|f(x) - f(y)|}{d(x, y)^\alpha}
\end{equation}
When evaluating functions defined on the real line $\mathbb{R}$, this functional space remains mathematically non-trivial strictly within the regime $0 < \alpha \leq 1$. For parameter values strictly less than $1$, the space of H\"older continuous functions exhibits highly desirable algebraic and analytic properties. A primary advantage is that a function's H\"older norm can be definitively characterized by the decay rate of its wavelet coefficients. Assuming a sufficiently regular wavelet basis $\{ \psi_{j,k} \}$ for $\mathbb{R}^n$---where the indices $j \in \mathbb{Z}$ and $k \in \mathbb{Z}^n$ correspond to the dyadic scale $2^{-j}$ and spatial translation, respectively---the H\"older norm is equivalent in magnitude to the following supremum:
\begin{equation}
\sup_{j,k} 2^{j(\alpha + 1/2)} |\langle f, \psi_{j,k} \rangle|
\end{equation}
This equivalence dictates that the ratio between this supremum and the classical H\"older norm is bounded from above and below by finite, absolute constants that are entirely independent of $f$. Conceptually, the inner products $\langle f, \psi_{j,k} \rangle$ quantify the localized variation of the function across all available scales.

While standard H\"older spaces are formulated within the context of a single metric space, numerous real-world applications involve data domains that are more accurately modeled as the product of multiple metric spaces. Bivariate or multi-dimensional datasets naturally emerge in various analytical contexts. For instance, analyzing transposable data matrices requires the simultaneous examination of underlying row and column structures. Co-clustering algorithms similarly seek meaningful, concurrent groupings across both dimensions of a data matrix, aligning perfectly with this multi-space mathematical perspective.

We focus our attention on the Cartesian product of two specific metric spaces, denoted as $(X, d_X)$ and $(Y, d_Y)$, to evaluate the regularity of a matrix function $S$ operating on $X \times Y$. For any $\alpha > 0$, a function $f$ is said to satisfy a \textit{mixed H\"older}($\alpha$) condition with constant $L = L(f, \alpha)$ if the mixed difference quotient is bounded by $L$:
\begin{equation}\label{eq_mixH}
\sup_{x \neq x', y \neq y'} \frac{|f(x,y) - f(x,y') - f(x',y) + f(x',y')|}{d_X(x,x')^\alpha d_Y(y,y')^\alpha} \leq L
\end{equation}
The mixed H\"older condition rigorously characterizes the coupled variation of $f$ across both domains \cite{ankenman2018mixed}. In Euclidean settings, this mixed difference quotient is not rotationally invariant, meaning its magnitude fluctuates significantly depending on the chosen coordinate frame. However, in the realm of applied data analysis, coordinate axes represent fundamental, intrinsic dimensions of the dataset (e.g., words versus documents in text mining, or time versus frequency within a spectrogram). In these contextually rich scenarios, the axes possess independent significance, making it mathematically logical to employ axis-dependent metrics such as the mixed H\"older norm.\\

\subsubsection{Tree Structure}
While the function spaces introduced previously were formulated using Euclidean geometry, the current work shifts to a framework built upon an entirely different structural foundation: abstract tree metrics. Let $X = \{x_1,\ldots,x_p\}$ and $Y = \{y_1,\ldots,y_n\}$ be discrete spaces of cardinalities $p$ and $n$, respectively. Conceptually, $X$ represents the row indices and $Y$ represents the column indices of a data matrix $S$, such that:
\begin{equation}
S(x_i,y_j) := S_{ij}
\end{equation}
denotes the $ij$-th entry of $S$. Without loss of generality, we assume both spaces are endowed with probability measures, such that the total measure of each space equals $1$. For any arbitrary subset $E$, its measure is denoted by $|E|$. We assume $X$ and $Y$ are equipped with multiscale partition trees and associated tree metrics $(\mathcal{T}_X, d_{\mathcal{T}_X})$ and $(\mathcal{T}_Y, d_{\mathcal{T}_Y})$, defined as follows:

\begin{definition}[Multiscale Partition Tree Structure $\mathcal{T}_X$]\label{def_tree}
Let $X = \{x_i\}_{i=1}^N$ be a set of sample points. A multiscale partition tree $\mathcal{T}_X = \left\{ I_k^\ell \right\}_{0 \le \ell \le L, k \in \mathcal{C}_k^\ell }$ is a hierarchical decomposition of the index set $\{1, \dots, N\}$ into disjoint subsets $I_k^\ell$, referred to as \emph{folders}. Here, $\ell = 0, \dots, L$ indexes the hierarchical level, $k$ indexes the individual folders at level $\ell$, and $\mathcal{C}_k^\ell$ denotes the index set of the child folders of $I_k^\ell$. The root of the tree is the universal set $I_0^0 = \{1, \dots, N\}$. The deepest level, $L$, consists of singleton folders $\{x_i\}$ for each individual element in $X$. Each folder $I_k^\ell$ at level $\ell$ is uniquely defined as the disjoint union of its children at level $\ell+1$:
\begin{equation}
I_k^\ell = \bigcup_{m \in \mathcal{C}_k^\ell} I_m^{\ell+1}
\end{equation}
Folders residing at the same level are strictly mutually exclusive: $I^\ell_k \cap I^\ell_{k'} = \emptyset$ for $k \neq k'$.
\end{definition}

\begin{definition}[\textbf{Tree metric}]\label{def_mixed_Holder}
Letting $I_{x,x'}$ denote the minimal enclosing folder within a tree $\mathcal{T}_X$ that contains both elements $x$ and $x'$, the associated tree distance $d_{\mathcal{T}_X}(x, x')$ is explicitly defined as:
\begin{equation}\label{eq_tree_d}
d_X(x, x') = \begin{cases} |I_{x,x'}|, & \text{if } x \neq x' \\ 0, & \text{if } x = x' \end{cases}
\end{equation}
Throughout this paper, we require the partition tree to be \emph{balanced}, establishing that there exist bounding constants $B_L$ and $B_U$ such that:
\begin{equation}\label{eq_balance}
0 < B_L \leq \frac{|I^{\ell+1}_m|}{|I^{\ell}_k|} \leq B_U < 1
\end{equation}
for any arbitrary folder $I_{k}^\ell \in \mathcal{T}_X$ and its direct child folder $I_m^{\ell+1}$.
\end{definition}

Prior literature \cite{ankenman2014geometry, coifman2011harmonic, ankenman2018mixed} establishes that wavelet-based characterizations of both standard and mixed H\"older spaces can be successfully translated to function spaces defined on individual partition trees and their Cartesian products. In this discrete paradigm, standard continuous wavelets are replaced with specialized Haar wavelets tailored to the geometry of the trees. Functions exhibiting mixed H\"older regularity with respect to these tree metrics demonstrate high compressibility, enabling precise mathematical reconstruction using a remarkably sparse selection of their tensor Haar coefficients. The concept of sparse grid reconstruction, standard in Euclidean settings, can similarly be adapted to this discrete framework.

In this work, we utilize the extended generalized Haar-Walsh transform (eGHWT) in conjunction with a best-basis selection algorithm \cite{saito2022eghwt} to extract the optimal wavelet basis tailored to the constructed discrete partition trees. Compared to using the Haar system alone, the eGHWT searches over a substantially richer dictionary of orthonormal bases — encompassing both Haar-like and Walsh-like functions at each level of the tree — enabling the best-basis algorithm to select the representation that achieves the sparsest expansion for the signal at hand \cite{saito2022eghwt}. Formal definitions of the Walsh system and the best-basis methodology are provided in Appendix \ref{sec_Walsh}.\\

\subsubsection{Wavelet Shrinkage}
The fundamental objective of statistical estimation here is the retrieval of the underlying signal $S$ from the corrupted observations $\widetilde{S}$ in \eqref{eq_model}. Successfully isolating this signal strictly requires structural assumptions. Parametric approaches constrain the signal to families governed by finite parameters, whereas nonparametric approaches impose less stringent conditions, often assuming the function merely belongs to a designated functional class.

To evaluate the efficacy of an estimator $\widehat{S}$, the minimax framework assesses the expected loss, or risk, commonly via the squared Euclidean metric $L(S, \widehat{S}) = \|S - \widehat{S}\|_2^2$. Because the observed data and the resulting estimator are stochastic, the risk fundamentally depends on the unknown true signal $S$:
\begin{equation}
R(S, \widehat{S}) = \mathbb{E}_S L(S, \widehat{S})
\end{equation}
Assuming $S$ resides within a specific functional family $\mathcal{F}$, the minimax paradigm identifies an optimal estimator $\widehat{S}^*$ that minimizes the maximum risk across all possible functions in that class:
\begin{equation}
\widehat{S}^* = \arg\min_{\widehat{S}} \sup_{S \in \mathcal{F}} R(S, \widehat{S})
\end{equation}
Rather than pursuing the exact minimax risk, robust statistical analyses often focus on finding estimators whose risk achieves the minimax rate, differing from the optimal minimax risk by at most a constant $C$ that is independent of both the sample size and the parameters defining the function space $\mathcal{F}$. 

While explicit linear estimators can achieve this rate when the exact parameters of a functional space are known a priori, this is rarely true in practice. There is a substantial need for adaptive estimators capable of near-optimal performance across an entire spectrum of functional spaces without requiring explicit knowledge of their underlying parameters.

\textit{Wavelet Shrinkage} ($\mathcal{WS}$) techniques introduced by Donoho and Johnstone provide a robust solution to this adaptability challenge \cite{donoho1994ideal, donoho1995adapting, donoho1996neo, donoho1995noising, donoho1995wavelet, donoho1998minimax}. These methods apply a thresholding mechanism that drives the empirical wavelet coefficients of the noisy signal toward zero using a non-linear soft-thresholding function:
\begin{equation}\label{eq_ws_shrinker}
\eta_t(x) = \text{sgn}(x)(|x| - t)_+
\end{equation}
The soft-thresholding function $\eta_t$ is not the minimax-optimal shrinker, but it is widely adopted due to its simplicity, computational efficiency, and near-optimal theoretical guarantees: the resulting estimator achieves an error rate matching the minimax rate up to a logarithmic factor $\log^a(n)$, where $a < 1$ is a parameter intrinsic to the functional class $\mathcal{F}$.

In this work, we apply the wavelet shrinkage procedure specifically using tree-based wavelet bases from eGHWT. The constructed tree structures $\mathcal{T}_X$ and $\mathcal{T}_Y$ define specialized eGHWT wavelets $\{\omega_I\}_{I \in \mathcal{T}_X}$ on $X$ and $\{\omega_J\}_{J\in \mathcal T_Y}$ on $Y$ \cite{coifman2011harmonic, gavish2012sampling, gavish2010multiscale}. The tensor product of these systems provides a basis for functions defined on $X \times Y$:
\begin{equation}\label{eq_tensor_haar}
\Phi_{IJ}(x, y) = \omega_I(x)\omega_J(y)
\end{equation}
These tensor functions form an orthonormal basis for any function on $X \times Y$. The discrete wavelet coefficient of $S$ projected onto the basis $\Phi_{IJ}$ is given by:
\begin{equation}\label{eq_ws_coeff}
\langle S, \Phi_{IJ} \rangle := \sum_{x \in I, y \in J} S(x,y)\Phi_{IJ}(x,y)
\end{equation}
To construct an estimator of $S$ from $\widetilde{S}$, the empirical matrix $\widetilde{S}$ is expanded into a two-dimensional wavelet series:
\begin{equation}\label{eq_ws_expand}
\widetilde{S}(x, y) = \sum_{\Phi_{IJ}} \langle \widetilde{S}, \Phi_{IJ} \rangle \Phi_{IJ}(x, y)
\end{equation}
When the elements of $Z$ are modeled as i.i.d. zero-mean Gaussian variables with variance $\sigma^2$, the observed wavelet coefficients follow a normal distribution centered precisely at the true coefficients. Modified coefficients are acquired by applying the soft-shrinker $\eta_t$ to these empirical observations, yielding the estimator $\widehat{C}_{\Phi_{IJ}} = \eta_t(\langle \widetilde{S}, \Phi_{IJ} \rangle)$. The resulting matrix estimate is subsequently constructed as:
\begin{equation}
\widehat{S}_{\eta_t}(x, y) = \sum_{\Phi_{IJ}} \widehat{C}_{\Phi_{IJ}} \Phi_{IJ}(x, y)
\end{equation}
Further details regarding the derivation of the optimal shrinker $\eta^*_t$ when $S$ exhibits a mixed H\"older($\alpha$) condition, along with associated theoretical guarantees, are provided in Section \ref{sec_ws}.
\subsection{Motivations}
Our goal is to recover the signal matrix $S$ from the noisy observation model \eqref{eq_model}, where $S$ is assumed to be both low-rank and satisfies a mixed H\"older regularity condition. This combined structural assumption arises naturally in many real-world datasets where the signal is governed by a small number of latent factors that vary smoothly along both the row and column dimensions. For example, in diffusion magnetic resonance imaging \cite{veraart2016denoising}, the signal matrix captures diffusion-weighted measurements across spatial voxels and gradient directions. The matrix is inherently low-rank due to the underlying tissue microstructure and exhibits smooth variation across neighboring voxels and diffusion encodings. In fetal electrocardiogram extraction from trans-abdominal maternal ECG \cite{su2019recovery, su2025data}, the signal is low-rank, driven by a few cardiac sources, and varies smoothly across time (columns) and cardiac cycles (rows). In hyperspectral imaging \cite{he2018hyperspectral}, the matrix $S$ records reflectance values across spatial pixels (rows) and spectral bands (columns), where a few material signatures contribute to a low-rank structure that varies smoothly across neighboring pixels and wavelengths.

When the row and column indices possess a natural sequential ordering, such as consecutive time steps or adjacent spatial locations, the mixed H\"older condition may hold with respect to this standard ordering. However, even in such settings, a sequential ordering may not capture the full structure of the data. For instance, in cardiac signal analysis \cite{su2019recovery, su2025data}, cardiac cycles that are far apart in time may exhibit highly similar morphology, representing non-local neighbors that a simple sequential ordering would fail to exploit. Similarly, in seismic data \cite{anvari2018enhancing}, spatially distant sensors may record similar waveforms due to shared geological structures. In these cases, learning the partition trees $\mathcal{T}_X$ and $\mathcal{T}_Y$ from the data allows the tree metric to group such non-local neighbors together, yielding a more compact representation under the mixed H\"older condition than any fixed sequential ordering. More broadly, when the intrinsic geometry is entirely unknown or non-sequential, tree estimation becomes indispensable. For example, in single-cell RNA sequencing \cite{linderman2022zero}, cells and genes lack a canonical ordering but exhibit hierarchical cluster structure, and the signal may exhibit mixed H\"older regularity with respect to the induced tree metrics once appropriate partition trees are learned.

These considerations expose the complementary limitations of spectral and spatial denoising methods when applied in isolation to high-dimensional, structured noise regimes.
Global approaches, such as optimal singular value shrinkage (eOptShrink), effectively regularize the singular values of a matrix but apply a uniform spectral correction that inherently fails to capture spatially or temporally localized patterns. More critically, while they correct singular values, they cannot correct the rotational bias present in the empirical singular vectors $\widetilde{\mathbf{u}}_i$ and $\widetilde{\mathbf{v}}_i$. This rotational perturbation is severely exacerbated in high-dimensional settings, distorting the fundamental signal representation and directly degrading the performance of downstream dimensionality reduction tasks, such as classification and clustering. 
Conversely, standard wavelet shrinkage ($\mathcal{WS}$) excels at adapting to local structures at multiple scales. However, current spatial methods generally presume i.i.d. Gaussian noise, failing to account for the complex, separable covariance structures common in real-world data as modeled in \eqref{colored and dependent noise model}. Furthermore, standard methods derive the underlying partition trees directly from the raw, corrupted data matrix (see Appendix \ref{sec_tree}). Under heavy noise regimes, this heuristic breaks down entirely; the noise disrupts the hierarchical clustering, destroying the geometric foundation required for an effective wavelet basis.
\subsection{Contributions}

The central theoretical contribution of this work is a spatially resolved perturbation theory for singular vectors under structured noise. We show that under the separable noise covariance model $Z = \mathcal{A}^{1/2}\mathcal{X}\mathcal{B}^{1/2}$, the variance of the singular vector perturbation projected onto any spatial patch is governed by the local noise energy on that patch. This is a fundamentally local quantity: different spatial regions experience different amounts of eigenvector distortion, dictated entirely by the noise covariance structure. To our knowledge, this spatially resolved characterization is new. Classical perturbation bounds \cite{davis1970rotation, wedin1972perturbation} provide only global worst-case control over singular subspace distances. Recent entry-wise analyses \cite{abbe2020entrywise} establish exact formulas for singular vector perturbations, but only under white noise, where the variance is spatially uniform. Our result generalizes this to separable colored noise, where the covariance structure makes the variance spatially non-uniform, and extends the analysis from individual entries to arbitrary spatial projections.

Building upon the established convergence rates for eOptShrink in \cite{su2025data}, our theoretical contributions proceed in four stages. First, we derive the first-order asymptotic perturbation expansion for singular vectors under the separable noise model. Second, we establish the exact leading-order formulas for the local variance of these perturbations projected onto arbitrary spatial patches, yielding the spatially resolved theory described above. The precise formulas are presented in Section~\ref{sec_result}. Third, leveraging this local variance characterization, we obtain precise variance formulas for the eGHWT coefficients of the residual error, which enable the construction of scale-adaptive thresholds for each basis element $\Phi_{IJ}$. Fourth, combining these results, we establish risk bounds for the e$\mathcal{OWS}$ estimator that quantify how the combined procedure improves upon standalone $\mathcal{OS}$ and $\mathcal{WS}$.

These theoretical results have an immediate algorithmic consequence. Existing optimal shrinkage methods correct singular values but leave the empirical singular vectors $\widetilde{\mathbf{u}}_i$ and $\widetilde{\mathbf{v}}_i$ untouched, despite the fact that these vectors carry substantial, spatially heterogeneous bias in high-dimensional settings. Because our theory precisely quantifies this bias at each scale and location, it enables scale-adaptive wavelet shrinkage that targets and removes the singular vector perturbation locally, achieving what uniform spectral methods cannot. To realize this idea, we propose the e$\mathcal{OWS}$ pipeline, which integrates three components into a unified methodology with rigorous theoretical guarantees: eOptShrink for singular value correction, the Questionnaire algorithm for learning partition trees from the denoised estimate, and scale-adaptive wavelet shrinkage via the eGHWT best basis for singular vector correction. The resulting estimator integrates optimally denoised singular values with geometrically corrected singular vectors, yielding a substantially superior approximation of the original signal matrix $S$.

\subsubsection{The e$\mathcal{OWS}$ Pipeline}
We focus on the regime above the BBP transition, where the signal singular vectors are detectable but carry substantial bias that spectral methods alone cannot remove. Guided by the above theoretical analysis, we propose an enhanced recovery methodology for signals $S$ that are simultaneously low-rank and mixed H\"older continuous, operating under a separable noise covariance model. Our method, termed \emph{extended optimal shrinkage and wavelet shrinkage} (e$\mathcal{OWS}$), combines three components into a unified pipeline.

\textit{Step 1: Singular value correction via eOptShrink.} The algorithm first computes an initial estimator via eOptShrink:
$\widehat{S}_{\widehat{r}^+} = \sum_{i=1}^{\widehat{r}^+} \widehat{d}_{e,i}\,\widetilde{\mathbf{u}}_i \widetilde{\mathbf{v}}_i^\top$, where $\widehat{r}^+$ is the estimated effective rank and $\widehat{d}_{e,i}$ is the corrected singular value. These estimators, defined formally in \eqref{eq_a1a2}, are shown to be consistent with explicit convergence rates in \cite{su2025data}. While eOptShrink effectively corrects the singular values, the empirical singular vectors $\widetilde{\mathbf{u}}_i$ and $\widetilde{\mathbf{v}}_i$ remain biased. Writing $\widetilde{\mathbf{u}}_i = \mathbf{u}_i + \Delta\mathbf{u}_i$ and $\widetilde{\mathbf{v}}_i = \mathbf{v}_i + \Delta\mathbf{v}_i$, the residual error between $\widehat{S}_{\widehat{r}^+}$ and the true truncated signal $S_{r^+} = \sum_{i=1}^{r^+} d_i\,\mathbf{u}_i \mathbf{v}_i^\top$ is:
\begin{align}\label{eq_error_vec}
\widehat{S}_{\widehat r^+} - S_{r^+} & \simeq \sum_{i=1}^{r^+} d_i \left( {\widetilde{\mathbf{u}}_i\widetilde{\mathbf{v}}_i^\top} - \mathbf{u}_i\mathbf{v}_i^\top \right) \nonumber \\
& = \sum_{i = 1}^{r^+} d_i \left( {\mathbf{u}_i\Delta\mathbf{v}_i^\top} + {\Delta\mathbf{u}_i\mathbf{v}_i^\top} + {\Delta\mathbf{u}_i \Delta\mathbf{v}_i^\top} \right)\,,
\end{align}
where $\simeq$ denotes asymptotic equality in the sense that the difference between the two sides vanishes as $n \to \infty$. The local variance theory developed in Section~\ref{sec_result} shows that this residual is spatially heterogeneous, with its magnitude on each patch determined by the local noise covariance. Steps 2 and 3 exploit this structure to correct the singular vector perturbation.

\textit{Step 2: Tree learning via the Questionnaire algorithm.} We assume there exist intrinsic partition tree structures, $\mathcal{T}_X$ and $\mathcal{T}_Y$, under which the signal matrix satisfies the mixed H\"older condition with respect to the induced tree metrics. In practical scenarios where these trees are not known \emph{a priori}, we construct them from the denoised estimate $\widehat{S}_{\textup{eOptShrink}}$ using the Questionnaire algorithm \cite{ankenman2014geometry, coifman2011harmonic} (detailed in Appendix \ref{sec_question}), where $\widehat{S}_{\textup{eOptShrink}}$ is the estimator of $S$ via eOptShrink. The Questionnaire algorithm exploits the coupled geometry between the rows and columns of a matrix to iteratively build hierarchical partition trees on both index sets. The key idea is that the structure of one dimension informs the organization of the other. Concretely, starting from an initial affinity on the column set $X$ constructed via a Gaussian kernel, a partition tree $\mathcal{T}_X$ is built by recursive spectral bipartitioning using the Fiedler vector of the associated graph Laplacian. This tree on $X$ then induces a multiscale metric on the row set $Y$: two rows $y_1, y_2 \in Y$ are considered similar if their corresponding row vectors $S(y_1, \cdot)$ and $S(y_2, \cdot)$ are close in the Earth Mover's Distance with respect to the tree metric on $X$. This dual affinity is used to construct a partition tree $\mathcal{T}_Y$ on $Y$, which in turn induces a refined metric on $X$, and the process repeats. Through this alternating refinement, the algorithm converges to a pair of trees $\mathcal{T}_X$ and $\mathcal{T}_Y$ that jointly capture the intrinsic hierarchical structure of both dimensions. 
By operating on the denoised estimate $\widehat{S}_{\textup{eOptShrink}}$ rather than the raw noisy matrix $\widetilde{S}$, the learned trees are substantially more reliable, particularly under heavy noise regimes where trees constructed from raw data would be severely distorted. A precise quantification of this improvement in terms of the entry-level variance reduction is provided after in Remark~\ref{rmk_variance_structure}(e).

The partition trees $\mathcal{T}_X$ and $\mathcal{T}_Y$ each generate a rich dictionary of orthonormal wavelet bases, rather than a single fixed basis. Each basis element is associated with a tile in a space-frequency plane: its spatial localization is determined by a node in the tree, while its frequency localization is determined by the oscillation index within that node. Since our signal matrix $S$ is two-dimensional, the tensor product of the row and column dictionaries gives rise to a collection of orthonormal bases indexed by tiles in a four-dimensional space-frequency domain: two spatial dimensions (one for $X$, one for $Y$) and two corresponding frequency dimensions. The eGHWT, combined with a best basis selection algorithm, identifies the optimal tiling $\mathcal{B}^*$ of this four-dimensional domain that yields the sparsest representation of the signal among all admissible tilings (see Appendix~\ref{sec_Walsh} for formal definitions). The resulting optimal basis consists of tensor elements $\Phi_{IJ} = \omega_I \omega_J$, where $\omega_I$ and $\omega_J$ are the best basis functions selected from the dictionaries on $\mathcal{T}_X$ and $\mathcal{T}_Y$, respectively.

\textit{Step 3: Singular vector correction via scale-adaptive wavelet shrinkage.} Using the eGHWT best basis $\{\Phi_{IJ}\}$, we apply wavelet shrinkage to the initial estimate $\widehat{S}_{\widehat r^+}$ with thresholds that adapt to the noise variance $\sigma_{\Phi_{IJ}}$ at each scale and basis element. The local variance formulas from Section~\ref{sec_result} provide the theoretical foundation for these scale-adaptive thresholds, ensuring that the shrinkage is calibrated to the spatially varying noise energy rather than relying on a single global noise parameter $\sigma$, yielding a superior approximation of the original signal matrix $S$.

\subsubsection{Applications and Experiments}
The practical relevance of optimal shrinkage and wavelet shrinkage has been demonstrated across a wide range of domains, as discussed in the Motivations section. To validate the e$\mathcal{OWS}$ pipeline in settings representative of these applications, we conduct comprehensive numerical experiments in Section \ref{sec_numerical}. First, we evaluate the method on simulated kernel matrices arising from acoustic wave propagation and sinusoidal models under three progressively more challenging noise models: white noise, separable covariance with spectral gaps, and separable covariance with oscillatory eigenvalue profiles. These kernel matrices naturally exhibit both low-rank and mixed H\"older structure, and are representative of interaction models widely used in computational electromagnetics and acoustics, such as wave scattering, antenna radiation, and harmonic signal analysis \cite{rokhlin1985rapid, greengard1987fast}. Second, we demonstrate the practical impact of e$\mathcal{OWS}$ on the problem of single-channel fetal ECG extraction from trans-abdominal maternal ECG recordings, a real-world application where optimal shrinkage has been shown to be a critical component \cite{su2019recovery, su2025data}. In this setting, the noise is neither white nor independent across segments, making it an ideal testbed for our scale-adaptive approach. We include these substantive real-data applications to demonstrate that the theoretical improvements translate into measurable gains in practice.

The remainder of the paper is organized as follows. Section \ref{sec_pre} reviews the requisite background for $\mathcal{OS}$ and $\mathcal{WS}$. Section \ref{sec_result} details the theoretical results concerning the first-order perturbation of singular vectors and the precise variance calculations for the eGHWT coefficients under a separable noise covariance model. Section \ref{sec_alg} formalizes the e$\mathcal{OWS}$ algorithm alongside its theoretical risk guarantees. Section \ref{sec_numerical} presents the numerical evaluations.
\section{Preliminaries}\label{sec_pre}
\subsection{Required Assumptions}

We begin by introducing the key analytical tools needed to state our assumptions and main results.

\begin{definition}[\textbf{Empirical spectral Distribution}]
    Define the empirical spectral distribution of an $n\times n$ symmetric matrix $H$ as 
    \begin{equation}
        \pi_H:= \frac{1}{n}\sum_{i=1}^n\delta_{\ell_i},
    \end{equation}
    where $\ell_1 \geq \ell_2 \geq \ldots\geq\ell_n $ are the eigenvalues of $H$ and $\delta$ is the Dirac delta measure.
\end{definition}

In the proportional regime where $p/n \to \beta$, the empirical spectral distributions of $ZZ^\top$ and $Z^\top Z$ converge to deterministic limiting measures, whose Stieltjes transforms play a central role in characterizing the asymptotic behavior of the noise.

\begin{definition}[\bf Integral transforms]
Denote the probability measures for the eigenvalues of $ZZ^\top$ and $Z^\top Z$ as $\mu_{1c}$ and $\mu_{2c}$, respectively, in the limit $n\to \infty$ with $p/n \to \beta$. Let $\rho_{1c}$ and $\rho_{2c}$ denote the densities of the absolutely continuous parts of $\mu_{1c}$ and $\mu_{2c}$, respectively. Denote $\lambda_+$ as the right edge of the support of $\rho_{1c}$ and $\rho_{2c}$.
Denote the $D$-transform as 
\begin{equation}\label{eq_defnt}
 \mathcal{T}(x):=xm_{1c}(x) m_{2c}(x)\,,
 \end{equation}
where $m_{1 c}$ and $m_{2 c}$ are the Stieltjes transforms
\begin{equation}\label{eq_defm1m2}
    m_{1 c}(x) = \int_0^{\lambda_+} \frac{d\mu_{1c}(t)}{t-x}\  \mbox{ and }\  m_{2c}(x) = \int_0^{\lambda_+} \frac{d\mu_{2c}(t)}{t-x}\,.
\end{equation}
\end{definition}

Our theoretical results hold with high probability in the following precise sense.

\begin{definition}[\bf Stochastic domination]\label{stoch_domination}
Let $\xi=\left(\xi^{(n)}(u):n\in\mathbb{N}, u\in U^{(n)}\right)$ and  $\zeta=\left(\zeta^{(n)}(u):n\in\mathbb{N}, u\in U^{(n)}\right)$
be two families of nonnegative random variables, where $U^{(n)}$ is a possibly $n$-dependent parameter set. We say $\xi$ is stochastically dominated by $\zeta$, uniformly in $u$, if for any fixed (small) $\epsilon>0$ and (large) $D>0$, 
$\sup_{u\in U^{(n)}}\mathbb{P}\left[\xi^{(n)}(u)>n^\epsilon\zeta^{(n)}(u)\right]\le n^{-D}$
for large enough $n \ge n_0(\epsilon, D)$, and we shall use the notation $\xi\prec\zeta$ or $\xi=O_\prec(\zeta)$. Throughout this paper, stochastic domination will always be uniform in all parameters that are not explicitly fixed, such as matrix indices, and $z$ that takes values in some compact set. Note that $n_0(\epsilon, D)$ may depend on quantities that are explicitly constant, such as $\tau$ in Assumption \ref{assum_main}. 
Moreover, we say that an event $\Xi$ holds with high probability if for any constant $D>0$, $\mathbb P(\Xi)\ge 1- n^{-D}$, when $n$ is sufficiently large.
\end{definition}

With these definitions in hand, we now state the assumptions on the model \eqref{eq_model}. Assumptions (i)--(v) govern the random matrix structure and ensure that the singular value shrinkage step is well-behaved, while Assumptions (vi)--(vii) impose regularity conditions on the signal matrix $S$ that enable the subsequent wavelet shrinkage step.

\begin{assumption}\label{assum_main}
Fix a small constant $0<\tau<1$. 
We need the following assumptions for the model \eqref{eq_model}:
\begin{itemize}
%\item[(i)] (Assumption on $x_{ij}$). Suppose $\mathcal{X}=[\varepsilon_{ij}]_{\tiny\substack{i=1,\ldots,p\\ j=1,\ldots,n}}\in \mathbb{R}^{p\times n}$ has independent entries and satisfies
%\begin{align}
%\max_{i,j} |\mathbb{E} \varepsilon_{ij}| \leq n^{-2-\tau}  \,, \  \max_{i,j}|\mathbb{E}  \varepsilon_{ij}^2 -n^{-1}| \leq  n^{-2-\tau}  \label{entry_assm1}
%\end{align}
%and
%\begin{equation}
%\max_{i,j}\mathbb{E}|\sqrt{n} \varepsilon_{ij}|^a \leq C
%\end{equation}
%for a constant $C>0$ and $a>4$.
\item[(i)] (Assumption on $\mathcal X$). Suppose the entries of \( \mathcal X \) are i.i.d. Gaussian random variables with mean zero and variance \( 1/n \).
\item[(ii)] (Assumptions on $p/n$). Define $\beta_n = p/n$. Assume $p=p(n)$ such that
\begin{equation}
\beta_n \to \beta \quad \text{for some } \beta \in (\tau, \tau^{-1}).
\end{equation}

\item[(iii)] (Assumption on $\mathcal A$ and $\mathcal B$). Assume $\mathcal A$ and $\mathcal B$ are deterministic symmetric positive semi-definite matrices with eigendecompositions
\begin{equation}
    \mathcal A = Q^a\Sigma^a(Q^a)^\top, \quad  \mathcal B = Q^b\Sigma^b(Q^b)^\top
\end{equation}
, where
$\Sigma^a=\textup{diag}(\sigma_1^a, \ldots, \sigma_p^a)$, and   $\Sigma^b=\textup{diag}( \sigma_1^b, \ldots,  \sigma_n^b)$. Without loss of generality, we let
$\sigma_1^a \geq \sigma_2^a \geq \ldots \geq \sigma_p^a \geq 0$,  $\sigma_1^b \geq \sigma_2^b \geq \ldots \geq \sigma_n^b \geq 0$, and $Q^a$, $Q^b$ are $p \times p$ and $n \times n$ orthogonal matrices, respectively.
Let $(\mathsf{M}_{1c}(z),\mathsf{M}_{2c}(z)) \in \mathbb{C}^+ \times \mathbb{C}^{+}$ as the unique solution to the following system of self-consistent equations
\begin{equation}
    \mathsf{M}_{1c}(z) = \beta_n 
    \int\frac{x}{-z[1+x\mathsf{M}_{2c}(z)]}\pi_{\mathcal A}(dx), \quad \mathsf{M}_{2c}(z) =  
    \int\frac{x}{-z[1+x\mathsf{M}_{1c}(z)]}\pi_{\mathcal B}(dx).
\end{equation}
We assume that for all sufficiently large $n$, 
\begin{equation}\label{ass3_eq1}
1+\mathsf{M}_{1c}(\lambda_+) \sigma_1^b \geq \tau\ \mbox{ and } \ 1+\mathsf{M}_{2c}(\lambda_+)  \sigma_1^a \geq \tau
\end{equation}
and 
\begin{equation}\label{ass3_eq2}
\sigma_1^a\vee  \sigma_1^b \le \tau^{-1}\ \mbox{ and } \ \pi_{\mathcal A}([0,\tau])\vee  \pi_{\mathcal B}([0,\tau]) \le 1 - \tau .
\end{equation}

\item[(iv)] (Assumption on the signal strength). %Define
%\begin{equation}
   % \phi_n:= n^{2/a-1/2}.  
%\end{equation} 
We assume \begin{equation}
      d_1 > d_2 > \cdots > d_r > 0
  \end{equation}
for some $r> 1$, $d_1<\tau^{-1}$, and the spectral gap $\min_{ij}|d_i-d_j| > \tau$. Denote a fixed value $\gamma>0$ as 
\begin{equation}\label{eq_defnalpha}
\gamma:=1/\sqrt{\mathcal{T}(\lambda_+)}\,.
\end{equation} 
We allow the singular values $d_k$ to depend on $n$ under the condition that there exists an integer $1< r^+\leq r$, called the effective rank of $\wt S$, such that
\begin{equation}\label{eq_assm_41}
d_k - \gamma > n^{\varepsilon-1/3}\ \mbox{ if and only if }\ 1\leq k\leq r^+\,,
\end{equation}
for some $\varepsilon>1/6$. 
\item[(v)] (Assumption on the distribution of singular vectors). 
Let $G_\ub^p \in \mathbb{R}^{p \times p}$ and $G_\vb^n \in \mathbb{R}^{n \times n}$ be two independent matrices with i.i.d entries distributed according to a fixed probability measure $\nu$ on $\mathbb{R}$ with mean zero and variance one, and satisfy the log-Sobolev inequality. For $i = 1,\ldots,p$ and $j = 1,\ldots,n $, we assume that the left and right singular vectors, $\ub_i \in \mathbb{R}^p$ and $\vb_j \in \mathbb{R}^n$, are the $i$-th column and the $j$-th column obtained from the Gram-Schmidt (or QR factorization) of $G_\ub^p$ and $G_\vb^n$ respectively.
\item[(vi)](Tree structure and Mixed Hölder condition on $S$). Denote its row and column spaces by \( X = \{x_1, \ldots, x_p\} \) and \( Y = \{y_1, \ldots, y_n\} \), where \( x_i \in \mathbb{R}^n \) is the \( i \)-th row of \( S \) and \( y_j \in \mathbb{R}^p \) is the \( j \)-th column of \( S \). $X$ and $Y$ are equipped with normalized counting measure, such that the measure for any single point $x\in X$ is $p^{-1}$ and the measure of any single point $y \in Y$ is $n^{-1}$. We assume that after the random generation of singular values and singular vectors described in (iv) and (v), the sets \( X \) and \( Y \) exhibit intrinsic hierarchical structures that can be represented by underlying trees \( \mathcal{T}_X \) and \( \mathcal{T}_Y \), respectively, such that $ S $ has a mixed Hölder($\alpha$) norm $ L = L(S, \alpha)$ for an $\alpha>0$ as defined in \eqref{eq_mixH}, with respect to the equipped trees \(\mathcal{T}_X\) and \(\mathcal{T}_Y\) and tree distances defined in~\eqref{eq_tree_d}.
\item[(vii)](Balanced partition tree)
We assume that the both partition trees $\mathcal T_X$ and $\mathcal T_Y$ from (vi) are balanced, meaning there are constants $ B_L $ and $ B_U $ such that \eqref{eq_balance} is satisfied for folders in $\mathcal T_X$ and $\mathcal T_Y$, respectively.
\end{itemize}
\end{assumption}

Assumptions (i)--(v) follow the setup in~\cite{su2025data} for deriving the theoretical guarantees of eOptShrink.

For Assumption~(i), \cite{su2025data} adopts a weaker condition requiring only that the entries of $\mathcal{X}$ have finite fourth moments. In this work, we assume that the entries of $\mathcal{X}$ are i.i.d.\ Gaussian to simplify the theoretical discussion. This assumption eliminates moment-dependent terms in Assumption~(iv) and in the subsequently introduced Theorem~\ref{thm_shrink} for the convergence guarantee of eOptShrink. Nevertheless, in our numerical experiments, we demonstrate that the proposed approach remains valid even when $\mathcal{X}$ is non-Gaussian, provided it satisfies weaker-moment conditions.

For Assumption~(iii), the self-consistent equations for $(\mathsf{M}_{1c}(z), \mathsf{M}_{2c}(z))$ arise naturally in random matrix theory as the fixed-point equations characterizing the limiting spectral distribution of $ZZ^\top$ and $Z^\top Z$ under the separable covariance model. The conditions in~\eqref{ass3_eq1} guarantee a regular square-root behavior of the spectral densities $\rho_{1c}$ and $\rho_{2c}$ near the bulk edge $\lambda_+$, rule out the existence of outlier eigenvalues of $ZZ^\top$, and ensure that $\mathcal{T}(\lambda_+)$ is well-defined. The conditions in~\eqref{ass3_eq2} ensure that the spectra of $\mathcal{A}$ and $\mathcal{B}$ cannot concentrate at zero.

For Assumption~(iv), compared to the original setup in~\cite{su2025data}, the condition is simplified to ensure the spectral gap $\min_{i,j}|d_i - d_j| > \tau$. A singular value $d_k$ fulfilling~\eqref{eq_assm_41} is considered a ``sufficiently strong'' signal: when a signal is sufficiently strong, it leads to a singular value outside the spectral bulk of $Z$, making it possible to study. Singular values that do not satisfy~\eqref{eq_assm_41} are considered ``weak'' and will stick to the bulk edge of $Z$, a phenomenon related to the BBP (Baik--Ben Arous--P\'{e}ch\'{e}) phase transition. This is why $r^+$ is called the effective rank of $S$.

For Assumption~(v), the randomness setup is the same as that in~\cite{su2025data} and others. The log-Sobolev inequality (see~\cite{anderson2010introduction}, Section~2.3.2 for the definition) implies that the entries of $\ub_i$ and $\vb_i$ have sub-Gaussian tails, which leads to the desired concentration property required in~\cite{su2025data}. These conditions facilitate variance calculation in the subsequent analysis of first-order singular vector perturbations.

For additional background and details of random matrix theory, we refer the reader to~\cite{benaych2012singular, nadakuditi2014optshrink, DY2019, su2025data}.

Assumptions~(vi) and (vii) encode the mixed H\"{o}lder regularity condition on the matrix $S$ with respect to the tree distances defined in~\eqref{eq_tree_d}, and the requirement that the partition trees $\mathcal{T}_X$ and $\mathcal{T}_Y$ are balanced in the sense of~\eqref{eq_balance}, respectively. These form the basis for the theoretical guarantees of $\mathcal{WS}$ established in~\cite{ankenman2018mixed}.

\subsection{Results of Optimal Shrinkage of Singular Values}\label{sec_eoptshrink}
Recall the definition of asymptotic loss and optimal shrinker provided in \cite[Definitions 1 and 2]{GD} . 
\begin{definition}[\bf Asymptotic loss]\label{def_1}
Let $\mathcal{L}:=\{L_{p,n}| p, n \in \mathbb{N}\}$ be a family of loss functions, where each $L_{p,n}: M_{p \times n} \times M_{p \times n} \rightarrow [0,\infty)$ is a loss function obeying  
that $\widehat S \to  L_{p,n}(S, \widehat S)$ is continuous and
$L_{p,n}(S,S)=0.$ Suppose $p=p(n)$ and $\lim_{n \to \infty} p(n)/n \to \beta>0$.  Let $\varphi: [0, \infty) \rightarrow [0, \infty)$ be a nonlinear function and consider $\widehat{S}_{\varphi}$ to be the singular shinkage estimate \eqref{eq_shrink}. When $\lim_{n \rightarrow \infty} L_{p,n}$ exists, we define the asymptotic loss of the shrinker $\varphi$ with respect to $L_{p, n}$ with the signal $\db=(d_1, \cdots, d_r)$ as 
$L_{\infty}(\varphi |\db) = \lim_{n \rightarrow \infty} L_{p,n}(S, \widehat{S}_{\varphi})$, where $S$ is defined in \eqref{eq_model}.
\end{definition}

\begin{definition}[\textbf{Optimal shrinker}]\label{def_2}  
Let $L_{\infty}$ and $\varphi$ be as defined in Definition \ref{def_1}. If a shrinker $\varphi^*$ has an asymptotic loss that satisfies 
$L_{\infty} (\varphi^*| \db) \leq L_{\infty} (\varphi | \db)$
for any other shrinker $\varphi$, any $r \geq 1$, and any $\db \in \mathbb{R}^r,$ then we say that $\varphi^*$ is uniquely asymptotically admissible (or simply ``optimal'') for the loss family $\mathcal{L}$ and that class of shrinkers. 
\end{definition} 

From now on, we denote the optimal shrinker of $\wt\lambda_i$ as 
\begin{equation}
\varphi^*_i :=\varphi^*(\wt\lambda_i).
\end{equation}
In \cite{GD}, Sections IV.~A and C, the optimal shrinkers under different loss functions were computed. When the loss function is the operator norm, the optimal shrinker was proved to be $\varphi_i^* = d_i$ when $p = n$ in \cite{GD}, Section IV.B. Later in \cite{leeb2020optimal}, Lemma 5.1, it is shown that $\varphi_i^* = d_i\sqrt{\frac{a_{1,i} \wedge a_{2,i}}{a_{1,i} \vee a_{2,i}}}$ when $p \neq n$, where $a_{1,i} := \lim_{n\to\infty} \langle \ub_i, \wt\ub_i \rangle^2$ and $a_{2,i} := \lim_{n\to\infty} \langle \vb_i, \wt\vb_i \rangle^2$. We list the results here for readers' convenience.
\begin{proposition}[\cite{GD,leeb2020optimal}]\label{prop_optimal_shrinker}
When $d_i \geq \alpha$, the optimal shrinker is 
$\varphi^*_i=d_i \sqrt{a_{1,i} a_{2,i}}$, $\varphi^*_i=d_i\sqrt{\frac{a_{1,i}\wedge a_{2,i}}{a_{1,i}\vee a_{2,i}}}$ and $\varphi^*_i=d_i\big(\sqrt{a_{1,i} a_{2,i}}-\sqrt{(1-a_{1,i})(1-a_{2,i})}\ \big)$ when the
Frobenius norm, operator norm, and nuclear norm are considered in the loss function, respectively. When $d_i < \gamma$, for any loss function, we have 
$\varphi^*_i=0$.
\end{proposition}

Based on Proposition \ref{prop_optimal_shrinker}, in \cite{su2025data} we estimated $d_i$, $a_{1,i}$ and $a_{2,i}$ using the eigenstructure of the noisy matrix $\wt S\wt S^\top$, and obtained estimated optimal shrinkers based on theoretical results on biased singular values and vectors, including their limiting behavior and associated convergence rates. We refer the reader to check Theorems 3.2 and 3.3 of \cite{su2025data} for the convergence of singular values $\wt \lambda_i \to \mathcal T^{-1}(d_i^{-2})$ and Theorems 3.4 and 3.5 of \cite{su2025data} for the convergence of inner product
$a_{1,i}\to a^{\infty}_{1,i}$ and $a_{2,i}\to a^{\infty}_{2,i}$, which form the foundation of the following eOptShrink algorithm.\\

\subsubsection{eOptShrink}
Define 
\begin{equation}
\mathbb O_+ := \{1,\cdots,r^+\}
\ \mbox{ and } \  
    \Delta(d_i) := |d_i-\gamma|^{1/2}\,.
\end{equation}
Define an estimator of the bulk edge
\begin{equation}\label{eq_hatlanbda_+}
    \widehat{\lambda}_{+}  := \wt\lambda_{\lfloor n^c\rfloor +1} + \frac{1}{2^{2/3}-1}\left( \wt\lambda_{\lfloor n^c\rfloor +1}-\wt\lambda_{2\lfloor n^c\rfloor +1} \right) \,.
\end{equation}
Then we estimate $r^+$ and set
\begin{equation}\label{eq_adrk0}
\widehat r^+ := \left|\{\wt\lambda_i|\wt \lambda_i>\widehat \lambda_{+}+n^{-1/3}\}\right|\,.
\end{equation}
The following theorems from \cite{su2025data} guarantee the performance of $\widehat r^+$. 
\begin{theorem}[Theorem 4.2 \cite{su2025data}]\label{thm_numbershrink2}
Suppose (i)-(v) of Assumption \ref{assum_main} hold true. Denote the event $\Xi(r^+):=\{\widehat{r}^+=r^+\}$ and let $0<c<1/2$. Then $\Xi(r^+)$ holds with high probability in the sense of Definition~\ref{stoch_domination}. 
\end{theorem}

Denote
\begin{equation}\label{eq_hatlanbda2}
\widehat{\lambda}_{j}  := \wt\lambda_{\lfloor n^c\rfloor+ \widehat r^++1} + \frac{1-\big(\frac{j-\widehat r^+-1}{\lfloor n^c\rfloor}\big)^{2/3}}{2^{2/3}-1}\left( \wt\lambda_{\lfloor n^c\rfloor+ \widehat r^++1}-\wt\lambda_{2\lfloor n^c\rfloor+\widehat r^++1} \right), 
\end{equation}
The estimated CDF of $\pi_{ZZ^\top}$ in \cite{su2025data} is
\begin{equation}\label{eq_imp}
    \widehat F_{\texttt{e}}(x) := \frac{1}{p-\widehat r^+} \left( \sum_{j=\widehat r^++1}^{\lfloor n^c\rfloor+\widehat r^+} \mathbbm 1(\widehat\lambda_{j} \le x)+\sum_{j=\lfloor n^c\rfloor+\widehat r^++1}^p \mathbbm 1(\wt\lambda_{j}\leq x)\right) \,.
\end{equation}

For $1\leq i \leq \widehat r^+$, 
denote the estimators of $m_{1c}(\wt\lambda_i)$ and $m_{2c}(\wt\lambda_i)$ as
\begin{align}\label{eq_mhat1mhat2}
& \widehat{m}_{e,1,i}:= \int \frac{d\widehat F_e(x)}{x-\wt\lambda_i} =  \frac{1}{p-\widehat r^+}\left( \sum_{j=\widehat r^++1}^{\lfloor n^c \rfloor+\widehat r^+} \frac{1}{\widehat\lambda_j-\wt\lambda_i} + \sum_{j=\lfloor n^c\rfloor+\widehat r^++1}^p \frac{1}{\wt\lambda_j-\wt\lambda_i}\right), \ 
\nonumber\\
& \widehat{m}_{e,2,i} := \frac{1-\beta_n}{\wt\lambda_i}+\beta_n \widehat{m}_{e,1,i}. 
\end{align}
For $1\leq i \leq \widehat r^+$, the estimator of the $D$-transform $\mathcal{T}(\wt\lambda_i)$ is
\begin{equation}
\widehat{\mathcal{T}}_{e,i}= \wt\lambda_i \widehat{m}_{e,1,i} \widehat{m}_{e,2,i},
\end{equation}
and the estimators of $d_i$ and inner products of the clean and noisy left singular vectors ${a}_{1,i}$ and right singular vectors ${a}_{2,i}$, are given by 
 \begin{equation}\label{eq_a1a2}
  \widehat{d}_{e,i}=\sqrt{\frac{1}{\widehat{\mathcal{T}}_{e,i}}}, \quad \widehat{a}_{e,1,i}=\frac{\widehat{m}_{e,1,i}}{\widehat{d}_{e,i}^2 \widehat{\mathcal{T}}'_{e,i}}, \ \mbox{ and } \  \widehat{a}_{e,2,i}=\frac{\widehat{m}_{e,2,i}}{ \widehat{d}_{e,i}^2\widehat{\mathcal{T}}'_{e,i}}.
 \end{equation}
As a result, we estimate the optimal shrinker $\varphi^*_i$ in Proposition \ref{prop_optimal_shrinker} by
\begin{equation}\label{eq_adsh}
\begin{array}{lr}
\displaystyle\widehat{\varphi}_{\texttt{e},i}=\widehat{d}_{\texttt{e},i} \sqrt{\widehat{a}_{\texttt{e},1,i}\widehat{a}_{\texttt{e},2,i}}, &\mbox{(Frobenius norm)}  \\
\displaystyle{\widehat{\varphi}_{\texttt{e},i}=\widehat d_{e,i}\sqrt{\frac{\widehat a_{e,1,i}\wedge \widehat a_{e,2,i}}{\widehat a_{e,1,i}\vee \widehat a_{e,2,i}}}}, &\mbox{(Operator norm)}\\
\displaystyle\widehat{\varphi}_{\texttt{e},i}=\widehat{d}_{\texttt{e},i}\Big(\sqrt{\widehat{a}_{\texttt{e},1,i} \widehat{a}_{\texttt{e},2,i}}-\sqrt{(1-\widehat{a}_{\texttt{e},1,i})(1-\widehat{a}_{\texttt{e},2,i})}\Big) \,, &\mbox{(Nuclear norm)}
\end{array}
\end{equation} 
for $1\leq i\leq \widehat r^+$, and $\widehat{\varphi}_{\texttt{e},i}=0$ otherwise. The following theorem provides the convergence guarantee:
\begin{theorem}[Theorem 4.4 of \cite{su2025data}]\label{thm_shrink} 
Suppose (i)-(v) of Assumption \ref{assum_main} hold true for some $\varepsilon>1/6$, and $c\in (0,1/2)$. 
For all three types of loss functions mentioned in Proposition \ref{prop_optimal_shrinker}, for $1\leq i \leq \widehat r^+$, conditional on $\Xi(r^+)$, we have
$|\varphi^*_i-\widehat{\varphi}_{\texttt{e},i}| \prec n^{-1/2}/\Delta(d_i)$.
\end{theorem}

%\subsection{Results from Mixed Hölder Condition}
%For a Haar function $ \phi_X(x) $ on $ X $, let $ I(\phi_X) $ denote the smallest folder containing the support of $ \phi_X $. Similarly, let $ J(\phi_Y) $ denote the smallest folder containing the support of a Haar function $ \phi_Y $ on $ Y $. For a tensor Haar function $\Phi_XY(x,y) = \phi_X(x)\phi_Y(y)$ on $ X \times Y $, we denote by $ R(\Phi_{IJ}) = I(\phi_{X}) \times J(\phi_Y) $ the smallest rectangle containing the support of $ \Phi_{IJ} $.
%The following result from \cite{ankenman2014geometry} characterizes the mixed Hölder($\alpha$) norm of a function $ S $ based on its tensor Haar coefficients $\langle S, \Phi \rangle $:

%\begin{theorem}\label{thm_holder_coeff}
%Suppose (ii) and (vi) of Assumption \ref{assum_main} hold.
%For any $\alpha > 0$, there is a constant $ C = C(B_L, B_U, \alpha) > 1 $ such that for any function $ S $ on $ X \times Y $,
%\begin{equation}
%\frac{1}{C} L(S, \alpha) \leq \sup_{\Phi_{IJ}} \frac{|\langle S, \Phi_{IJ} \rangle|}{|R(\Phi_{IJ})|^{\alpha+1/2}} \leq L(S, \alpha),
%\end{equation}
%where the supremum is over all Haar functions $\Phi_{IJ}(x,y)$.
%\end{theorem}
\subsection{Results of Wavelet Shrinkage}\label{sec_ws}
We introduce the theoretical results of wavelet shrinkage in 
this section. Let $S$ be a mixed Hölder($\alpha$) matrix with 
norm $L$ as defined in Definition~\ref{def_mixed_Holder}, and 
be contaminated by additive noise $Z$ as in 
model~\eqref{eq_model}. When $\alpha$ and $L$ are not 
specified, in~\cite{donoho1995adapting, donoho1996neo, 
ankenman2018mixed}, the \textit{Wavelet Shrinkage} scheme is 
applied to construct an estimator from $\wt S$ over all 
$\alpha > 0$. We note that 
in~\cite{donoho1994ideal, donoho1995adapting, 
donoho1996neo, ankenman2018mixed}, the results are stated 
under $\mathcal{O}(1)$ entry-wise noise variance, whereas our 
framework assumes $\mathcal{O}(n^{-1})$ variance as in 
Assumption~\ref{assum_main}\,(i). The results below have been 
restated to match the scaling and conventions of our framework.

Expand $S$ in a two-dimensional Haar series ${\Phi_{IJ}}$ as defined in \eqref{eq_ws_expand}.
Since the Haar transform is orthogonal, each observed Haar coefficient $\langle \wt S, \Phi_{IJ} \rangle$ is normally distributed around the true coefficient $\langle S, \Phi_{IJ} \rangle$, with variance $\sigma^2 / n$.
To have an estimator $\mathcal{WS}$ of $\langle S, \Phi_{IJ} \rangle$ from $\langle \wt S, \Phi_{IJ} \rangle$,
Donoho and Johnstone considered the loss function \begin{equation} \frac{\mathbb{E}_{\theta}[(\widehat{\theta} - \theta)^{2}]}{\delta + m_{\delta}^{\ell}(\theta)}\,,
\end{equation}
where $\widehat{\theta}$ is the estimator of $\theta$, $\delta \in (0, 1)$, $\ell \in (0, 2)$, and  $m_{\delta}^{\ell}(\theta)$ is defined by
\begin{equation}
m_{\delta}^{\ell}(\theta) := t^{\ell} \min \left\{ \left(\frac{\theta}{t}\right)^{2}, 1 \right\}\,, 
\end{equation}
with $t = \sqrt{2\log(1/\delta)}$. 
The loss function penalizes errors made at small values of $\theta$ more than those made at large values of $\theta$. The wavelet shrinker defined in \eqref{eq_ws_shrinker} is consistent with the structure of the loss function, as it moves \( x \) toward zero by an amount \( \min\{t, |x|\} \), thereby applying stronger shrinkage to smaller values of \( x \).
The following corollary from \cite{ankenman2018mixed} shows that the shrinker yields an upper-bound of the loss function:

\begin{corollary} [Corollary 1 \cite{ankenman2018mixed}]\label{col_wsh}
 Suppose $\wt \theta \sim N(\theta, \sigma^2)$, where $\sigma^2$ is assumed to be known. $\ell \in (0,2)$. Define the estimator of $\theta$ to be $\widehat{\theta} = \eta_{\sigma t}(Y)$, where $t = \sqrt{2 \log(1/\delta)}$ and $\delta \in (0, 1)$. Then there is a constant $C > 0$ such that for all $\delta$ sufficiently small,
\begin{equation}
\sup_{\theta} \frac{\mathbb{E}_{\theta}[(\widehat{\theta} - \theta)^{2}]}{\delta + m_{\delta}^{\ell}(\theta / \sigma)} \leq C \sigma^2 \log(1/\delta)^{1 - \ell/2}.
\end{equation}
\end{corollary}
Construct the estimator of wavelet coefficient $\langle \wt S, \Phi_{IJ} \rangle$ by the shrinker $\eta_{\frac{t \sigma}{\sqrt{n}}}$ such that  \begin{equation}\label{eq_ws_shrinker2}
\widehat{C}_{\Phi_{IJ}} := \eta_{\frac{t\sigma}{\sqrt{n}}}(\langle \wt S, \Phi_{IJ} \rangle)
\end{equation}
for nonconstant tensor Haar functions $\Phi_{IJ}$, and $\widehat{C}_{1} = \langle \wt S, 1 \rangle$ and construct the estimator using the shrinker as
\begin{equation}
\widehat S_{t}(x, y) = \sum_{\Phi_{IJ}} \widehat{C}_{\Phi_{IJ}} \Phi_{IJ}(x, y). 
\end{equation}
A common threshold choice is $t = t^* = \sqrt{2\log(pn)}$. Since the wavelet coefficients of the noise are i.i.d. Gaussian, and the maximum of $N = pn$ i.i.d. standard Gaussian random variables satisfies $\mathbb{E}\max_{i \leq N} |g_i| = \sqrt{2\log N}(1 + o(1))$ as $N \to \infty$ (see, e.g., \cite{vershynin2026high}, Exercise 2.38), the constant $\sqrt{2}$ is asymptotically sharp. This threshold ensures that, with high probability, all noise-dominated coefficients are suppressed while signal-dominated components are retained. When \( c = \frac{2\alpha}{2\alpha + 1} \) and \( \ell = 2(1 - c) \), the following theorem from~\cite{ankenman2018mixed} follows directly from Corollary~\ref{col_wsh} and provides a theoretical guarantee for the performance of \( \widehat{S}_{t^*} \).
\begin{theorem}\label{thm_ws_bound}
Suppose (vi)-(vii) of Assumption \ref{assum_main} hold, $S$ has mixed H\"older($\alpha$) condition with constant $L$ as defined in \eqref{eq_mixH} for some $\alpha>0$, and entries of $Z$ are i.i.d. normally distributed with mean $0$ and variance $\sigma^2/n$. With $t^* = \sqrt{2\log(pn)}$, we have:\\ 
(1)(Section 4.1 \cite{ankenman2018mixed}. Mean squared error of $\widehat S_{t^*}$.) 
\begin{align}\label{eq_ws_f}
    &\mathbb{E}_S\norm{S - \widehat S_{t^*}}_{2,n}^2 \nonumber \\
    &\leq C\log^{2\alpha/(2\alpha+1)}(pn)(\sigma^2L^{1/\alpha}/(pn^2))^{2\alpha/(2\alpha+1)}\log_{B_U^{-1}}(pn^2L/\sigma^2)(1+o(1)),
\end{align}
with $C = C(B_L,B_U,\alpha)$.\\
(2)(Section 4.2 \cite{ankenman2018mixed}. Pointwise squared error of $\widehat S_{t^*}$.)\\
\noindent
For any $(x_0,y_0)\in X \times Y$, we have
\begin{align}\label{eq_ws_pt}
     &\mathbb{E}(S(x_0,y_0) - \widehat S_{t^*}(x_0,y_0))^2 \nonumber\\
     & \leq C\log^{2\alpha/(2\alpha+1)}(pn)(\sigma^2L^{1/\alpha}/(pn^2))^{2\alpha/(2\alpha+1)}\log^2_{B_U^{-1}}(pn^2L/\sigma^2)(1+o(1)),
\end{align}
with $C = C(B_L,B_U,\alpha)$.  
\end{theorem}

\begin{remark}
(\textit{a}) Note that Corollary \ref{col_wsh} does not depend on the choice of orthogonal basis, therefore the same asymptotic result as in Theorem~\ref{thm_ws_bound} holds when the basis is replaced by the optimal tensor elements $\Phi_{IJ} = \omega_I \omega_J$ selected via the eGHWT best basis algorithm. As demonstrated in \cite{saito2022eghwt}, the eGHWT best basis yields sparser representations and greater energy concentration compared to the standard Haar basis, which motivates our adoption of the eGHWT in this work.
(\textit{b}) Both right-hand sides of \eqref{eq_ws_f} and \eqref{eq_ws_pt} converge to zero as \( n \to \infty \). The shrinker defined in \eqref{eq_ws_shrinker} does not require the parameters \( L \) or \( \alpha \), and depends only on the noise level \( \sigma \), making it particularly suitable for data-driven applications where such parameters may be difficult to estimate.
However, the entries of the noise matrix \( Z \) are typically not i.i.d.\ Gaussian, for instance, the correlated and dependent noise model in \eqref{colored and dependent noise model}. In such settings, a single global parameter \( \sigma \) is insufficient to adapt to the varying noise levels across different scales in the tree structure of \( S \).
Moreover, the tree structures \( \mathcal{T}_X \) and \( \mathcal{T}_Y \) are often constructed directly from the raw data matrix. When the noise \( Z \) is sufficiently strong, it can distort these trees and, as a result, degrade the quality of the associated wavelet basis. These two limitations motivate our combined approach: we first apply eOptShrink to correct the biased singular values and vectors, producing a denoised estimate with reduced noise. The partition trees $\mathcal{T}_X$ and $\mathcal{T}_Y$ are then constructed from this denoised estimate rather than the raw data, yielding more reliable tree structures. Finally, we apply wavelet shrinkage via the eGHWT best basis with scale-adaptive thresholds that account for the residual colored noise structure, rather than relying on a single global $\sigma$.
\end{remark}

\section{Proposed Algorithm}\label{sec_alg}

Our proposed extended Optimal and Wavelet Shrinkage (e$\mathcal{OWS}$) methodology operates in a sequential, plug-in manner. We first construct the estimators $\widehat{r}^+$ and $\widehat{d}_{e,i}$, $\widehat{a}_{e,1,i}$, and $\widehat{a}_{e,2,i}$ for $i = 1,\dots, \widehat{r}^+$, based on \eqref{eq_a1a2}. Let $\widehat{S}_{\textup{eOptShrink}} = \sum_{i=1}^{\widehat{r}^+}\widehat{\varphi}_{e,i} \wt\ub_i\wt\vb_i^\T$ serve as the initial estimator by eOptShrink, where the shrinker $\widehat{\varphi}_{e,i}$ is defined in \eqref{eq_adsh}. We subsequently extract the empirical noise matrix as $\widehat{Z} := \wt{S} - \widehat{S}_{\textup{eOptShrink}}$. 

To adaptively shrink the wavelet coefficients, we must estimate their localized variances. For each estimated rank $i \in \{1, \dots, \widehat{r}^+\}$, we define the empirical resolvent sums:
\begin{align}
    \widehat{\mathcal{S}}_i &:= \frac{1}{np} \sum_{k=\widehat{r}^+ + 1}^p \frac{\widehat{d}_{e,i}^2 + \widetilde\lambda_k}{(\widetilde\lambda_i - \widetilde\lambda_k)^2}\,,\\[4pt]
    \widehat{\mathcal{R}}_i &:= \frac{\widehat{d}_{e,i}^2}{n}\left(\frac{1}{p}\sum_{k=\widehat{r}^+ + 1}^p \frac{1}{\widetilde\lambda_i - \widetilde\lambda_k}\right)^{\!2}\,.
\end{align}
Using these sums and the extracted noise matrix, we define the per-component variance estimates:
\begin{align}
    \widehat{\mathcal{V}}_{1,i} &:= \frac{1}{n}\bigg[\widehat{\mathcal{R}}_i\cdot\frac{\|\omega_I^\top\widehat{Z}\|_2^2}{n} + \widehat{\mathcal{S}}_i\cdot\frac{\|\widehat{Z}\|_F^2}{pn}\bigg]\,,\\[4pt]
    \widehat{\mathcal{V}}_{2,i} &:= \frac{1}{p}\bigg[\widehat{\mathcal{R}}_i'\cdot\frac{\|\widehat{Z}\omega_J\|_2^2}{p} + \widehat{\mathcal{S}}_i\cdot\frac{\|\widehat{Z}\|_F^2}{pn}\bigg]\,,
\end{align}
where $\widehat{\mathcal{R}}_i' := \widehat{d}_{e,i}^2(\frac{1}{n}\sum_{k=\widehat{r}^++1}^n \frac{1}{\widetilde\lambda_i - \widetilde\lambda_k'})^2$ is the right-side analog using the eigenvalues of $\widetilde{S}^\top\widetilde{S}$, and $\norm{\cdot}_F$ is the Frobenius norm. The data-driven variance estimator for the wavelet coefficient associated with the eGHWT basis $\Phi_{IJ} = \omega_I\omega_J^\top$ is then:
\begin{equation}\label{eq_sigma+phi2_final}
    \widehat\sigma_{\Phi_{IJ}}^2 := \sum_{i=1}^{\widehat{r}^+} \widehat{d}_{e,i}^2 \Big(\widehat{\mathcal{V}}_{1,i} + \widehat{\mathcal{V}}_{2,i} + np\,\widehat{\mathcal{V}}_{1,i}\,\widehat{\mathcal{V}}_{2,i}\Big)\,.
\end{equation}
Based on this estimated variance, we modify the standard wavelet shrinker from \eqref{eq_ws_shrinker} to adapt dynamically to each coefficient in the eGHWT basis, taking the form $\eta_{t\widehat{\sigma}_{\Phi_{IJ}}}(\langle \widehat{S}_{\textup{eOptShrink}}, \Phi_{IJ} \rangle)$. To determine the optimal threshold scaling factor $t$, we apply the universal thresholding principle discussed in Section~\ref{sec_ws}. As shown in Theorem~\ref{thm_variance}, the variation in the estimated coefficients is fundamentally driven by the noise matrix $Z$. Because $Z$ is assumed to be sub-Gaussian, the resulting standardized noise coefficients are also sub-Gaussian. Thus, their maximum absolute value across the $pn$ entries is tightly bounded by $\sqrt{2\log(pn)}$ with high probability. We therefore set the universal scaling factor to $\tau^* = \sqrt{2\log(pn)}$, yielding the localized wavelet shrinker $\eta_{\Phi_{IJ}}^* := \eta_{\tau^*\widehat\sigma_{\Phi_{IJ}}}$. The complete procedure for generating the final estimator $\widehat{S}_{e\mathcal{OWS}}$ is summarized in Algorithm~\ref{alg_eqOS}.

\begin{algorithm}[hbt!]
\begin{algorithmic}
\caption{e$\mathcal{OWS}$}\label{alg_eqOS}
\STATE \textbf{Input}: $\widetilde{S} = \sum_{i=1}^{p \wedge n} \sqrt{\widetilde{\lambda}_i} \, \widetilde{\ub}_i \widetilde{\vb}_i^\top$; a constant $c = \min\left(\frac{1}{2.01}, \frac{1}{\log(\log n)}\right)$; the choice of loss function (Frobenius, operator, or nuclear norm).
\vspace{2mm}
\STATE \textbf{Compute}: 
\begin{itemize}
\item[(i)] The estimator $\widehat{r}^+$ for the effective rank via \eqref{eq_adrk0}.
\item[(ii)] The estimators $\widehat{d}_{e,i}$ via \eqref{eq_a1a2} and the shrinkers $\widehat\varphi_{e,i}$ via \eqref{eq_adsh} for $i = 1,\dots, \widehat{r}^+$.
\item[(iii)] The intermediate matrices $\widehat{S}_{\widehat r^+} = \sum_{i = 1}^{\widehat{r}^+} {\widehat d_{e,i}} \wt\ub_i\wt\vb_i^\T$, $\widehat S_{\textup{eOptShrink}} = \sum_{i = 1}^{\widehat{r}^+} {\widehat\varphi_{e,i}} \wt\ub_i\wt\vb_i^\T$, and the estimated noise matrix $\widehat{Z} = \wt{S} - \widehat S_{\textup{eOptShrink}}$.
\item[(iv)] The row and column partition trees $\mathcal{T}_X$ and $\mathcal{T}_Y$ derived from $\widehat S_{\textup{eOptShrink}}$ using the Questionnaire approach detailed in Appendix \ref{sec_question}. 
\item[(v)] The eGHWT best basis $\{\Phi_{IJ}\}_{I\in \mathcal{T}_X, J\in \mathcal{T}_Y}$ induced by $\mathcal{T}_X$ and $\mathcal{T}_Y$, and the corresponding coefficients $\widehat{C}_{\Phi_{IJ}} = \langle \widehat{S}_{\widehat r^+}, \Phi_{IJ} \rangle$.
\item[(vi)] The estimated variance $\widehat\sigma^2_{\Phi_{IJ}}$ for each eGHWT coefficient via \eqref{eq_sigma+phi2_final}. 
\item[(vii)] The universal threshold scaling factor $\tau^* = \sqrt{2\log(pn)}$. 

\end{itemize}
\STATE \textbf{Output}: 
  The spatially shrunk matrix:
    $$ \widehat{S}_{e\mathcal{OWS}} = \sum_{\Phi_{IJ}} \eta_{\tau^*\widehat\sigma_{\Phi_{IJ}}}(\widehat{C}_{\Phi_{IJ}})\Phi_{IJ}\,. $$
\end{algorithmic}
\end{algorithm}

\section{THEORETICAL GUARANTEES}\label{sec_result}
All proofs of the lemmas, propositions, and theorems stated in this section are provided in Appendix~\ref{sec_proofs}.

To rigorously quantify the statistics of orthogonal perturbations under the separable covariance structure of \eqref{colored and dependent noise model}, we embed the observed data matrix into a continuous-time stochastic process. We parameterize the noise evolution by introducing a fictitious time parameter $t \in [0, 1/n]$, defining the matrix-valued stochastic process:
\begin{equation}
    \widetilde{S}(t) = S + Z(t)\,,
\end{equation}
with the initial condition $\widetilde{S}(0) = S$. Here, the noise increment $dZ(t) = \mathcal{A}^{1/2} d\mathcal{X}(t) \mathcal{B}^{1/2}$ is driven by a standard matrix-valued Brownian motion $\mathcal{X}(t) \in \mathbb{R}^{p \times n}$. 

Denote the singular value decomposition of the unperturbed signal matrix as:
\begin{equation}\label{eq_svd_S}
    S = U \Sigma V^\top = U_r \Sigma_r V_r^\top + U_c \Sigma_c V_c^\top\,,
\end{equation}
where $U_r = [\mathbf{u}_1, \ldots, \mathbf{u}_r]$ and $V_r = [\mathbf{v}_1, \ldots, \mathbf{v}_r]$ span the rank-$r$ signal subspace associated with the nonzero singular values $\Sigma_r = \textup{diag}(d_1,\ldots,d_r)$. The matrices $U_c=[\mathbf{u}_{r+1}, \ldots, \mathbf{u}_p]$ and $V_c = [\mathbf{v}_{r+1}, \ldots, \mathbf{v}_n]$ span the orthogonal noise subspaces associated with $\Sigma_c = 0$.

Under the noise perturbation $Z(t)$, the subspace decomposition evolves smoothly as $\widetilde{S}(t) = \widetilde{U}(t) \widetilde{\Sigma}(t) \widetilde{V}(t)^\top$. We systematically partition the time-dependent left and right singular matrices into the rank-$r$ signal subspace and the orthogonal noise subspace:
\begin{align}
    \widetilde{U}(t) &= \big[ \widetilde{U}_r(t) \mid \widetilde{U}_c(t) \big] = \big[ \widetilde{\mathbf{u}}_1(t), \dots, \widetilde{\mathbf{u}}_r(t) \mid \widetilde{\mathbf{u}}_{r+1}(t), \dots, \widetilde{\mathbf{u}}_p(t) \big]\,, \\
    \widetilde{V}(t) &= \big[ \widetilde{V}_r(t) \mid \widetilde{V}_c(t) \big] = \big[ \widetilde{\mathbf{v}}_1(t), \dots, \widetilde{\mathbf{v}}_r(t) \mid \widetilde{\mathbf{v}}_{r+1}(t), \dots, \widetilde{\mathbf{v}}_n(t) \big]\,.
\end{align}
The time-dependent singular values are collected in $\widetilde{\Sigma}(t) = \textup{diag}(\widetilde{d}_1(t), \ldots, \widetilde{d}_{p \wedge n}(t))$. We enforce boundary conditions such that at $t=0$, the process perfectly matches the unperturbed SVD: $\widetilde{\mathbf{u}}_i(0) \equiv \mathbf{u}_i$ and $\widetilde{\mathbf{v}}_i(0) \equiv \mathbf{v}_i$ for all $i$, with $\widetilde{d}_i(0) = d_i$ for $1 \le i \le r$, and $\widetilde{d}_i(0) = 0$ for $i > r$. At the terminal time $t = 1/n$, the process reaches the fully perturbed empirical quantities: $\widetilde{\mathbf{u}}_i(1/n) \equiv \widetilde{\mathbf{u}}_i$, $\widetilde{\mathbf{v}}_i(1/n) \equiv \widetilde{\mathbf{v}}_i$, and $\widetilde{d}_i^2(1/n) = \widetilde{\lambda}_i$. 

We track the variations of the perturbed singular vectors for $i = 1, \ldots, r$. By applying It\^o's lemma and building upon the first-order algebraic perturbation framework introduced in \cite{liu2008first}, we derive the exact stochastic differential equations governing these vectors.

\begin{lemma}\label{lemma_exact_sde}
Suppose (i)-(v) of Assumption \ref{assum_main} hold. Let $\widetilde{S}(t)$ evolve as defined above, partitioned dynamically into its signal and noise subspaces. For each signal component $i \in \{1, \dots, r\}$, define the continuous-time resolvent matrices governing the spectral gaps:
\begin{align*}
    \widetilde{D}_i^r(t) &:= \textup{diag}\left( \frac{1}{\widetilde{d}_i^2(t) - \widetilde{d}_k^2(t)} \right)_{k=1, \dots, r}^{k \neq i} \in \mathbb{R}^{r \times r} \quad (\text{with } 0 \text{ at the } i\text{-th diagonal entry})\,, \\
    \widetilde{D}_i^{c,u}(t) &:= \left(\widetilde{d}_i^2(t) I_{p-r} - \widetilde{\Sigma}_c(t)\widetilde{\Sigma}_c(t)^\top\right)^{-1} \in \mathbb{R}^{(p-r) \times (p-r)}\,, \\
    \widetilde{D}_i^{c,v}(t) &:= \left(\widetilde{d}_i^2(t) I_{n-r} - \widetilde{\Sigma}_c(t)^\top\widetilde{\Sigma}_c(t)\right)^{-1} \in \mathbb{R}^{(n-r) \times (n-r)}\,.
\end{align*}

The exact stochastic differentials of the perturbed singular vectors for $i = 1,\ldots,r$ are given by:
\begin{align}
    d\widetilde{\mathbf{u}}_i(t) &= d\widetilde{\mathbf{u}}_i^r(t) + d\widetilde{\mathbf{u}}_i^c(t) + \mathcal{O}(dt)\,, \label{eq_exact_du_split} \\
    d\widetilde{\mathbf{v}}_i(t) &= d\widetilde{\mathbf{v}}_i^r(t) + d\widetilde{\mathbf{v}}_i^c(t) + \mathcal{O}(dt)\,, \label{eq_exact_dv_split}
\end{align}
where the variations constrained strictly within the rank-$r$ signal subspace are:
\begin{align}
    d\widetilde{\mathbf{u}}_i^r(t) &:= \widetilde{d}_i(t) \widetilde{U}_r(t) \widetilde{D}_i^r(t) \widetilde{U}_r(t)^\top (dZ(t)) \widetilde{\mathbf{v}}_i(t) + \widetilde{U}_r(t) \widetilde{D}_i^r(t) \widetilde{\Sigma}_r(t) \widetilde{V}_r(t)^\top (dZ(t))^\top \widetilde{\mathbf{u}}_i(t)\,, \label{eq_exact_du_r} \\
    d\widetilde{\mathbf{v}}_i^r(t) &:= \widetilde{d}_i(t) \widetilde{V}_r(t) \widetilde{D}_i^r(t) \widetilde{V}_r(t)^\top (dZ(t))^\top \widetilde{\mathbf{u}}_i(t) + \widetilde{V}_r(t) \widetilde{D}_i^r(t) \widetilde{\Sigma}_r(t) \widetilde{U}_r(t)^\top (dZ(t)) \widetilde{\mathbf{v}}_i(t)\,, \label{eq_exact_dv_r}
\end{align}
and the variations spilling orthogonally into the noise subspace are:
\begin{align}
    d\widetilde{\mathbf{u}}_i^c(t) &:= \widetilde{d}_i(t) \widetilde{U}_c(t) \widetilde{D}_i^{c,u}(t) \widetilde{U}_c(t)^\top (dZ(t)) \widetilde{\mathbf{v}}_i(t) + \widetilde{U}_c(t) \widetilde{D}_i^{c,u}(t) \widetilde{\Sigma}_c(t) \widetilde{V}_c(t)^\top (dZ(t))^\top \widetilde{\mathbf{u}}_i(t)\,, \label{eq_exact_du_c} \\
    d\widetilde{\mathbf{v}}_i^c(t) &:= \widetilde{d}_i(t) \widetilde{V}_c(t) \widetilde{D}_i^{c,v}(t) \widetilde{V}_c(t)^\top (dZ(t))^\top \widetilde{\mathbf{u}}_i(t) + \widetilde{V}_c(t) \widetilde{D}_i^{c,v}(t) \widetilde{\Sigma}_c(t)^\top \widetilde{U}_c(t)^\top (dZ(t)) \widetilde{\mathbf{v}}_i(t)\,. \label{eq_exact_dv_c}
\end{align}
\end{lemma}

This continuous-time framework allows us to rigorously isolate the first-order perturbation approximations $d\widetilde{\mathbf{u}}_i^r(t)$ and $d\widetilde{\mathbf{u}}_i^c(t)$, which are strictly linear with respect to the noise increment $dZ(t)$, from the higher-order It\^o drift $\mathcal{O}(dt)$. By integrating the exact stochastic differentials over the continuous path $t \in [0, 1/n]$, we recover the singular vectors of the final observed matrix. The total perturbations can thus be explicitly decomposed into signal-space and noise-space contributions:
\begin{equation}\label{eq_order}
    \Delta \mathbf{u}_i = \Delta \mathbf{u}_i^r + \Delta \mathbf{u}_i^c + \mathcal{O}(n^{-1})\,, \quad \Delta \mathbf{v}_i = \Delta \mathbf{v}_i^r + \Delta \mathbf{v}_i^c + \mathcal{O}(n^{-1})\,,
\end{equation}
where the cumulative variations constrained within the signal subspaces are $\Delta \mathbf{u}_i^r := \int_0^{1/n} d\widetilde{\mathbf{u}}_i^r(t)$ and $\Delta \mathbf{v}_i^r := \int_0^{1/n} d\widetilde{\mathbf{v}}_i^r(t)$, and the cumulative variations spilling into the orthogonal noise subspaces are $\Delta \mathbf{u}_i^c := \int_0^{1/n} d\widetilde{\mathbf{u}}_i^c(t)$ and $\Delta \mathbf{v}_i^c := \int_0^{1/n} d\widetilde{\mathbf{v}}_i^c(t)$. The $\mathcal{O}(n^{-1})$ residual arises directly from integrating the $\mathcal{O}(dt)$ higher-order It\^o correction terms over the length-$1/n$ integration path.

The following theorem establishes the precise leading-order variance of the singular vector perturbations and their projections onto the spatial wavelet domain. 

We define two resolvent integrals over the evolving noise spectrum. For each signal component $i \in \{1,\dots,r^+\}$, let $\widetilde{d}_k(t)$ denote the $k$-th singular value of the continuously perturbed matrix $\widetilde{S}(t) = S + Z(t)$. Define:
\begin{align}
    \mathcal{S}_i &:= \mathbb{E}\left[ \int_0^{1/n} \frac{1}{p} \sum_{j=1}^{p-r} \frac{d_i^2 + \widetilde{d}_{r+j}(t)^2}{\big(\widetilde{d}_i(t)^2 - \widetilde{d}_{r+j}(t)^2\big)^2}\,dt \right]\,,\label{eq_Si}\\[4pt]
    \mathcal{R}_i &:= \mathbb{E}\left[ \int_0^{1/n} \widetilde{d}_i(t)^2 \left(\frac{1}{p}\sum_{j=1}^{p-r} \frac{1}{\widetilde{d}_i(t)^2 - \widetilde{d}_{r+j}(t)^2}\right)^2\,dt \right]\,.\label{eq_Ri}
\end{align}

\begin{theorem}\label{thm_variance}
Suppose (i)--(v) of Assumption~\ref{assum_main} hold. Let $S_{r^+} = \sum_{i=1}^{r^+} d_i \mathbf{u}_i \mathbf{v}_i^\top$ and $\widehat{S}_{r^+} = \sum_{i=1}^{r^+} d_i \widetilde{\mathbf{u}}_i \widetilde{\mathbf{v}}_i^\top$. For $i = 1,\ldots,r^+$, the following hold.

\begin{itemize}
    \item[\textup{(i)}] \textbf{Noise subspace dominance.} The signal-subspace variation satisfies $\mathbb{E}[(\Delta\mathbf{u}_i^r)_a^2] = \mathcal{O}(n^{-2})$ for all $1 \leq a \leq p$, uniformly in $\Delta(d_i)$. The total entry-wise variance is therefore determined by the noise-subspace component:
    \begin{equation}\label{eq_gap_v}
        \mathbb{E}[(\Delta \mathbf{u}_i)_a^2] = \mathbb{E}[(\Delta \mathbf{u}_i^c)_a^2] + \mathcal{O}(n^{-2})\,, \quad
        \mathbb{E}[(\Delta \mathbf{v}_i)_b^2] = \mathbb{E}[(\Delta \mathbf{v}_i^c)_b^2] + \mathcal{O}(n^{-2})\,.
    \end{equation}
 
    \item[\textup{(ii)}] \textbf{Local variance of singular vector perturbations.} By part~(i), the total perturbation variance $\mathbb{E}[(\omega_I^\top\Delta\mathbf{u}_i)^2]$ equals $\mathbb{E}[(\omega_I^\top\Delta\mathbf{u}_i^c)^2]$ up to $\mathcal{O}(n^{-2})$. The exact leading-order variance of this dominant noise-subspace component is: for any unit vectors $\omega_I \in \mathbb{R}^p$ and $\omega_J \in \mathbb{R}^n$,
    \begin{equation}\label{eq_patchwise_u}
        \mathbb{E}\!\big[(\omega_I^\top \Delta\mathbf{u}_i)^2\big] = \frac{\textup{Tr}(\mathcal{B})}{n}\bigg[\mathcal{R}_i\,\omega_I^\top\mathcal{A}\,\omega_I + \frac{\mathcal{S}_i\,\textup{Tr}(\mathcal{A})}{p}\bigg] + \mathcal{O}\!\big(n^{-2}\Delta(d_i)^{-1}\big)\,,
    \end{equation}
    \begin{equation}\label{eq_patchwise_v}
        \mathbb{E}\!\big[(\omega_J^\top \Delta\mathbf{v}_i)^2\big] = \frac{\textup{Tr}(\mathcal{A})}{p}\bigg[\mathcal{R}_i'\,\omega_J^\top\mathcal{B}\,\omega_J + \frac{\mathcal{S}_i\,\textup{Tr}(\mathcal{B})}{n}\bigg] + \mathcal{O}\!\big(n^{-2}\Delta(d_i)^{-1}\big)\,,
    \end{equation}
    where $\mathcal{R}_i'$ is defined analogously to $\mathcal{R}_i$ with the roles of $p$ and $n$ interchanged. Each formula consists of a \emph{spatially non-uniform} term governed by $\omega_I^\top\mathcal{A}\,\omega_I$ (or $\omega_J^\top\mathcal{B}\,\omega_J$) through $\mathcal{R}_i$, and a \emph{spatially uniform} term involving only the traces through $\mathcal{S}_i$.
   \item[\textup{(iii)}] \textbf{Wavelet coefficient variance.} Given the eGHWT spatial basis function $\Phi_{IJ}(x,y) = \omega_I(x)\omega_J(y)$, the variance of the wavelet coefficient error admits the decomposition:
    \begin{equation}\label{eq_sigma_exact_integral}
        \sigma_{\Phi_{IJ}}^2 := \mathbb{E}\!\left[\langle \widehat{S}_{r^+} - S_{r^+},\,\Phi_{IJ}\rangle^2\right] = \sum_{i=1}^{r^+} d_i^2 \Big( \mathcal{V}_{1,i} + \mathcal{V}_{2,i} + \mathcal{V}_{3,i} \Big) + \mathcal{O}\!\left(n^{-5/2}\Delta(d_{r^+})^{-1}\right),
    \end{equation}
    where:
    \begin{align}
        \mathcal{V}_{1,i} &= \frac{\textup{Tr}(\mathcal{B})}{n^2}\bigg[\mathcal{R}_i\,\omega_I^\top\mathcal{A}\,\omega_I + \frac{\mathcal{S}_i\,\textup{Tr}(\mathcal{A})}{p}\bigg]\,,\label{eq_V1}\\[4pt]
        \mathcal{V}_{2,i} &= \frac{\textup{Tr}(\mathcal{A})}{p^2}\bigg[\mathcal{R}_i'\,\omega_J^\top\mathcal{B}\,\omega_J + \frac{\mathcal{S}_i\,\textup{Tr}(\mathcal{B})}{n}\bigg]\,,\label{eq_V2}\\[4pt]
        \mathcal{V}_{3,i} &= np\,\mathcal{V}_{1,i}\,\mathcal{V}_{2,i}\,.\label{eq_V3}
    \end{align}
\end{itemize}
\end{theorem}
 
\begin{remark}\label{rmk_variance_structure}
Several aspects of Theorem~\ref{thm_variance} merit discussion.
\begin{itemize}
    \item[\textup{(a)}] \textbf{Noise subspace dominance.} Part~(i) shows that the signal-subspace variation $\Delta\mathbf{u}_i^r$ contributes $\mathcal{O}(n^{-2})$ in entry-wise variance, uniformly in $\Delta(d_i)$. This is because the signal resolvent $\widetilde{D}_i^r$ has eigenvalues $(\widetilde{d}_i^2 - \widetilde{d}_j^2)^{-1}$ for $j \leq r$, $j \neq i$, which are bounded by the inter-signal gaps $\mathcal{O}(1)$, and the signal row projection satisfies $\|\mathbf{e}_a^\top\widetilde{U}_r\|^2 = \mathcal{O}(n^{-1})$ by Theorem~\ref{thm_bbp_alignment}. The perturbation of $\widetilde{\mathbf{u}}_i$ is therefore overwhelmingly determined by how $\mathbf{u}_i$ leaks into the empirical noise subspace, not by how it mixes with other signal vectors. Part~(ii) provides the exact leading-order formula for this dominant component.
 
    \item[\textup{(b)}] \textbf{Two-term structure: spatially non-uniform and uniform contributions.} Part~(ii) reveals that the patchwise variance $\mathbb{E}[(\omega_I^\top\Delta\mathbf{u}_i)^2]$ consists of two terms. The first, governed by the squared Stieltjes transform $\mathcal{R}_i$ and the local noise energy $\omega_I^\top\mathcal{A}\,\omega_I$, is \emph{spatially non-uniform}: different patches $\omega_I$ experience different perturbation variances depending on the noise covariance structure. The second, governed by the resolvent $\mathcal{S}_i$ and the global averages $\frac{\textup{Tr}(\mathcal{A})}{p}$, $\frac{\textup{Tr}(\mathcal{B})}{n}$, is \emph{spatially uniform} and contributes equally across all patches.
 
    When the noise is white ($\mathcal{A} = I$, $\mathcal{B} = \sigma^2 I$), $\omega_I^\top\mathcal{A}\,\omega_I = 1 = \textup{Tr}(\mathcal{A})/p$ and $\textup{Tr}(\mathcal{B})/n = \sigma^2$, so both terms collapse into $\sigma^2(\mathcal{R}_i + \mathcal{S}_i)$, which is independent of $\omega_I$ and spatially uniform. For general $\mathcal{A}$, the $\mathcal{R}_i$ term breaks this uniformity through the factor $\omega_I^\top\mathcal{A}\,\omega_I$, making the perturbation variance spatially adaptive to the noise covariance.
 
    \item[\textup{(c)}] \textbf{Behavior of $\mathcal{R}_i$ and $\mathcal{S}_i$ across signal regimes.} The relative importance of the two terms depends on the signal strength. In the strong signal regime where $\Delta(d_i) = \mathcal{O}(1)$, both $\mathcal{R}_i$ and $\mathcal{S}_i$ are $\mathcal{O}(n^{-1})$, so the spatially non-uniform and uniform contributions are comparable and neither can be neglected. When the signal is close to the BBP transition and $\Delta(d_i) \to 0$, the Stieltjes transform derivative $|m'|$ diverges as $\mathcal{O}(\Delta(d_i)^{-1})$ due to the square-root edge of the Marchenko--Pastur law, inflating $\mathcal{S}_i$ to $\mathcal{O}(n^{-1}\Delta(d_i)^{-1})$, while the squared Stieltjes transform $m^2$ remains $\mathcal{O}(1)$, keeping $\mathcal{R}_i = \mathcal{O}(n^{-1})$. In this regime, the spatially uniform term dominates the variance.
 
    \item[\textup{(d)}] \textbf{Validity of the expansion and comparison with classical thresholds.} The leading-order variance components $\mathcal{V}_{1,i}$ and $\mathcal{V}_{2,i}$ are $\mathcal{O}(n^{-2})$ for strong signal and $\mathcal{O}(n^{-2}\Delta(d_i)^{-1})$ near BBP. The remainder $\mathcal{O}(n^{-2}\Delta(d_i)^{-1})$ in part~(ii) arises from the diagonal Isserlis correction (which involves $\frac{1}{p}\sum D_j^2\omega_I^\top\mathcal{A}\omega_I$, strictly lower order than the corresponding $\frac{\mathcal{S}_i\textup{Tr}(\mathcal{A})}{p}$ term). For the wavelet coefficient variance in part~(iii), the sum over $i = 1,\ldots,r^+$ is dominated by the weakest signal component, yielding the remainder $\mathcal{O}(n^{-5/2}\Delta(d_{r^+})^{-1})$, which is a factor $n^{-1/2}$ smaller than the leading term and confirms the expansion is valid.
 
    The coefficient-level noise standard deviation $\sigma_{\Phi_{IJ}}$ is much smaller than the raw noise level $\sigma$, because eOptShrink has already removed the bulk of the noise at the spectral level. Consequently, the adaptive wavelet threshold $\sigma_{\Phi_{IJ}}\sqrt{2\log(pn)}$ applied to the residual is strictly smaller than the Donoho--Johnstone universal threshold $\sigma\sqrt{2\log(pn)}$ that would be applied to the raw data. The e$\mathcal{O}$WS pipeline operates above the BBP threshold --- since eOptShrink requires spectral detectability --- but exploits the spatial regularity of the signal through wavelet shrinkage to recover structure that spectral methods alone cannot resolve.
 
    \item[\textup{(e)}] \textbf{Improved tree learning from the denoised estimate.} Theorem~\ref{thm_variance} provides a quantitative justification for running the Questionnaire algorithm on the initial denoised estimate rather than the raw noisy matrix $\widetilde{S}$. The Earth Mover's Distance affinities in the Questionnaire algorithm are computed via folder averages at multiple scales of the partition tree (Appendix~\ref{sec_question}, \eqref{eq:emd}). The noise variance of a folder average for a folder $I$ is of order $\sigma^2_{\mathrm{entry}} / |I|$, where $\sigma^2_{\mathrm{entry}}$ denotes the entry-level variance. For the raw noisy matrix $\widetilde{S}$, the entry-level noise variance is $\mathcal{O}(n^{-1})$. After eOptShrink, the residual error is dominated by the singular vector perturbations. By part~(ii) with $\omega_I = \mathbf{e}_a$, the entry-level perturbation variance is $\mathcal{O}(n^{-1})$ for strong signal and $\mathcal{O}(n^{-1}\Delta(d_i)^{-1})$ near BBP. The entries of the residual matrix involve cross-products $d_i(\Delta\mathbf{u}_i)_a(\mathbf{v}_i)_b$, and since $\mathbb{E}[(\mathbf{v}_i)_b^2] = \mathcal{O}(n^{-1})$ by Assumption~(v), the effective entry-level variance of the residual is $\mathcal{O}(n^{-2})$ for strong signal and $\mathcal{O}(n^{-2}\Delta(d_i)^{-1})$ near BBP. Compared to the raw noise level $\mathcal{O}(n^{-1})$, this represents at least an $\mathcal{O}(n^{-1})$ reduction in variance, making the geometric affinities significantly more accurate, particularly at fine scales where $|I|$ is small.
\end{itemize}
\end{remark}

In practice, the exact continuous-time trajectory of the eigenvalues $\widetilde{d}_k(t)$, the true singular values $d_i$, and the population spatial covariances $\mathcal{A}$ and $\mathcal{B}$ are fundamentally unobservable. We only have access to the final noisy data matrix at $t = 1/n$. However, under the high-dimensional asymptotic framework where $p \asymp n$, the theoretical variance $\sigma_{\Phi_{IJ}}^2$ can be tightly approximated using solely observable quantities. This approximation is justified by three concentration phenomena.

\emph{Concentration of resolvent integrals.}
The theoretical resolvent integrals $\mathcal{S}_i$ and $\mathcal{R}_i$ are defined as expectations of path integrals over $t \in [0, 1/n]$. In the continuous-time noise model, the support of the empirical spectral distribution expands as $t$ grows, and the spectral gap between the spiked singular value $\widetilde{d}_i(t)$ and the bulk edge monotonically shrinks. Consequently, both the integrand of $\mathcal{S}_i$ and the integrand of $\mathcal{R}_i$ are monotonically increasing in $t$ and reach their maxima at the observable endpoint $t = 1/n$. By replacing the path integrals with static evaluations at $t = 1/n$ and dividing by $n$ to account for the path length, the empirical resolvent sums $\widehat{\mathcal{S}}_i$ and $\widehat{\mathcal{R}}_i$ serve as tight upper bounds for the true integrals $\mathcal{S}_i$ and $\mathcal{R}_i$.

\emph{Concentration of quadratic forms.}
The theoretical variance components depend on the unobservable population spatial quantities $\omega_I^\top\mathcal{A}\,\omega_I \cdot \frac{\textup{Tr}(\mathcal{B})}{n}$ and $\frac{\textup{Tr}(\mathcal{A})}{p} \cdot \frac{\textup{Tr}(\mathcal{B})}{n}$. These are estimated by the empirical quadratic forms $\frac{\|\omega_I^\top\widehat{Z}\|_2^2}{n}$ and $\frac{\|\widehat{Z}\|_F^2}{pn}$, respectively. We verify the first; the others follow identically.

Recall the separable noise model $Z = \mathcal{A}^{1/2}\mathcal{X}\mathcal{B}^{1/2}$, where $\mathcal{X} \in \mathbb{R}^{p \times n}$ contains independent sub-Gaussian entries with zero mean and unit variance. Expanding the squared norm and applying the identity $\mathbb{E}[\mathcal{X}^\top M\mathcal{X}] = \textup{Tr}(M)\,I_n$ for any deterministic matrix $M$:
\begin{align}
    \mathbb{E}\!\left[\|\omega_I^\top Z\|_2^2\right] &= \mathbb{E}\!\left[\textup{Tr}\!\left(\omega_I^\top\mathcal{A}^{1/2}\mathcal{X}\mathcal{B}\mathcal{X}^\top\mathcal{A}^{1/2}\omega_I\right)\right] \nonumber\\
    &= \textup{Tr}\!\Big((\omega_I^\top\mathcal{A}\,\omega_I)\,\mathcal{B}\Big) = \textup{Tr}(\mathcal{B})\,\big(\omega_I^\top\mathcal{A}\,\omega_I\big)\,.
\end{align}
By the Hanson--Wright inequality, the empirical quadratic form $\frac{1}{n}\|\omega_I^\top Z\|_2^2$ concentrates around $\frac{\textup{Tr}(\mathcal{B})}{n}(\omega_I^\top\mathcal{A}\,\omega_I)$ with relative fluctuation $\mathcal{O}_\prec(n^{-1/2})$. Thus:
\begin{equation}
    \frac{\|\omega_I^\top\widehat{Z}\|_2^2}{n} = \frac{\textup{Tr}(\mathcal{B})}{n}\,\big(\omega_I^\top\mathcal{A}\,\omega_I\big) + \mathcal{O}_\prec(n^{-1/2})\,.
\end{equation}
Analogously, the column-wise form concentrates as $\frac{\|\widehat{Z}\omega_J\|_2^2}{p} = \frac{\textup{Tr}(\mathcal{A})}{p}(\omega_J^\top\mathcal{B}\,\omega_J) + \mathcal{O}_\prec(n^{-1/2})$, and the Frobenius norm concentrates as $\frac{\|\widehat{Z}\|_F^2}{pn} = \frac{\textup{Tr}(\mathcal{A})}{p}\cdot\frac{\textup{Tr}(\mathcal{B})}{n} + \mathcal{O}_\prec(n^{-1/2})$.

The leading-order variance components $\mathcal{V}_{1,i}$ and $\mathcal{V}_{2,i}$ are each $\mathcal{O}(n^{-2}\Delta(d_i)^{-1})$. The resolvent upper bounds introduce a multiplicative overestimate that is bounded and controlled. The quadratic form substitution introduces an additive error of $\mathcal{O}_\prec(n^{-1/2})$ in each spatial factor. Since the spatial factors multiply the $\mathcal{O}(n^{-2}\Delta(d_i)^{-1})$ variance components, the total plug-in error is $\mathcal{O}_\prec(n^{-5/2}\Delta(d_i)^{-1})$, which is a factor $n^{-1/2}$ smaller than the leading term. For the cross-term $\mathcal{V}_{3,i} = np\,\mathcal{V}_{1,i}\mathcal{V}_{2,i}$, the relative error of the product is bounded by the sum of the individual relative errors, so the combined fluctuation remains $\mathcal{O}_\prec(n^{-1/2})$. Consequently, the empirical substitution error across all three variance components is absorbed into the higher-order residuals.

Using these properties, we define the data-driven variance estimator $\widehat{\sigma}_{\Phi_{IJ}}^2$ exactly as formulated in \eqref{eq_sigma+phi2_final}. Because $\widehat{\mathcal{S}}_i$ and $\widehat{\mathcal{R}}_i$ upper bound their respective resolvent integrals, $\widehat{\sigma}_{\Phi_{IJ}}^2$ acts as a tight, conservative upper bound for $\sigma_{\Phi_{IJ}}^2$, ensuring we do not under-shrink the noise. The following proposition establishes that this empirical estimator is asymptotically consistent up to the leading order.

\begin{proposition}\label{prop_variance_bound}
Suppose (i)-(v) of Assumption \ref{assum_main} hold. Define the theoretical projection 
variance $\sigma_{\Phi_{IJ}}^2$ as in \eqref{eq_sigma_exact_integral} and the empirical 
estimator $\widehat{\sigma}_{\Phi_{IJ}}^2$ as in \eqref{eq_sigma+phi2_final}. Both 
quantities are $\mathcal{O}(n^{-2}\Delta(d_{r^+})^{-1})$. Then 
$\widehat{\sigma}_{\Phi_{IJ}}^2$ provides a computable, data-driven, and tight upper 
bound for the true projection variance:
\begin{equation}
    \sigma_{\Phi_{IJ}}^2 \le \widehat{\sigma}_{\Phi_{IJ}}^2 
    + \mathcal{O}_\prec\left( n^{-5/2}\Delta(d_{r^+})^{-1} \right)\,.
\end{equation}
\end{proposition}
\begin{remark}\label{rmk_upper_bound_choice}
While Proposition \ref{prop_variance_bound} establishes $\widehat{\sigma}_{\Phi_{IJ}}^2$ as a tight upper bound, a more direct evaluation of the theoretical variance is mathematically possible. Because the unperturbed signal singular values $d_i$ can be reliably estimated via $\widehat{d}_{e,i}$, one could theoretically compute an approximated numerical integral of the exact $\mathcal{S}_i$ and $\mathcal{R}_i$ over the path $t \in [0, 1/n]$. However, to simplify the procedure, we opt for the upper-bound approach in our algorithm. The bound $\widehat{\sigma}_{\Phi_{IJ}}^2$ offers a simple, closed-form estimator that is instantaneously computable without the need for numerical integration. Furthermore, adopting a tight but conservative upper bound is practically advantageous for the subsequent shrinkage steps, ensuring that the algorithm strictly avoids under-shrinking and robustly suppresses residual noise.
\end{remark}

\begin{theorem}\label{thm_eows_bound}
Suppose Assumption~\ref{assum_main} holds. 
We have:\\
(1) (Mean squared error of $\widehat{S}_{e\mathcal{OWS}}$.)
\begin{equation}\label{eq_ws_f2}
    \mathbb{E}_S\|S - \widehat{S}_{e\mathcal{OWS}}\|_{2,n}^2 \leq C\log^{2\alpha/(2\alpha+1)}(n)\big(L^{1/\alpha}n^{-4}\Delta(d_{r^+})^{-1}\big)^{2\alpha/(2\alpha+1)}\log_{B_U^{-1}}\big(n^4 L\,\Delta(d_{r^+})\big)(1+o(1))\,,
\end{equation}
with $C = C(\mathcal{A},\mathcal{B},B_L,B_U,\alpha)$.\\
(2) (Pointwise squared error of $\widehat{S}_{e\mathcal{OWS}}$.)\\
\noindent
For any $(x_0,y_0)\in X \times Y$, we have
\begin{equation}\label{eq_ws_pt2}
     \mathbb{E}\big(S(x_0,y_0) - \widehat{S}_{e\mathcal{OWS}}(x_0,y_0)\big)^2 \leq C\log^{2\alpha/(2\alpha+1)}(n^2)\big(L^{1/\alpha}n^{-4}\Delta(d_{r^+})^{-1}\big)^{2\alpha/(2\alpha+1)}\log^2_{B_U^{-1}}\big(n^4 L\,\Delta(d_{r^+})\big)(1+o(1))\,,
\end{equation}
with $C = C(\mathcal{A},\mathcal{B},B_L,B_U,\alpha)$.  
\end{theorem}
\begin{corollary}\label{cor_rate_comparison}
Under the conditions of Theorem~\ref{thm_eows_bound}, the following hold.
\begin{itemize}
    \item[\textup{(i)}] \textbf{Improvement over eOptShrink.} By Proposition~\ref{prop_variance_bound}, the coefficient-level noise variance after eOptShrink is $\sigma_{\Phi_{IJ}}^2 = \mathcal{O}(n^{-2}\Delta(d_{r^+})^{-1})$, which replaces the $\mathcal{O}(n^{-1})$ entry-level noise variance of the raw data. For optimal-shrinkage-based approaches, the loss functions converge to $\mathcal{O}(1)$ constants, and the corresponding MSE scales as $\mathcal{O}(n^{-2})$. When $\alpha > 1/2$ and $n$ is sufficiently large, the RHS of \eqref{eq_ws_f2} is $o(n^{-2})$, so e$\mathcal{OWS}$ strictly outperforms eOptShrink alone.
    \item[\textup{(ii)}] \textbf{Improvement over wavelet shrinkage.} The wavelet-shrinkage-only rates in \eqref{eq_ws_f} and \eqref{eq_ws_pt} use the raw noise variance $\mathcal{O}(n^{-1})$, yielding an effective noise scaling of $n^{-2}$ in the risk bound (after the $\|\cdot\|_{2,n}^2$ normalization). In contrast, the e$\mathcal{OWS}$ rates in \eqref{eq_ws_f2} and \eqref{eq_ws_pt2} use the reduced variance $\mathcal{O}(n^{-2}\Delta(d_{r^+})^{-1})$, yielding an effective noise scaling of $n^{-4}\Delta(d_{r^+})^{-1}$. When $p \asymp n$, this gives a factor of $n^{-2}\Delta(d_{r^+})$ improvement in the noise term, demonstrating that the combined pipeline substantially outperforms either method in isolation.
\end{itemize}
\end{corollary}

\section{Numerical Experiments}\label{sec_numerical}
We evaluated the performance of e\( \mathcal{OWS} \) through (1) simulated kernel matrices under various noise models, and (2) the single-channel fetal ECG (fECG) extraction task using a semi-real database. 

\subsection{Recovery of Kernel}
For simulated kernel matrices, we consider different types of noise.
Suppose $\mathcal X\in \mathbb{R}^{p\times n}$ has i.i.d. entries from a Student's t-distribution with $10$ degrees of freedom.
Set $\mathcal A = 
 Q_{\mathcal A}D_{\mathcal A}Q_{\mathcal A}^{T}\in \mathbb{R}^{p\times p}$, where $D_{\mathcal A} = \textup{diag} \{\ell_1,\ell_2,\ldots, \ell_p\}$, $Q_{\mathcal A}\in O(p)$ is generated by the QR decomposition of a random $p \times p$ matrix independent of $\mathcal X$. 
The same method is used to generate $\mathcal B=Q_{\mathcal B}D_{\mathcal B}Q_{\mathcal B}^{T}\in \mathbb{R}^{n\times n}$, which is independent of $\mathcal A$ and $\mathcal X$. The noise matrix is given by \( Z = \frac{1}{L} \mathcal{A}^{1/2} \mathcal{X} \mathcal{B}^{1/2} \), where \( L \) is a normalization factor chosen such that \( \| Z \|_{2} = \sqrt{p} \). Under this normalization, the average variance of the entries of \( Z \) is \( 1/n \).
We consider three types of noise constructions. The first one is the white noise (called TYPE1 below), where $D_{\mathcal A} = I_p$ and $D_{\mathcal B} = I_n$. 
The second one has a separable covariance structure (called TYPE2 below) with a gap in the limiting distribution, where $D_{\mathcal A}=\textup{diag}\Big\{\sqrt{1+9\times \frac{1}{p}}, \sqrt{1+9\times \frac{2}{p}}, \cdots, \sqrt{1+9\times \frac{p-1}{p}},\sqrt{10}\Big\}$ and 
$D_{\mathcal B}=\textup{diag}\Big\{ \sqrt{10+\frac{1}{n}},\sqrt{10+\frac{2}{n}},\cdots, \sqrt{10+\frac{\lfloor n/4\rfloor}{n}},$
$ \sqrt{0.3}, \cdots, \sqrt{0.3}, \sqrt{0.3} \Big\}$.
The third one (called TYPE3 below) has a more complicated separable covariance structure with $D_{\mathcal A}=\textup{diag}\Big\{\exp(\frac{1}{p}), \exp(\frac{2}{p}) \cdots, \exp(\frac{p-1}{p}),$ 
$ \exp{(1)}\Big\}$ and 
$D_{\mathcal B}=\textup{diag}\Big\{ 1.1+\sin(4\pi(\frac{1}{n})),1.1+\sin(4\pi(\frac{2}{n})),\cdots, 1.1+\sin(4\pi(\frac{n-1}{n})),\ 1.1+\sin(4\pi) \Big\}$.
The signal matrices \( S = \sum_{i=1}^r d_i \ub_i \vb_i^\top \) are generated from (1) acoustic wave propagation; (2) simulated fECG. Detailed descriptions of each data model are provided in the following sections.  

For each numerical simulation, we compare the performance of eOptShrink, \( \mathcal{WS} \), and e\( \mathcal{OWS} \). For \( \mathcal{WS} \), since the noise matrix \( Z \) is normalized such that the average variance of its entries is \( 1/n \), we apply \( \mathcal{WS} \) under the assumption that the wavelet coefficients follow a Gaussian distribution \( \mathcal{N}(0, 1/n) \), such that the shrinker is set to \( \eta_{\sqrt{2\log(pn)/n}} \). For eOptShrink, we select the shrinker that minimizes the asymptotic Frobenius loss.

For all estimators $\widehat S$,
the absolute performance of the relative to the ground truth $S$ is quantified via the scaled mean squared error (MSE):
\begin{equation}\label{eq_mse}
\| \widehat{S} - S \|_{2,n}^2 := \frac{\mathcal{L}_n^{\mathrm{fro}}(\widehat{S}, S)}{pn}
\end{equation}
where $\mathcal{L}_n^{\mathrm{fro}}(\widehat{S}, S) := \| \widehat{S} - S \|_F^2 = \sum_{i,j} ( \widehat{S}_{ij} - S_{ij} )^2$ denotes the standard Frobenius loss.
Let the estimator have SVD \( \widehat{S} = \sum_{i=1}^{p \wedge n} \widehat{d}_i \widehat{\ub}_i \widehat{\vb}_i^\top \) with singular values ordered as \( \widehat{d}_1 \geq \widehat{d}_2 \geq \cdots \geq \widehat{d}_{p \wedge n} \), we assess the recovery of the singular vectors corresponding to the smallest retained signal component \( \widehat{d}_{\widehat{r}^+} \). Specifically, we compute the projection of the estimated left singular vector \( \widehat{\ub}_{\widehat{r}^+} \) onto the true subspace \( U_{\widehat{r}^+} := [\ub_1, \ldots, \ub_{\widehat{r}^+}] \), quantified by \( \| \langle U_{\widehat{r}^+}, \widehat{\ub}_{\widehat{r}^+} \rangle \|_F \), and similarly, the projection of the estimated right singular vector \( \widehat{\vb}_{\widehat{r}^+} \) onto the subspace \( V_{\widehat{r}^+} := [\vb_1, \ldots, \vb_{\widehat{r}^+}] \), given by \( \| \langle V_{\widehat{r}^+}, \widehat{\vb}_{\widehat{r}^+} \rangle \|_F \). We refer to these as the “left inner product” and “right inner product,” respectively. These metrics serve as indicators of how accurately each estimator recovers the singular vectors associated with the smallest detectable signal component.

We perform simulations with increasing matrix sizes \( n = 256, 512, 1024,\) \(2048\), and \(4096 \), and for each \( n \), we run 10 independent trials. To assess comparative performance, we report the median with interquartile range error bars. 
Paired \( t \)-tests are conducted to evaluate statistical significance. A result is considered statistically significant if the \( p \)-value is less than 0.005.\\

\subsubsection{Acoustic Waves}
Consider an acoustic wave propagation scenario between two point clouds embedded in \( \mathbb{R}^3 \). The acoustic interaction between source points $\{x_i\}_{i=1}^p$ and target points $\{y_j\}_{j=1}^n$ is modeled by the kernel
\begin{equation}
S_{ij} = \frac{C\cos(2\pi \nu \|x_i - y_j\|)}{\|x_i - y_j\|},
\end{equation}
where \( \nu \) is the frequency parameter. This kernel corresponds to the real part of the Green's function of the three-dimensional Helmholtz equation with wavenumber \( k = 2\pi \nu \), and $C$ is a constant representing the amplitude. It captures oscillatory wave behavior with spatial decay and is commonly used to model time-harmonic acoustic interactions in free space.

\begin{figure}[ht]
    \centering
    \includegraphics[width=1.0\linewidth]{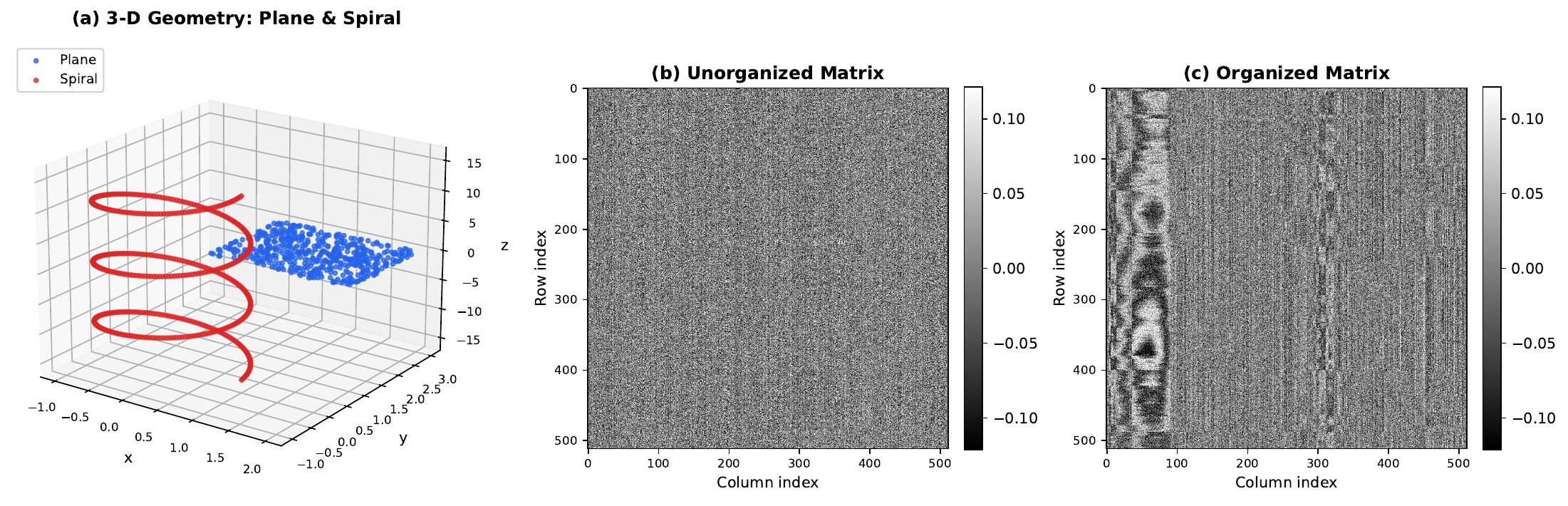}
    \caption{Helmholtz kernel matrix between a planar point cloud and a helix in $\mathbb{R}^3$.
    (a)~The three-dimensional geometry: $p=512$ points sampled uniformly on a flat plane (blue)
    and $n=512$ points along a helix of radius $r=1$ (red).
    (b)~The resulting kernel matrix $K_{ij} = C\cos(2\pi\nu\lVert x_i - y_j\rVert) / \lVert x_i - y_j\rVert$.
    (c)~The same matrix with rows and columns organized by the row and column trees
    based on the affinity from Algorithm~\ref{alg:questionnaire},
    revealing the smooth oscillatory pattern induced by the pairwise distances
    between the two point sets.}
    \label{fig:spiral_plan}
\end{figure}

We examine the scenario with equispaced helical source points \( \{x_i\}_{i=1}^p \) (red) and non-equispaced target points \( \{y_j\}_{j=1}^n \) (blue) randomly and uniformly sampled on a planar sheet, as illustrated in Figure~\ref{fig:spiral_plan}(a). Since \( \{x_i\}_{i=1}^p \) and \( \{y_j\}_{j=1}^n \) are well-separated, the kernel naturally satisfies the mixed H\"{o}lder($\alpha$) condition and is low-rank due to the regularity in the target and source domains. We set $\nu = 1$ and choose $C$ to normalize the matrix such that $\norm{S}_F^2 = 150$, with $p = n$ for $n = 256, 512, 1024, 2048$, and $4096$. Noise matrices $Z$ of TYPE1, TYPE2, and TYPE3 are added to the clean kernel $S$ to form the noisy observation $\wt{S} = S + Z$.

As illustrated in Figure~\ref{fig:spiral_plan}(b), the raw noisy kernel matrix $\wt{S}$ appears unstructured due to the random sampling scheme. Applying the Questionnaire algorithm (Algorithm~\ref{alg:questionnaire}) on $\widehat S_{\textup{eOptShrink}}$to reorganize the rows and columns reveals the underlying smooth oscillatory pattern induced by the pairwise distances between the two point sets, as shown in Figure~\ref{fig:spiral_plan}(c). This reorganization is crucial for wavelet-based denoising: in the organized matrix, the local smoothness of the kernel concentrates the signal energy into a small number of large wavelet coefficients, whereas in the unorganized matrix the energy is dispersed across many coefficients, making it difficult to separate signal from noise via thresholding.

\begin{figure}
    \centering\includegraphics[width=0.80\linewidth]{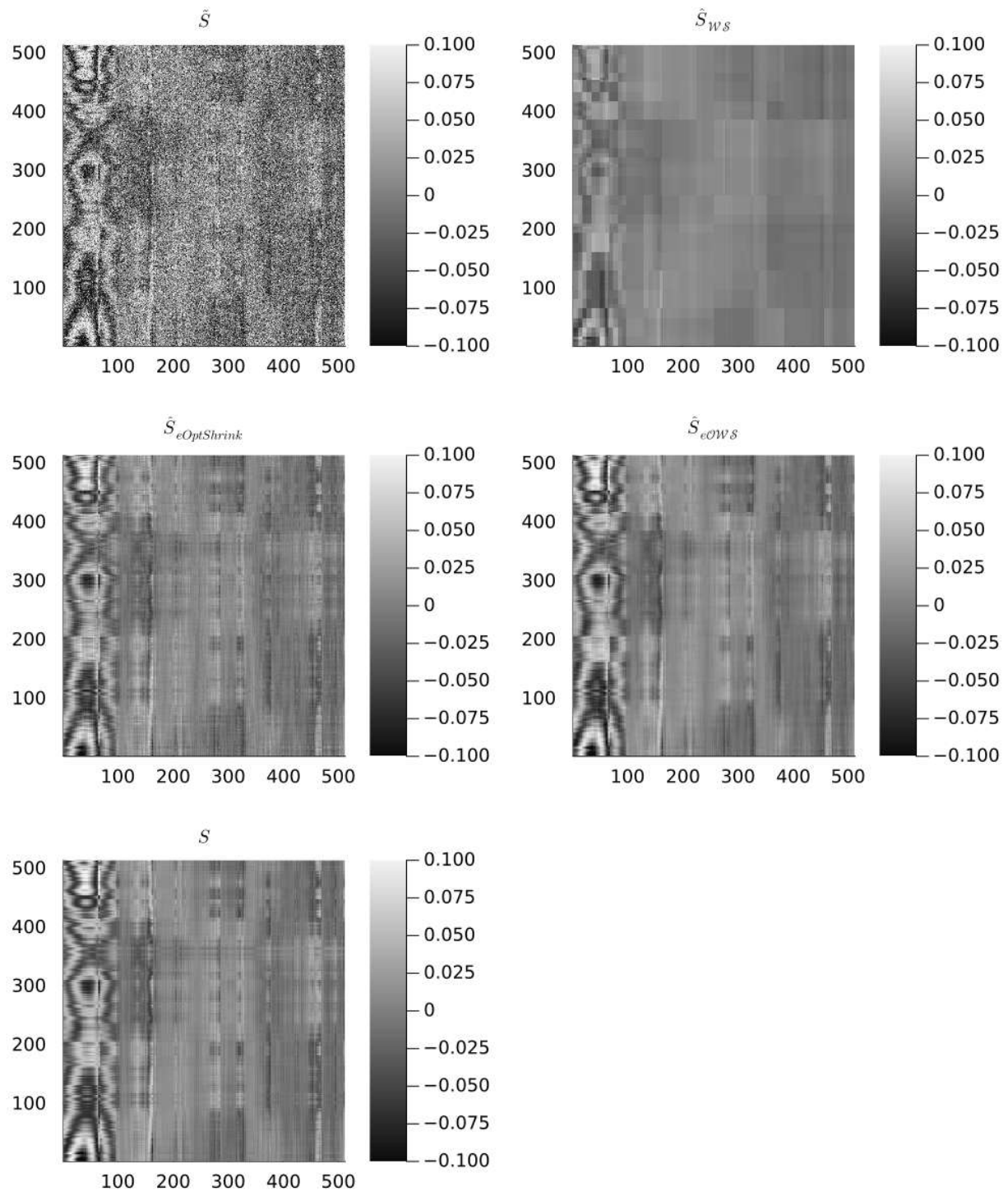}
    \caption{Comparison of denoised Helmholtz kernels. Top-left: the noisy kernel \( \widetilde{S} \) contaminated with TYPE3 noise. Top-right: denoised result using \( \mathcal{WS} \). Middle-left: denoised result using eOptShrink. Middle-right: denoised result using e\( \mathcal{OWS} \). Bottom-left: the clean ground-truth kernel \( S \).
}\label{fig:helmholtz-comparison}
\end{figure}

We apply the proposed Algorithm~\ref{alg_eqOS} to construct the estimator \( \widehat{S}_{e\mathcal{OWS}} \) from the noisy matrix \( \widetilde{S} \). Figure~\ref{fig:helmholtz-comparison} shows a comparison of our estimator \( \widehat{S}_{e\mathcal{OWS}} \) with the $\mathcal{WS}$ estimator \( \widehat{S}_{\mathcal{WS}} \), the eOptShrink estimator \( \widehat{S}_{\mathrm{eOptShrink}} \), the noisy observation \( \widetilde{S} \), and the ground truth matrix \( S \), for TYPE3 noise with $p = n = 512$. In all panels, rows and columns are reorganized by applying the Questionnaire algorithm (Algorithm~\ref{alg:questionnaire}) to $\widehat{S}_{\mathrm{eOptShrink}}$, so that the local smooth structure of each matrix is made visible.

The estimator \( \widehat{S}_{\mathcal{WS}} \) tends to oversmooth certain regions due to the use of a global variance estimate under the assumption of i.i.d.\ Gaussian noise, whereas the actual noise matrix \( Z \) is correlated and satisfies a higher-moment condition rather than a Gaussianity assumption. Both eOptShrink and our proposed e\( \mathcal{OWS} \) produce more accurate results, with e\( \mathcal{OWS} \) demonstrating superior recovery by better capturing the local piecewise-smooth structure of the signal.

\begin{figure}
\centering    \includegraphics[width=0.8\linewidth]{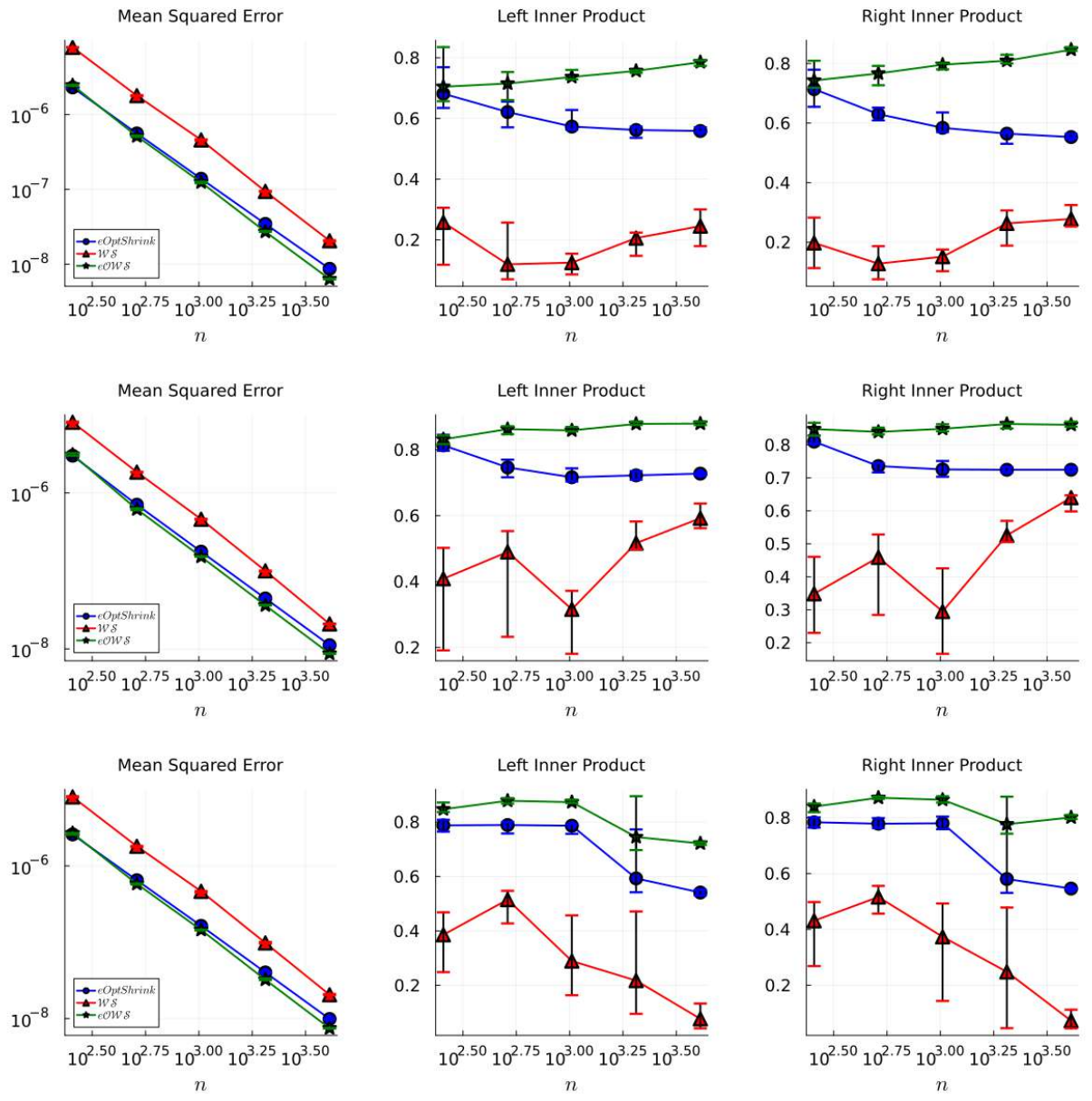}
\caption{Comparison of the denoising performance of eOptShrink, $\mathcal{WS}$, and e$\mathcal{OWS}$ on the Helmholtz kernel. The top, middle, and bottom rows correspond to the contamination with TYPE1, TYPE2, and TYPE3 noise, respectively. The first column displays the MSE with the $y$-axis on a log scale, the second column shows the left inner product, and the third column presents the right inner product. The $x$-axis in all plots represents the value of $n$ on a log scale. The blue circles mark the results of eOptShrink, the red triangles represent $\mathcal{WS}$, and the green stars denote e$\mathcal{OWS}$.
}
\label{fig:acoustic_curve}
\end{figure}

Figure~\ref{fig:acoustic_curve} shows the denoising performance of eOptShrink, $\mathcal{WS}$, and e$\mathcal{OWS}$ on the Helmholtz kernel under TYPE1, TYPE2, and TYPE3 noise. Across all values of $n$, e$\mathcal{OWS}$ achieves the lowest MSE and the highest left and right inner products compared to the other two approaches, with statistical significance. For all methods, the MSE decreases toward zero as $n$ increases, with e$\mathcal{OWS}$ exhibiting a faster convergence rate.
The values of the left and right inner products converge to $1$ as $n$ increases for TYPE1 and TYPE2 noise. Note that for TYPE3 noise, a drop in inner product values is observed at $n = 2048$. This is because, as $n$ grows, the estimated effective rank $\widehat r^+$ becomes more accurate, leading to the identification of a weaker signal component from the data matrix. Nevertheless, the inner product values ultimately converge to $1$ as $n$ increases further.

\subsubsection{Sinusoidal Waves}
In this experiment, we construct the signal matrix \( S = UDV^\top \in \mathbb{R}^{n \times 2n} \), where the entries of \( U \) are defined by sine waves:
\[
U_{ij} = \sin(2\pi j x_i), \quad \text{for } i = 1, \ldots, n, \; j = 1, \ldots, 10,
\]
with \( x_i \sim \mathcal{U}[0,1] \) sampled independently. This gives \( U \in \mathbb{R}^{n \times 10} \). The entries of \( V \in \mathbb{R}^{2n \times 10} \) are defined by cosine waves:
\[
V_{k\ell} = \cos(2\pi \ell y_k), \quad \text{for } k = 1, \ldots, 2n, \; \ell = 1, \ldots, 10,
\]
where \( y_k \sim \mathcal{U}[0,1] \) are also sampled independently. The diagonal matrix \( D = \mathrm{diag}(1, 2, \ldots, 10) \) contains the singular values. Therefore, the signal matrix \( S \) has approximate rank 10 and exhibits an underlying smooth structure, provided that \( \{x_i\}_{i=1}^n \) and \( \{y_k\}_{k=1}^{2n} \) are well organized.

\begin{figure}
    \centering\includegraphics[width=0.8\linewidth]{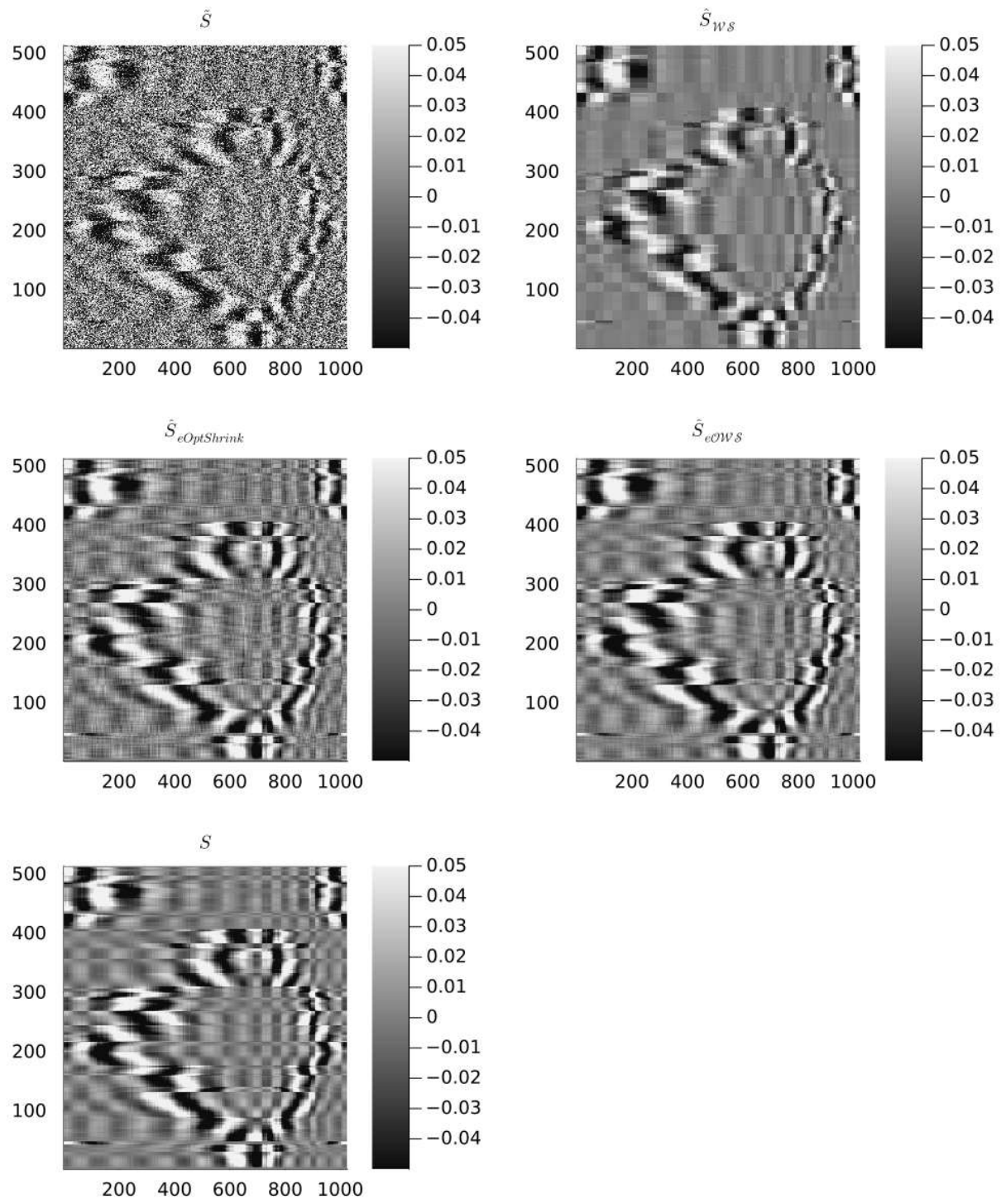}
    \caption{Comparison of denoised Sinusoidal waves. Top-left: the noisy kernel \( \widetilde{S} \) contaminated with TYPE2 noise. Top-right: denoised result using \( \mathcal{WS} \). Middle-left: denoised result using eOptShrink. Middle-right: denoised result using e\( \mathcal{OWS} \). Bottom-left: the clean ground-truth kernel \( S \).
}\label{fig:sin-comparison}
\end{figure}

\begin{figure}
\centering    \includegraphics[width=0.8\linewidth]{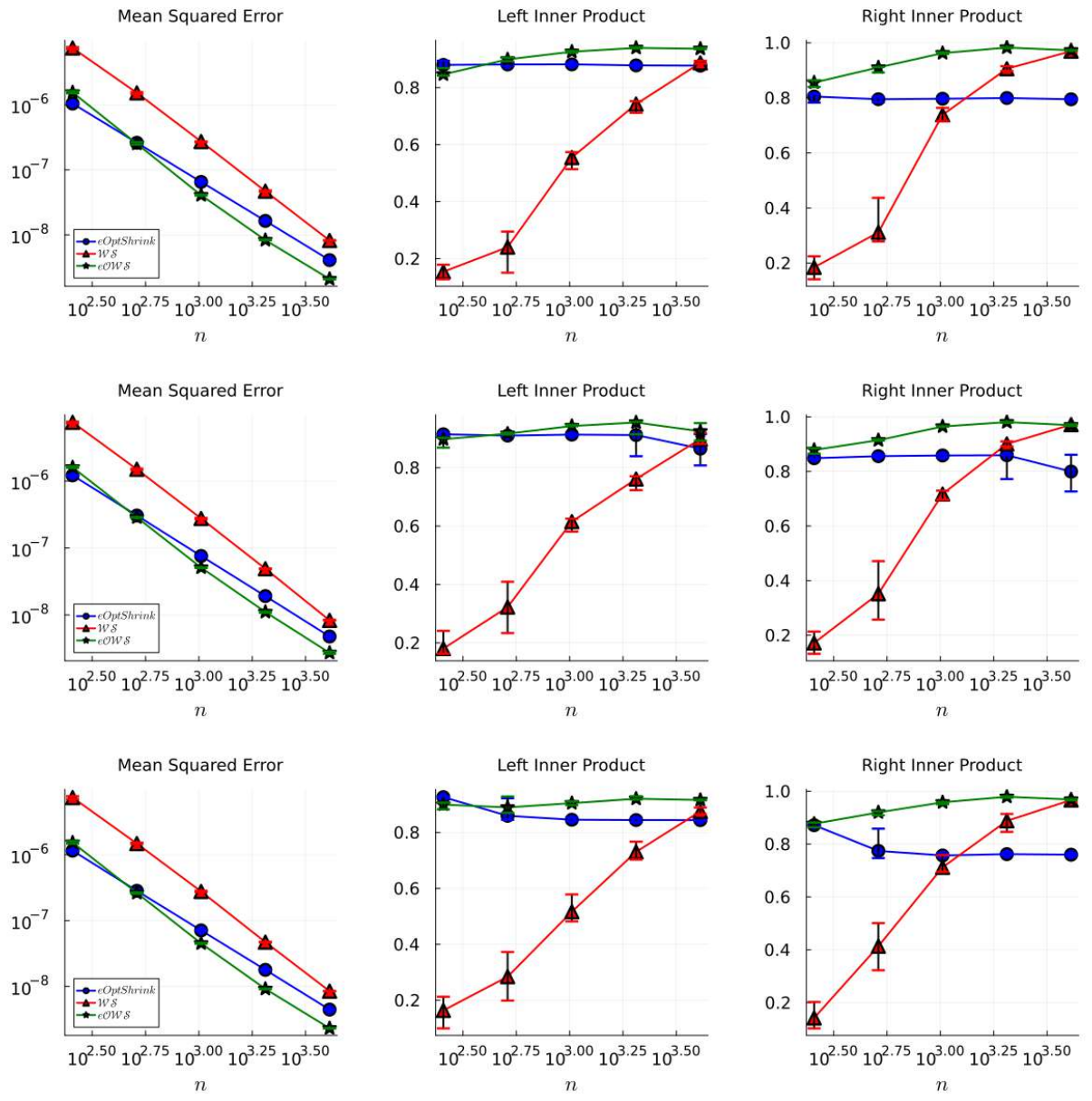}
\caption{Comparison of the denoising performance of eOptShrink, $\mathcal{WS}$, and e$\mathcal{OWS}$ on the matrix of sinusoidal waves. The top, middle, and bottom rows correspond to the contamination with TYPE1, TYPE2, and TYPE3 noise, respectively. The first column displays the MSE with the $y$-axis on a log scale, the second column shows the left inner product, and the third column presents the right inner product. The $x$-axis in all plots represents the value of $n$ on a log scale. The blue circles mark the results of eOptShrink, the red triangles represent $\mathcal{WS}$, and the green stars denote e$\mathcal{OWS}$.
}
\label{fig:sin_curve}
\end{figure}

Similar to the previous experiment, we apply the proposed Algorithm~\ref{alg_eqOS} to construct the estimator \( \widehat{S}_{e\mathcal{OWS}} \) from the noisy matrix \( \widetilde{S} \). Figure~\ref{fig:sin-comparison} presents a comparison between our estimator \( \widehat{S}_{e\mathcal{OWS}} \), the denoised matrix \( \widehat{S}_{\mathcal{WS}} \) obtained via \( \mathcal{WS} \), the eOptShrink estimator \( \widehat{S}_{eOptShrink} \), the noisy observation \( \widetilde{S} \), and the ground truth matrix \( S \), under TYPE2 noise with \( n = 512 \). The rows and columns in all matrices are reorganized using Algorithm~\ref{alg:questionnaire} applied to \( \widehat{S}_{eOptShrink} \). As shown, e\( \mathcal{OWS} \) demonstrates superior recovery performance by more effectively capturing the local piecewise-smooth structure of the signal.
Figure~\ref{fig:sin_curve} shows the denoising performance of eOptShrink, $\mathcal{WS}$, and e$\mathcal{OWS}$ on the sinusoidal waves under TYPE1, TYPE2, and TYPE3 noise. As $n$ grows, e$\mathcal{OWS}$ achieves the lowest MSE and the highest left and right inner products compared to the other two approaches, with statistical significance. 

\subsection{Fetal ECG extraction problem}
In our previous work \cite{su2019recovery}, $\mathcal{OS}$ was identified as a critical step in recovering the maternal and fetal ECG (mECG and fECG) when only one or two ta-mECG channels are available. Later, in \cite{su2025data}, we demonstrated that eOptShrink outperforms $\mathcal{OS}$ in morphology recovery by accounting for a more general separable covariance structure in the noise matrix.

The algorithm consists of two main steps. The first step is primarily designed for two-channel recordings and is omitted when only a single channel is available. The second step comprises two substeps: Step 2-1 detects maternal R peaks from the single-channel ta-mECG and is independent of $\mathcal{OS}$; in Step 2-2, the fECG is treated as noise and the mECG as the target signal. $\mathcal{OS}$ or eOptShrink is then applied to recover the mECG from the ta-mECG.
After subtracting the estimated mECG from the ta-mECG, Steps 2-1 and 2-2 are repeated on the residual signal to detect fetal R-peak locations and recover the fECG. Since both the fECG (when treated as noise) and the background noise are not white and may exhibit dependencies across segments, eOptShrink demonstrates superior performance compared to $\mathcal{OS}$ \cite{su2025data}.
We refer the reader to \cite{su2019recovery, su2020robust, su2025data} for additional details on the ECG recovery framework.

\subsubsection*{Semi-Real Simulated Database}
\label{Section more numerics simulated fECG}
To demonstrate that e$\mathcal{OWS}$ provides a further improvement in morphology recovery performance, we consider a semi-real simulated database following the same way detailed in \cite{su2019recovery}.
The data is constructed from the Physikalisch-Technische Bundesanstalt (PTB) Database \cite{bousseljot1995nutzung} (\url{https://physionet.org/physiobank/database/ptbdb/}), abbreviated as \texttt{PTBDB}. The database contains 549 recordings from 290 subjects (one to five recordings for one subject) aged 17 to 87 with the mean age 57.2. 52 out of 290 subjects are healthy. 
Each recording includes simultaneously recorded conventional 12 lead and the Frank lead ECGs. Each signal is digitalized with the sampling rate 1000 Hz. More technical details can be found online. Take 57-second Frank lead ECGs from a healthy recording, denoted as $V_x(t)$, $V_y(t)$ and $V_z(t)$ at time $t\in \mathbb{R}$, as the maternal vectocardiogram (VCG). Take $(\theta_{xy},\theta_z) = (\frac{\pi}{4},\frac{\pi}{4})$, and the simulated mECG is created by $\text{mECG}(t)=(V_x(t)  \cos{\theta_{xy}} + V_y(t)  \sin{\theta_{xy}}) \cos{\theta_{z}} + V_z(t)  \sin{\theta_{z}}$. We create {40} mECGs. The {\em simulated fECG} of healthy fetus are created from another 40 recordings from healthy subjects, where $114$-second V2 and V4 recordings are taken. The simulated and simulated fECG come from different subjects. The simulated fECG then are resampled at $500$ Hz.
As a result, the simulated fECG has about double the heart rate compared with the simulated mECG if we consider both signals sampled at $1000$ Hz. The amplitude of the simulated fECG is normalized to the same level of simulated mECG and then multiplied by $0<R<1$ shown in the second column of Table \ref{table3} to make the amplitude relationship consistent with the usual situation of real ta-mECG signals. We generate {40} simulated healthy fECGs.
The clean simulated ta-mECG is generated by directly summing simulated mECG and fECG. 
We then create a simulated noise starting with a random vector $\textbf{x} = (x_1,x_2,x_3,\ldots)$ with i.i.d entries with student t-10 distribution. The noise is then created and denotes as $\textbf{z}$ with the entries $z_i = (1+0.5\sin((i \mod 500)/500))(x_{i}+x_{i+1}+x_{i+2})$. The final simulated ta-mECG is generated by adding the created noise to the clean simulated ta-mECG according to the desired SNR ratio shown in the first column of Table \ref{table3}. As a result, for each combination of the fECG amplitude $R$ and SNR, we acquire {$40$} recordings of {$57$} seconds simulated ta-mECG signals with the sampling rate $1000$ Hz, and each recording has $5.7\times 10^4$ data points. 

We apply $\mathcal{WS}$, eOptShrink, and $e\mathcal{OWS}$ to each recording in the simulated database in Step 2-2 and compare its performance with $\mathcal{WS}$ and eOptShrink. We apply \( \mathcal{WS} \) under the assumption that the wavelet coefficients follow a Gaussian distribution \( \mathcal{N}(0, 1/n) \), such that the shrinker is set to \( \eta_{\sqrt{2\log(pn)/n}} \). eOptShrink is applied using the shrinker with respect to the Frobenius norm loss. 
For each recording, the root mean squared error (RMSE) between the reconstructed fECG cycle and the corresponding ground truth cycle is computed when the detected fetal R-peak falls within a $50$ms matching window of the true R-peak location. These RMSEs are then aggregated across 40 recordings. 

We report the mean~$\pm$~standard deviation, as presented in Table~\ref{table3}, and perform paired $t$-tests to assess the statistical significance of performance differences between methods, considering $p < 0.005$ as significant. Across varying overall SNR levels and fECG amplitudes, e$\mathcal{OWS}$ consistently achieves lower RMSEs with statistical significance. 

\begin{figure}
\centering    \includegraphics[width=1\linewidth]{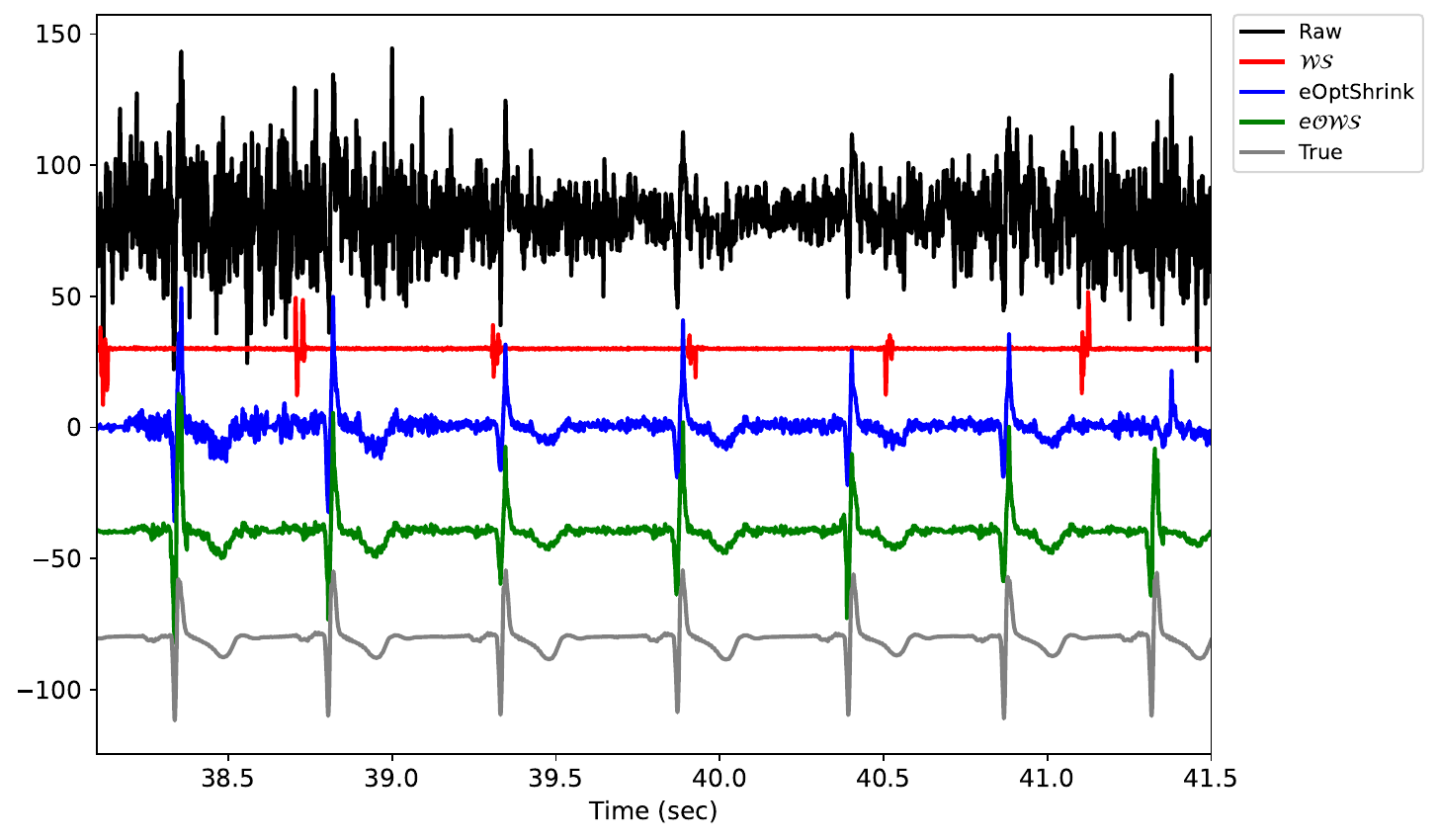}
    \caption{fECG reconstruction from a segment of the semi-real database. The black line represents the raw fECG. The red, blue, and green lines correspond to the reconstructed fECG obtained using $\mathcal{WS}$, eOptShrink, and e$\mathcal{OWS}$, respectively. The gray line indicates the ground truth fECG.}
    \label{fig:fECG_rec}
\end{figure}

Figure~\ref{fig:fECG_rec} shows a segment of raw fECG (obtained by subtracting the reconstructed maternal ECG from the original noisy ta-mECG), the reconstructed fECG from this segment, as well as the ground truth ECG, from one subject in the simulated database. 
In the raw fECG, the segment from 38s to 39.5s exhibits higher noise variation, while the segment from 39.5s to 40.5s is less noisy and clearly displays a fetal ECG cycle. All three methods, $\mathcal{WS}$, eOptShrink, and e$\mathcal{OWS}$, are applied to denoise this segment. Due to its inaccurate estimation of the noise level, $\mathcal{WS}$ oversmooths the entire segment, resulting in no discernible fECG morphology. eOptShrink provides a more accurate reconstruction of the fECG morphology; however, in the 38s to 39.5s portion, the recovered signal contains stronger residual noise compared to the 39.5s to 40.5s portion, due to increased noise in the associated singular vectors.
In contrast, e$\mathcal{OWS}$ achieves further denoising and yields a clearer morphological structure across both the noisy and cleaner segments.

\begin{table}[htb!]
\centering
{
\begin{tabular}{|c|c||c|c|c|}
\hline
SNR & $R$ & $\mathcal{WS}$ & eOptShrink\cite{su2025data} & e$\mathcal{OWS}$\\
\hline
\multirow{3}{*}{$1$ dB} & 1/4 &108.83 $\pm$ 21.63 &     44.30 $\pm$ 16.99 &  $43.31^{\star}$  $\pm$ 18.93\\
& 1/6 & 72.84 $\pm$ 15.72& 30.74 $\pm$ 12.68 & $29.83^{\star}$  $\pm$ 14.10\\
& 1/8 & 54.44 $\pm$ 11.07& 23.79 $\pm$ 10.03 & $23.13^{\star}$ $\pm$ 10.26\\
\hline
\multirow{3}{*}{$0$ dB} & 1/4 &109.44 $\pm$ 23.34&     45.81 $\pm$ 17.03 &  $43.27^{\star}$  $\pm$ 18.54\\
& 1/6 & 73.17 $\pm$ 14.64& 31.87 $\pm$ 12.11 & $30.40^{\star}$  $\pm$ 12.33\\
& 1/8 &54.66 $\pm$ 11.15& 24.65 $\pm$ 9.57 &$ 23.37^{\star}$ $\pm$ 10.16\\
\hline
\multirow{3}{*}{$-1$ dB} & 1/4 &110.79 $\pm$ 23.3 &     49.50 $\pm$ 18.51 &  $46.50^{\star}$ $\pm$ 19.02\\
& 1/6 & 73.91$\pm$ 16.50& 34.35 $\pm$ 13.50 & $32.53^{\star}$ $\pm$ 13.79\\
& 1/8 & 55.32 $\pm$ 13.08& 26.62 $\pm$ 10.93 & $25.55^{\star}$ $\pm$ 11.88\\
\hline
\end{tabular}
}
\caption{The comparison of RMSE for the fECG morphology of different algorithms applied to the simulated ta-mECG database. $R$ is the simulated fECG amplitude. All results are presented as mean $\pm$ standard deviation. 
The star next to the mean stands for the smaller RMSE with statistical significance when comparing with $\mathcal{WS}$ and eOptShrink by the paired t-test with $p<0.005$.}\label{table3}
\end{table}

\section{Discussion and Future Directions}
In this paper, we proposed the e$\mathcal{OWS}$ pipeline for recovering signal matrices that are simultaneously low-rank and mixed H\"older continuous from observations corrupted by noise with separable covariance structure. The method combines three components: eOptShrink for correcting biased singular values, the Questionnaire algorithm for learning partition trees from the denoised estimate, and scale-adaptive wavelet shrinkage via the eGHWT best basis for correcting the residual singular vector perturbation. On the theoretical side, we derived the first-order asymptotic perturbation expansion for singular vectors under the separable covariance model, obtained precise variance formulas for the eGHWT coefficients under colored noise, and established risk bounds for the e$\mathcal{OWS}$ estimator. Numerical experiments on simulated kernel matrices and a semi-real fetal ECG extraction task demonstrated that e$\mathcal{OWS}$ consistently outperforms standalone optimal shrinkage and wavelet shrinkage methods, with statistically significant improvements in both MSE and singular vector recovery.

The current pipeline is initialized via SVD and therefore inherits the BBP phase transition: signal components with singular values below the threshold $\gamma = 1/\sqrt{\mathcal{T}(\lambda_+)}$ are not recoverable by eOptShrink, as they remain embedded in the noise bulk. The Questionnaire algorithm, operating on the denoised estimate $\widehat{S}_{\textup{eOptShrink}}$, is similarly limited, since the sub-threshold components are absent from the input matrix. It is worth noting that the spectral and wavelet approaches are governed by fundamentally different detection mechanisms. The BBP transition determines whether individual singular vectors can be estimated consistently: a signal component is spectrally detectable if and only if its singular value exceeds the threshold $\gamma$, regardless of the spatial structure of the corresponding singular vectors. In contrast, direct wavelet shrinkage on the noisy matrix $\widetilde{S}$ operates on individual wavelet coefficients $\langle \widetilde{S}, \Phi_{IJ}\rangle$, each with noise variance $\sigma^2_{\Phi_{IJ}}$, and can detect signal whenever a true coefficient $\langle S, \Phi_{IJ}\rangle$ exceeds $\sigma_{\Phi_{IJ}}\sqrt{2\log(pn)}$. For a smooth signal whose energy concentrates at coarse scales, the coarse-scale wavelet coefficients can aggregate contributions from multiple singular components, and may therefore exceed the wavelet detection threshold even when no individual singular value exceeds the BBP threshold. The two thresholds thus operate on fundamentally different decompositions of the signal: spectral detection is per-component and geometry-blind, while wavelet detection is per-coefficient and exploits spatial regularity. Neither uniformly dominates the other.

A natural question, then, is whether the geometric structure of the signal could be exploited \emph{prior} to SVD to push beyond the BBP barrier. In principle, if the partition trees $\mathcal{T}_X$ and $\mathcal{T}_Y$ could be reliably estimated from the raw noisy matrix $\widetilde{S}$, wavelet shrinkage could be applied before SVD, effectively reducing the noise level and potentially lifting sub-threshold signal components above the BBP transition. However, the Questionnaire algorithm relies on Earth Mover's Distance affinities computed via folder averages over the partition tree. When applied directly to the raw noisy matrix, these averages are heavily contaminated, which can severely distort the learned tree structure and degrade the quality of the subsequent wavelet basis (Remark~\ref{rmk_variance_structure}(e)). This is precisely the motivation for the ``denoise first, then learn geometry'' ordering adopted in our pipeline.

Nevertheless, one could envision an iterative scheme that alternates between geometric and spectral denoising:
\begin{enumerate}
    \item Apply the Questionnaire algorithm to learn approximate trees from the raw noisy matrix.
    \item Apply $\mathcal{WS}$ to reduce noise.
    \item Perform SVD and $\mathcal{OS}$ on the wavelet-shrunk matrix, noting that the residual noise after wavelet shrinkage is no longer i.i.d.\ but correlated and dependent, which falls within the scope of our separable covariance framework.
    \item Refine the partition trees from the improved estimate and repeat.
\end{enumerate}
In each iteration, the improved estimate would yield more reliable trees, which in turn would produce more effective wavelet shrinkage, progressively reducing the noise level. If the noise reduction is sufficient to lift previously sub-threshold signal components above the BBP transition, additional singular components could be recovered in subsequent iterations. This is analogous in spirit to the well-known observation in the sparse PCA literature that when structural information about singular vectors is available, such as sparsity, entrywise thresholding before PCA can enable nontrivial estimation in regimes where standard spectral methods fail. In our setting, the mixed H\"older regularity serves as the structural prior, and wavelet shrinkage is the mechanism for exploiting it.

A rigorous analysis of this iterative scheme, including conditions under which convergence is guaranteed and sub-threshold recovery is achievable, is an interesting direction for future work and is beyond the scope of the current paper, which focuses on the single-pass pipeline and its theoretical guarantees above the BBP transition.
\appendix

\subsection{Proofs of the Main Results}\label{sec_proofs}

Before presenting the proofs of our main results, we first establish several fundamental properties regarding the Baik-Ben Arous-P\'ech\'e (BBP) phase transition under our high-dimensional asymptotic framework. We rely on the exact locations of the outlier eigenvalues, the rigidity of the noise bulk edge, and the geometric alignment of the empirical singular vectors. The following foundational results are adapted from \cite{su2025data}.

Note that the original results established in \cite{su2025data} were derived under a broader distribution of $\mathcal{X}$. In this work, we specifically simplify our mathematical framework by assuming $\mathcal{X}$ consists of i.i.d. Gaussian entries. Due to this, the following theorems are tailored strictly to this simplified setting, and we omit their proofs.

Let $\mathcal{T}(x)$ denote the $D$-transform (or rectangular Stieltjes transform) associated with the noise spectral distribution. We define the theoretical outlier location mapping and the asymptotic alignment factors as:
\begin{equation}
    \theta(x) := \mathcal{T}^{-1}(x^{-2})\,, \quad a_1(x) := \frac{m_{1c}(\theta(x))}{x^2 \mathcal{T}'(\theta(x))}\,, \quad a_2(x) := \frac{m_{2c}(\theta(x))}{x^2 \mathcal{T}'(\theta(x))}\,.
\end{equation}

The first theorem provides the sharp concentration rates for both the spiked signal singular values (outliers) and the leading noise singular values (bulk edge).

\begin{theorem}[Theorem 3.2 of \cite{su2025data}]\label{thm_bbp_eigenvalues}
Suppose (i)-(v) of Assumption \ref{assum_main} hold. For the signal components $1 \leq i \leq r^+$, the outlier empirical eigenvalues concentrate around their theoretical limits:
\begin{equation}
    |\widetilde{\lambda}_i - \theta(d_i)| \prec  n^{-1/2} \Delta(d_i)\,.
\end{equation}
Furthermore, for any fixed integer $\varpi > r^+$, the leading non-outlier eigenvalues adhere rigidly to the right edge of the noise bulk $\lambda_+$:
\begin{equation}
    |\widetilde{\lambda}_i - \lambda_+| \prec n^{-2/3}\,, \quad \text{for } r^+ + 1 \leq i \leq \varpi\,.
\end{equation}
\end{theorem}

To analyze the spatial geometry, we denote the projection operators onto the empirical signal subspace $\mathbf{A} = \{1, \dots, r^+\}$ for the left and right singular vectors as:
\begin{equation}
    \mathcal{P}_{\mathbf{A}} := \sum_{k=1}^{r^+} \widetilde{\mathbf{u}}_k \widetilde{\mathbf{u}}_k^\top \quad \text{and} \quad \mathcal{P}_{\mathbf{A}}' := \sum_{k=1}^{r^+} \widetilde{\mathbf{v}}_k \widetilde{\mathbf{v}}_k^\top\,.
\end{equation}

The following theorem quantifies the macroscopic overlap between the unperturbed population singular vectors and this empirical signal subspace.

\begin{theorem}[Theorem 3.5 of \cite{su2025data}]\label{thm_bbp_alignment}
Suppose (i)-(v) of Assumption \ref{assum_main} hold. For any $i,j \in \mathbf{A}$, the empirical signal subspace aligns with the true singular vectors up to the deterministic factors $a_1$ and $a_2$, bounded by:
\begin{equation}
    \big| \langle \mathbf{u}_i, \mathcal{P}_{\mathbf{A}} \mathbf{u}_j \rangle - \delta_{ij} a_1(d_i) \big| \vee \big| \langle \mathbf{v}_i, \mathcal{P}_{\mathbf{A}}' \mathbf{v}_j \rangle - \delta_{ij} a_2(d_i) \big| \prec n^{-1/2}\Delta(d_i)^{-1}\,.
\end{equation}
\end{theorem}

Finally, we require a sharp bound on the delocalization of the noise subspace. The true signal vectors must remain nearly orthogonal to the empirical noise singular vectors.

\begin{theorem}[Theorem 3.4 of \cite{su2025data}]\label{thm_bbp_delocalization}
Suppose (i)-(v) of Assumption \ref{assum_main} hold. Define the edge scaling factors $\eta_i := n^{-3/4} + n^{-5/6}i^{1/3}$ and $\varkappa_i := i^{2/3}n^{-2/3}$. For any sufficiently small constant $c > 0$ and $r^+ + 1 \leq i \leq cn$, the overlap between the true signal vectors $j \in \{1, \dots, r^+\}$ and the empirical noise vectors is bounded by:
\begin{equation}
    |\langle \mathbf{u}_j, \widetilde{\mathbf{u}}_i \rangle|^2 \vee |\langle \mathbf{v}_j, \widetilde{\mathbf{v}}_i \rangle|^2 \prec \frac{n^{-1} + \eta_i \sqrt{\varkappa_i}}{\Delta(d_j)^4 + \varkappa_i}\,.
\end{equation}
\end{theorem}

\subsubsection{Proof of Lemma \ref{lemma_exact_sde}}
\begin{proof}
Let $\widetilde{M}(t) = \widetilde{S}(t)\widetilde{S}(t)^\top$ and let $\widetilde{\lambda}_i(t) = \widetilde{d}_i(t)^2$ be the eigenvalue associated with the exact left singular vector $\widetilde{\mathbf{u}}_i(t)$. The fundamental eigenvalue equation at time $t$ is:
\begin{equation}\label{eq_eigen_u_t}
\widetilde{M}(t)\widetilde{\mathbf{u}}_i(t) = \widetilde{\lambda}_i(t)\widetilde{\mathbf{u}}_i(t).
\end{equation}
Applying It\^o's product rule, the differential of $\widetilde{M}(t)$ is given by $d\widetilde{M}(t) = (dZ(t))\widetilde{S}(t)^\top + \widetilde{S}(t)(dZ(t))^\top + \mathcal{O}(dt)$, where the $\mathcal{O}(dt)$ term encapsulates the quadratic variation of the noise. Applying It\^o's product rule to \eqref{eq_eigen_u_t} and isolating the leading-order martingale components yields:
\begin{equation}\label{eq_diff_eigen_t}
(d\widetilde{M}(t))\widetilde{\mathbf{u}}_i(t) + \widetilde{M}(t) d\widetilde{\mathbf{u}}_i(t) = d\widetilde{\lambda}_i(t) \widetilde{\mathbf{u}}_i(t) + \widetilde{\lambda}_i(t) d\widetilde{\mathbf{u}}_i(t) + \mathcal{O}(dt).
\end{equation}

Since the exact singular vectors $\widetilde{U}(t) = [\widetilde{\mathbf{u}}_1(t), \ldots, \widetilde{\mathbf{u}}_p(t)]$ form an orthonormal basis for $\mathbb{R}^p$ at any time $t$, we project the differential onto this basis: $d\widetilde{\mathbf{u}}_i(t) = \sum_{j=1}^p c_{ij}(t) \widetilde{\mathbf{u}}_j(t)$, with $c_{ij}(t) = \widetilde{\mathbf{u}}_j(t)^\top d\widetilde{\mathbf{u}}_i(t)$. 

For $j \neq i$, multiplying \eqref{eq_diff_eigen_t} from the left by $\widetilde{\mathbf{u}}_j(t)^\top$ and utilizing the orthogonality and eigenvalue properties yields:
\begin{equation}
\widetilde{\mathbf{u}}_j(t)^\top (d\widetilde{M}(t))\widetilde{\mathbf{u}}_i(t) + \widetilde{\lambda}_j(t) c_{ij}(t) = \widetilde{\lambda}_i(t) c_{ij}(t) + \mathcal{O}(dt).
\end{equation}
Rearranging to solve for $c_{ij}(t)$ and applying the SVD identities $\widetilde{S}(t)^\top \widetilde{\mathbf{u}}_i(t) = \widetilde{d}_i(t) \widetilde{\mathbf{v}}_i(t)$ and $\widetilde{\mathbf{u}}_j(t)^\top \widetilde{S}(t) = \widetilde{d}_j(t) \widetilde{\mathbf{v}}_j(t)^\top$, we obtain the scalar coefficients:
\begin{equation}\label{eq_cij_simp_t}
c_{ij}(t) = \frac{\widetilde{d}_i(t) \widetilde{\mathbf{u}}_j(t)^\top (dZ(t)) \widetilde{\mathbf{v}}_i(t) + \widetilde{d}_j(t) \widetilde{\mathbf{v}}_j(t)^\top (dZ(t))^\top \widetilde{\mathbf{u}}_i(t)}{\wt\lambda_i(t) - \wt\lambda_j(t)} + \mathcal{O}(dt).
\end{equation}
For $j=i$, applying It\^o's lemma to the constraint $\|\widetilde{\mathbf{u}}_i(t)\|^2 = 1$ demonstrates that the martingale part of $c_{ii}(t)$ must vanish, leaving purely a drift term $c_{ii}(t) = \mathcal{O}(dt)$.

We now construct the full vector $d\widetilde{\mathbf{u}}_i(t)$ by partitioning the summation $\sum_{j \neq i} c_{ij}(t) \widetilde{\mathbf{u}}_j(t)$ into the signal subspace ($j \le r$) and the noise subspace ($j > r$).

\textbf{Signal Subspace Component $d\widetilde{\mathbf{u}}_i^r(t)$:} 
For indices $j \le r$ (with $j \neq i$), we collect the basis vectors into the matrix $\widetilde{U}_r(t)$. By defining the $r \times r$ diagonal matrix $\widetilde{D}_i^r(t)$ with elements $1/(\wt\lambda_i(t) - \wt\lambda_j(t))$ and a zero at the $i$-th position, the scalar summation over $j \le r$ is mathematically equivalent to the block-matrix multiplication:
\begin{equation}
\widetilde{d}_i(t) \widetilde{U}_r(t) \widetilde{D}_i^r(t) \widetilde{U}_r(t)^\top (dZ(t)) \widetilde{\mathbf{v}}_i(t) + \widetilde{U}_r(t) \widetilde{D}_i^r(t) \widetilde{\Sigma}_r(t) \widetilde{V}_r(t)^\top (dZ(t))^\top \widetilde{\mathbf{u}}_i(t),
\end{equation}
which exactly establishes \eqref{eq_exact_du_r}.

\textbf{Noise Subspace Component $d\widetilde{\mathbf{u}}_i^c(t)$:}
For indices $j > r$, we collect the corresponding basis vectors into the $p \times (p-r)$ matrix $\widetilde{U}_c(t)$. We define $\widetilde{D}_i^c(t)$ as the $(p-r) \times (p-r)$ diagonal matrix containing the denominators $1/(\widetilde{d}_i(t)^2 - \widetilde{d}_j(t)^2)$. The noise singular values $\widetilde{d}_j(t)$ natively form the diagonal matrix $\widetilde{\Sigma}_c(t)$. Consequently, the summation over the noise indices factors perfectly into the symmetric block-matrix form:
\begin{equation}
\widetilde{d}_i(t) \widetilde{U}_c(t) \widetilde{D}_i^c(t) \widetilde{U}_c(t)^\top (dZ(t)) \widetilde{\mathbf{v}}_i(t) + \widetilde{U}_c(t) \widetilde{D}_i^c(t) \widetilde{\Sigma}_c(t) \widetilde{V}_c(t)^\top (dZ(t))^\top \widetilde{\mathbf{u}}_i(t),
\end{equation}
which exactly establishes \eqref{eq_exact_du_c}.

Appending the cumulative $\mathcal{O}(dt)$ drift vector from $c_{ii}(t)$ and the cross-variations concludes the proof for $d\widetilde{\mathbf{u}}_i(t)$. By symmetry, applying the identical It\^o differentiation procedure to the right eigenvalue equation $\widetilde{S}(t)^\top\widetilde{S}(t) \widetilde{\mathbf{v}}_i(t) = \widetilde{\lambda}_i(t)\widetilde{\mathbf{v}}_i(t)$ yields the symmetric decomposition for the right singular vector, establishing \eqref{eq_exact_dv_r} and \eqref{eq_exact_dv_c}.
\end{proof}

\subsection{Proof of Theorem~\ref{thm_variance}}
\begin{proof}[Proof of Theorem~\ref{thm_variance}]

\textbf{Step 1: Perturbation decomposition and noise subspace dominance.}

We begin by expanding the rank-1 estimation error. For each signal component $i$, the outer-product perturbation decomposes as:
\begin{equation}\label{eq_outerproduct_expand}
    \widetilde{\mathbf{u}}_i\widetilde{\mathbf{v}}_i^\top - \mathbf{u}_i\mathbf{v}_i^\top = \Delta\mathbf{u}_i\,\mathbf{v}_i^\top + \mathbf{u}_i\,\Delta\mathbf{v}_i^\top + \Delta\mathbf{u}_i\,\Delta\mathbf{v}_i^\top\,,
\end{equation}
where $\Delta\mathbf{u}_i := \widetilde{\mathbf{u}}_i - \mathbf{u}_i$ and $\Delta\mathbf{v}_i := \widetilde{\mathbf{v}}_i - \mathbf{v}_i$. Projecting onto the separable wavelet basis element $\Phi_{IJ} = \omega_I\omega_J^\top$ via the trace inner product, the total wavelet coefficient residual becomes:
\begin{equation}
    \langle \widehat{S}_{r^+} - S_{r^+},\,\Phi_{IJ}\rangle = \sum_{i=1}^{r^+} d_i\,\mathcal{E}_i\,,
\end{equation}
where the projection error associated with the $i$-th component is:
\begin{equation}
    \mathcal{E}_i := (\omega_I^\top\Delta\mathbf{u}_i)(\mathbf{v}_i^\top\omega_J) + (\omega_I^\top\mathbf{u}_i)(\Delta\mathbf{v}_i^\top\omega_J) + (\omega_I^\top\Delta\mathbf{u}_i)(\Delta\mathbf{v}_i^\top\omega_J)\,.
\end{equation}

By Lemma~\ref{lemma_exact_sde}, each singular vector perturbation decomposes into orthogonal subspace components:
\begin{equation}
    \Delta\mathbf{u}_i = \Delta\mathbf{u}_i^c + \Delta\mathbf{u}_i^r + \mathcal{O}(n^{-2})\,,
\end{equation}
and analogously for $\Delta\mathbf{v}_i$, where the $\mathcal{O}(n^{-2})$ residual arises from integrating the higher-order It\^o drift over the path $[0,1/n]$. We now establish that $\Delta\mathbf{u}_i^r$ is negligible.

Applying It\^o's isometry to the exact differential \eqref{eq_exact_du_r} for the signal-subspace component, the entry-wise quadratic variation involves the signal row projection $\|\mathbf{e}_a^\top\widetilde{U}_r(t)\|_2^2 = \mathcal{O}(n^{-1})$ by the delocalization of the signal vectors from Assumption~(v) and Theorem~\ref{thm_bbp_alignment}, and the signal resolvent $\widetilde{D}_i^r(t)$ whose eigenvalues $(\widetilde{d}_i^2 - \widetilde{d}_j^2)^{-1}$ for $j \leq r$, $j \neq i$ are bounded by the inter-signal gaps $\mathcal{O}(1)$. Since $\|\mathcal{A}\|_2, \|\mathcal{B}\|_2 = \mathcal{O}(1)$ and the signal vectors are unit-norm:
\begin{equation}\label{eq_ur_order}
    \mathbb{E}[(\Delta \mathbf{u}_i^r)_a^2] = \mathcal{O}(n^{-2})\,,
\end{equation}
uniformly in $\Delta(d_i)$. This establishes part~(i), and we may restrict the leading-order analysis to:
\begin{equation}\label{eq_Eic}
    \mathcal{E}_i^c := (\omega_I^\top\Delta\mathbf{u}_i^c)(\mathbf{v}_i^\top\omega_J) + (\omega_I^\top\mathbf{u}_i)(\Delta\mathbf{v}_i^{c\top}\omega_J) + (\omega_I^\top\Delta\mathbf{u}_i^c)(\Delta\mathbf{v}_i^{c\top}\omega_J)\,,
\end{equation}
with the signal-subspace contributions absorbed into the $\mathcal{O}(n^{-2})$ remainder.

\medskip
\textbf{Step 2: Squaring the error and elimination of cross-terms.}

Squaring the wavelet coefficient residual and taking expectations:
\begin{equation}\label{eq_squared_expansion}
    \mathbb{E}\!\left[\langle \widehat{S}_{r^+}-S_{r^+},\,\Phi_{IJ}\rangle^2\right] = \sum_{i=1}^{r^+} d_i^2\,\mathbb{E}\big[(\mathcal{E}_i^c)^2\big] + \sum_{i\neq j} d_i d_j\,\mathbb{E}\big[\mathcal{E}_i^c\,\mathcal{E}_j^c\big] + \mathcal{O}\!\big(n^{-5/2}\Delta(d_{r^+})^{-1}\big)\,.
\end{equation}

\emph{Inter-rank cross-terms ($i \neq j$).}
Expanding $\mathbb{E}[\mathcal{E}_i^c\,\mathcal{E}_j^c]$ generates products of noise-driven perturbations from distinct signal components. At leading order, the quadratic terms are proportional to $\mathbf{v}_i^\top\mathbf{v}_j = 0$ and $\mathbf{u}_i^\top\mathbf{u}_j = 0$, and therefore vanish by the orthogonality of distinct signal vectors. The fourth-order cross-terms involve paired covariances that likewise vanish. Hence the inter-rank sum is zero at leading order.

\emph{Intra-rank squared error ($i = j$).}
Expanding $(\mathcal{E}_i^c)^2$ produces three squared variance components and three cross-products:
\begin{alignat}{2}
    &\mathcal{V}_{1,i} := \mathbb{E}\!\big[(\omega_I^\top\Delta\mathbf{u}_i^c)^2(\mathbf{v}_i^\top\omega_J)^2\big]\,,\quad &&\mathcal{V}_{2,i} := \mathbb{E}\!\big[(\omega_I^\top\mathbf{u}_i)^2(\Delta\mathbf{v}_i^{c\top}\omega_J)^2\big]\,,\\[2pt]
    &\mathcal{V}_{3,i} := \mathbb{E}\!\big[(\omega_I^\top\Delta\mathbf{u}_i^c)^2(\Delta\mathbf{v}_i^{c\top}\omega_J)^2\big]\,,
\end{alignat}
and the cross-products $\mathcal{T}_{1,i}$, $\mathcal{T}_{2,i}$, $\mathcal{T}_{3,i}$.

\begin{itemize}
    \item \emph{Left--right interaction ($\mathcal{T}_{1,i}$).}
    $\mathcal{T}_{1,i} = 2(\mathbf{v}_i^\top\omega_J)(\omega_I^\top\mathbf{u}_i)\,\mathbb{E}[(\omega_I^\top\Delta\mathbf{u}_i^c)(\Delta\mathbf{v}_i^{c\top}\omega_J)]$ couples the left and right noise-subspace perturbations. Applying It\^o's product formula, its expectation reduces to the integral of the cross-variation rate, which contains $\widetilde{U}_c(t)^\top(\mathcal{A}-I)\widetilde{\mathbf{u}}_i(t) = \mathcal{O}_\prec(n^{-1/2})$ by isotropy. Multiplying by the projection magnitudes $(\mathbf{v}_i^\top\omega_J)(\omega_I^\top\mathbf{u}_i) = \mathcal{O}(n^{-1})$ from Assumption~(v) and integrating over $[0,1/n]$:
    \begin{equation}
        \mathbb{E}[\mathcal{T}_{1,i}] = \mathcal{O}(n^{-3})\,.
    \end{equation}

    \item \emph{Cubic terms ($\mathcal{T}_{2,i}$ and $\mathcal{T}_{3,i}$).}
    Define $X_t := \omega_I^\top\Delta\widetilde{\mathbf{u}}_i^c(t)$ and $Y_t := \Delta\widetilde{\mathbf{v}}_i^{c\top}(t)\omega_J$. By the multidimensional It\^o formula:
    \begin{equation}
        d(X_t^2 Y_t) = 2X_t Y_t\,dX_t + X_t^2\,dY_t + Y_t\,d\langle X,X\rangle_t + 2X_t\,d\langle X,Y\rangle_t\,.
    \end{equation}
    The $dX_t$ and $dY_t$ terms are local martingales and vanish under expectation. Since $Y_t$ is a martingale with $\mathbb{E}[Y_t] = 0$, the finite-variation terms contribute only through the covariance between $Y_t$ and the quadratic variation rate, both $\mathcal{O}(n^{-1/2})$. Integrating over $[0,1/n]$ and multiplying by $\mathbf{v}_i^\top\omega_J = \mathcal{O}(n^{-1/2})$:
    \begin{equation}
        \mathbb{E}[\mathcal{T}_{2,i}] = \mathcal{O}(n^{-5/2})\,.
    \end{equation}
    By left-right symmetry, $\mathbb{E}[\mathcal{T}_{3,i}] = \mathcal{O}(n^{-5/2})$.
\end{itemize}

All cross-products lie within the remainder, and:
\begin{equation}\label{eq_squared_marginal}
    \mathbb{E}\!\left[\langle \widehat{S}_{r^+}-S_{r^+},\,\Phi_{IJ}\rangle^2\right] = \sum_{i=1}^{r^+} d_i^2 \big(\mathcal{V}_{1,i} + \mathcal{V}_{2,i} + \mathcal{V}_{3,i}\big) + \mathcal{O}\!\big(n^{-5/2}\Delta(d_{r^+})^{-1}\big)\,.
\end{equation}

\medskip
\textbf{Step 3: Evaluating the variance components via It\^o isometry and Isserlis' theorem.}

\emph{Computation of $\mathbb{E}[(\omega_I^\top\Delta\mathbf{u}_i^c)^2]$.}
Substituting the exact differential from \eqref{eq_exact_du_c} and writing $dZ(t) = \mathcal{A}^{1/2}\,d\mathcal{X}(t)\,\mathcal{B}^{1/2}$, we define the driver matrices:
\begin{align}
    M_{1,i}(t) &:= \widetilde{d}_i(t)\,\widetilde{U}_c(t)\widetilde{D}_i^{c,u}(t)\widetilde{U}_c(t)^\top \;\in\mathbb{R}^{p\times p}\,,\\
    M_{2,i}(t) &:= \widetilde{U}_c(t)\widetilde{D}_i^{c,u}(t)\widetilde{\Sigma}_c(t)\widetilde{V}_c(t)^\top \;\in\mathbb{R}^{p\times n}\,,
\end{align}
where $\widetilde{D}_i^{c,u}(t) = (\widetilde{d}_i(t)^2 I_{p-r} - \widetilde{\Sigma}_c(t)\widetilde{\Sigma}_c(t)^\top)^{-1}$ is the noise-subspace resolvent from Lemma~\ref{lemma_exact_sde}. Projecting onto $\omega_I$:
\begin{equation}
    d\!\big(\omega_I^\top\widetilde{\mathbf{u}}_i^c(t)\big) = \big(\omega_I^\top M_{1,i}(t)\big)\,dZ(t)\,\widetilde{\mathbf{v}}_i(t) + \big(\omega_I^\top M_{2,i}(t)\big)\,dZ(t)^\top\widetilde{\mathbf{u}}_i(t)\,.
\end{equation}

By It\^o's isometry, the expected quadratic variation decomposes into two diagonal terms and a cross-term. The cross-term contains the factor $\widetilde{U}_c(t)^\top\mathcal{A}\,\widetilde{\mathbf{u}}_i(t) = \widetilde{U}_c(t)^\top(\mathcal{A}-I)\widetilde{\mathbf{u}}_i(t)$, which is $\mathcal{O}_\prec(n^{-1/2})$ entry-wise by isotropy. After weighting by the resolvent and integrating over $[0,1/n]$, this contributes $\mathcal{O}(n^{-2}\Delta(d_i)^{-1})$, absorbed into the remainder. The diagonal terms give:
\begin{equation}\label{eq_ito_expanded}
    \mathbb{E}\!\big[(\omega_I^\top\Delta\mathbf{u}_i^c)^2\big] = \int_0^{1/n}\!\mathbb{E}\!\Big[\underbrace{\big(\omega_I^\top M_{1,i}\mathcal{A}\,M_{1,i}^\top\omega_I\big)\big(\widetilde{\mathbf{v}}_i^\top\mathcal{B}\,\widetilde{\mathbf{v}}_i\big)}_{\text{Term 1}} + \underbrace{\big(\omega_I^\top M_{2,i}\mathcal{B}\,M_{2,i}^\top\omega_I\big)\big(\widetilde{\mathbf{u}}_i^\top\mathcal{A}\,\widetilde{\mathbf{u}}_i\big)}_{\text{Term 2}}\Big]\,dt + \mathcal{O}\!\big(n^{-2}\Delta(d_i)^{-1}\big)\,.
\end{equation}

The signal vector quadratic forms concentrate by the Hanson--Wright inequality:
\begin{equation}\label{eq_signal_conc}
    \widetilde{\mathbf{v}}_i^\top\mathcal{B}\,\widetilde{\mathbf{v}}_i = \frac{\textup{Tr}(\mathcal{B})}{n} + \mathcal{O}_\prec(n^{-1/2})\,,\qquad \widetilde{\mathbf{u}}_i^\top\mathcal{A}\,\widetilde{\mathbf{u}}_i = \frac{\textup{Tr}(\mathcal{A})}{p} + \mathcal{O}_\prec(n^{-1/2})\,.
\end{equation}

\emph{Evaluation of Term 2.} We have $M_{2,i} = \widetilde{U}_c\widetilde{D}_i^{c,u}\widetilde{\Sigma}_c\widetilde{V}_c^\top$. By isotropy of $\widetilde{V}_c$ (independent of $\widetilde{U}_c$), the inner matrix concentrates: $\widetilde{V}_c^\top\mathcal{B}\,\widetilde{V}_c = \frac{\textup{Tr}(\mathcal{B})}{n}I_{n-r} + \mathcal{O}_\prec(n^{-1/2})$. The remaining form $\omega_I^\top\widetilde{U}_c(\cdots)\widetilde{U}_c^\top\omega_I$ is degree 2 in $\widetilde{U}_c$ and concentrates by the Hanson--Wright inequality:
\begin{equation}\label{eq_term2}
    \omega_I^\top M_{2,i}\mathcal{B}\,M_{2,i}^\top\omega_I = \frac{\textup{Tr}(\mathcal{B})}{n}\cdot\frac{1}{p}\sum_{j=1}^{p-r} \frac{\widetilde{d}_{r+j}(t)^2}{(\widetilde{d}_i(t)^2 - \widetilde{d}_{r+j}(t)^2)^2} + \mathcal{O}_\prec\!\big(n^{-1/2}\Delta(d_i)^{-1}\big)\,.
\end{equation}

\emph{Evaluation of Term 1.} Since $M_{1,i} = \widetilde{d}_i\widetilde{U}_c\widetilde{D}_i^{c,u}\widetilde{U}_c^\top$, the quadratic form $\omega_I^\top M_{1,i}\mathcal{A}\,M_{1,i}^\top\omega_I$ is degree 4 in $\widetilde{U}_c$. Since $\widetilde{D}_i^{c,u} = (\widetilde{d}_i^2 I_{p-r} - \widetilde{\Sigma}_c\widetilde{\Sigma}_c^\top)^{-1}$ is diagonal in the eigenbasis of $\widetilde{\Sigma}_c\widetilde{\Sigma}_c^\top$ with eigenvalues $D_j := (\widetilde{d}_i^2 - \widetilde{d}_{r+j}^2)^{-1}$, $j = 1,\ldots,p-r$, we expand in this basis:
\begin{equation}
    \omega_I^\top M_{1,i}\mathcal{A}\,M_{1,i}^\top\omega_I = \widetilde{d}_i^2\sum_{j,k=1}^{p-r} D_j D_k\,(\widetilde{\mathbf{u}}_{r+j}^\top\omega_I)(\widetilde{\mathbf{u}}_{r+k}^\top\omega_I)(\widetilde{\mathbf{u}}_{r+j}^\top\mathcal{A}\,\widetilde{\mathbf{u}}_{r+k})\,.
\end{equation}
We evaluate using Isserlis' theorem, treating the noise eigenvectors as independent isotropic vectors with entry variance $1/p$.

\textit{Off-diagonal ($j \neq k$).} Conditioning on $\widetilde{\mathbf{u}}_{r+k}$ and using $\mathbb{E}[\widetilde{\mathbf{u}}_{r+j}\widetilde{\mathbf{u}}_{r+j}^\top] = \frac{1}{p}I_p$:
\begin{equation}
    \mathbb{E}\big[(\widetilde{\mathbf{u}}_{r+j}^\top\omega_I)(\widetilde{\mathbf{u}}_{r+k}^\top\omega_I)(\widetilde{\mathbf{u}}_{r+j}^\top\mathcal{A}\,\widetilde{\mathbf{u}}_{r+k})\big] = \frac{\omega_I^\top\mathcal{A}\,\omega_I}{p^2}\,.
\end{equation}
Summing over $j \neq k$:
\begin{equation}\label{eq_offdiag}
    \text{Off-diagonal} = \frac{\widetilde{d}_i^2}{p^2}\bigg[\bigg(\sum_j D_j\bigg)^{\!2} - \sum_j D_j^2\bigg]\omega_I^\top\mathcal{A}\,\omega_I\,.
\end{equation}

\textit{Diagonal ($j = k$).} For $\mathbf{g}\sim\mathcal{N}(0,\frac{1}{p}I)$, Isserlis' theorem gives three pairings:
\begin{equation}
    \mathbb{E}[(\mathbf{g}^\top\omega_I)^2(\mathbf{g}^\top\mathcal{A}\,\mathbf{g})] = \frac{\textup{Tr}(\mathcal{A}) + 2\,\omega_I^\top\mathcal{A}\,\omega_I}{p^2}\,.
\end{equation}
Summing over $j$:
\begin{equation}\label{eq_diag}
    \text{Diagonal} = \frac{\widetilde{d}_i^2}{p^2}\sum_j D_j^2\,\big(\textup{Tr}(\mathcal{A}) + 2\,\omega_I^\top\mathcal{A}\,\omega_I\big)\,.
\end{equation}

Combining \eqref{eq_offdiag} and \eqref{eq_diag}:
\begin{equation}\label{eq_term1_full}
    \mathbb{E}\big[\omega_I^\top M_{1,i}\mathcal{A}\,M_{1,i}^\top\omega_I\big] = \frac{\widetilde{d}_i^2}{p^2}\bigg[\bigg(\bigg(\sum_j D_j\bigg)^{\!2} + \sum_j D_j^2\bigg)\omega_I^\top\mathcal{A}\,\omega_I + \sum_j D_j^2\,\textup{Tr}(\mathcal{A})\bigg]\,.
\end{equation}

\emph{Combining Terms 1 and 2.} Multiplying Term 1 by $\frac{\textup{Tr}(\mathcal{B})}{n}$ and Term 2 by $\frac{\textup{Tr}(\mathcal{A})}{p}$, we collect terms by their spatial dependence. The off-diagonal Isserlis contribution gives the dominant $\omega_I^\top\mathcal{A}\omega_I$ coefficient:
\begin{equation}
    \frac{\textup{Tr}(\mathcal{B})}{n}\cdot\frac{\widetilde{d}_i^2}{p^2}\bigg(\sum_j D_j\bigg)^{\!2}\omega_I^\top\mathcal{A}\,\omega_I\,,
\end{equation}
where $(\sum D_j)^2/p^2 = (\frac{1}{p}\sum D_j)^2 = m(\widetilde{d}_i^2)^2 = \mathcal{O}(1)$. The diagonal Isserlis correction $\frac{\widetilde{d}_i^2}{p^2}\sum D_j^2\omega_I^\top\mathcal{A}\omega_I$ is $\mathcal{O}(n^{-1}\Delta(d_i)^{-1})$, which is strictly lower order than the off-diagonal term $\mathcal{O}(1)$, and is absorbed into the remainder after integration.

The $\textup{Tr}(\mathcal{A})$ coefficient combines the Term 1 diagonal and all of Term 2:
\begin{equation}
    \frac{\textup{Tr}(\mathcal{B})}{n}\cdot\frac{1}{p}\sum_j D_j^2\big(\widetilde{d}_i^2 + \widetilde{d}_{r+j}^2\big)\cdot\frac{\textup{Tr}(\mathcal{A})}{p}\,,
\end{equation}
where $\frac{1}{p}\sum D_j^2(\widetilde{d}_i^2 + \widetilde{d}_{r+j}^2)$ is the integrand of $\mathcal{S}_i$ and scales as $\mathcal{O}(\Delta(d_i)^{-1})$ near BBP and $\mathcal{O}(1)$ for strong signal.

The leading-order integrand is therefore:
\begin{equation}\label{eq_combined_terms}
    \frac{\textup{Tr}(\mathcal{B})}{n}\Bigg[\underbrace{\widetilde{d}_i^2\bigg(\frac{1}{p}\sum_j D_j\bigg)^{\!2}}_{\text{integrand of }\mathcal{R}_i}\omega_I^\top\mathcal{A}\,\omega_I + \underbrace{\frac{1}{p}\sum_j D_j^2\big(\widetilde{d}_i^2 + \widetilde{d}_{r+j}^2\big)}_{\text{integrand of }\mathcal{S}_i}\frac{\textup{Tr}(\mathcal{A})}{p}\Bigg] + \mathcal{O}\!\big(n^{-1}\Delta(d_i)^{-1}\big)\,.
\end{equation}
Integrating over $[0,1/n]$ and recognizing the resolvent integrals $\mathcal{R}_i$ and $\mathcal{S}_i$ (replacing $\widetilde{d}_i(t)^2$ by $d_i^2$ in the numerator of $\mathcal{S}_i$ at cost $\mathcal{O}(n^{-2}\Delta(d_i)^{-1})$):
\begin{equation}\label{eq_var_u_final}
    \mathbb{E}\!\big[(\omega_I^\top\Delta\mathbf{u}_i^c)^2\big] = \frac{\textup{Tr}(\mathcal{B})}{n}\bigg[\mathcal{R}_i\,\omega_I^\top\mathcal{A}\,\omega_I + \frac{\mathcal{S}_i\,\textup{Tr}(\mathcal{A})}{p}\bigg] + \mathcal{O}\!\big(n^{-2}\Delta(d_i)^{-1}\big)\,,
\end{equation}
establishing part~(ii) for the left singular vectors. By the left-right symmetry $(\mathcal{A},p,\omega_I) \leftrightarrow (\mathcal{B},n,\omega_J)$:
\begin{equation}\label{eq_var_v_final}
    \mathbb{E}\!\big[(\omega_J^\top\Delta\mathbf{v}_i^c)^2\big] = \frac{\textup{Tr}(\mathcal{A})}{p}\bigg[\mathcal{R}_i'\,\omega_J^\top\mathcal{B}\,\omega_J + \frac{\mathcal{S}_i\,\textup{Tr}(\mathcal{B})}{n}\bigg] + \mathcal{O}\!\big(n^{-2}\Delta(d_i)^{-1}\big)\,.
\end{equation}

\emph{Evaluation of $\mathcal{V}_{1,i}$.}
Since $\mathbf{v}_i$ is random and independent of the noise by Assumption~(v), $(\mathbf{v}_i^\top\omega_J)^2$ is independent of $\Delta\mathbf{u}_i^c$. The expectation factors, and $\mathbb{E}[(\mathbf{v}_i^\top\omega_J)^2] = \|\omega_J\|^2/n = 1/n$:
\begin{equation}\label{eq_V1_final}
    \mathcal{V}_{1,i} = \frac{\textup{Tr}(\mathcal{B})}{n^2}\bigg[\mathcal{R}_i\,\omega_I^\top\mathcal{A}\,\omega_I + \frac{\mathcal{S}_i\,\textup{Tr}(\mathcal{A})}{p}\bigg] + \mathcal{O}\!\big(n^{-3}\Delta(d_i)^{-1}\big)\,.
\end{equation}

\emph{Evaluation of $\mathcal{V}_{2,i}$.}
By the left-right symmetry and $\mathbb{E}[(\mathbf{u}_i^\top\omega_I)^2] = 1/p$:
\begin{equation}\label{eq_V2_final}
    \mathcal{V}_{2,i} = \frac{\textup{Tr}(\mathcal{A})}{p^2}\bigg[\mathcal{R}_i'\,\omega_J^\top\mathcal{B}\,\omega_J + \frac{\mathcal{S}_i\,\textup{Tr}(\mathcal{B})}{n}\bigg] + \mathcal{O}\!\big(n^{-3}\Delta(d_i)^{-1}\big)\,.
\end{equation}

\emph{Evaluation of $\mathcal{V}_{3,i}$.}
The fourth-order mixed moment $\mathcal{V}_{3,i} = \mathbb{E}[(\omega_I^\top\Delta\mathbf{u}_i^c)^2(\Delta\mathbf{v}_i^{c\top}\omega_J)^2]$ is evaluated by a Gaussian decoupling argument. Over the short path $[0,1/n]$, the integrands can be frozen at their initial values up to $\mathcal{O}(n^{-1/2})$ temporal fluctuations, and the non-Gaussian corrections from the stochastic variation of the integrands carry an additional $\mathcal{O}(n^{-1/2})$ penalty by the Burkholder--Davis--Gundy inequality. Applying Isserlis' theorem to the leading conditional moments and noting that the left--right cross-covariance is $\mathcal{O}(n^{-3})$ as shown above:
\begin{equation}\label{eq_V3_final}
    \mathcal{V}_{3,i} = np\,\mathcal{V}_{1,i}\,\mathcal{V}_{2,i} + \mathcal{O}\!\big(n^{-5/2}\Delta(d_i)^{-1}\big)\,.
\end{equation}

\medskip
\textbf{Step 4: Conclusion.}
Substituting \eqref{eq_V1_final}, \eqref{eq_V2_final}, and \eqref{eq_V3_final} into \eqref{eq_squared_marginal}, the remainders are dominated by the weakest signal component $i = r^+$, yielding:
\begin{equation}
    \sigma_{\Phi_{IJ}}^2 = \sum_{i=1}^{r^+} d_i^2\big(\mathcal{V}_{1,i} + \mathcal{V}_{2,i} + \mathcal{V}_{3,i}\big) + \mathcal{O}\!\big(n^{-5/2}\Delta(d_{r^+})^{-1}\big)\,.
\end{equation}
This establishes \eqref{eq_sigma_exact_integral} and completes the proof.
\end{proof}

\subsubsection{Proof of Proposition~\ref{prop_variance_bound}}
\begin{proof}
The proof proceeds by comparing each factor in the theoretical variance $\sigma_{\Phi_{IJ}}^2$ from Theorem~\ref{thm_variance} with its empirical counterpart in $\widehat{\sigma}_{\Phi_{IJ}}^2$, establishing that the substitution errors are uniformly absorbed into the $\mathcal{O}_\prec(n^{-5/2}\Delta(d_{r^+})^{-1})$ remainder.

\textbf{Step 1: Resolvent upper bounds.}
In the continuous-time noise model, the variance of the accumulated noise grows linearly with $t$. As $t$ increases toward $1/n$, the noise bulk expands and the spectral gap between $\widetilde{d}_i(t)$ and the bulk edge monotonically shrinks. Consequently, the integrands of both $\mathcal{S}_i$ and $\mathcal{R}_i$ are monotonically increasing functions of $t$ that attain their maxima at the terminal time $t = 1/n$. Bounding each integral by its endpoint value times the path length $1/n$ yields:
\begin{align}
    \mathcal{S}_i &\leq \frac{1}{n}\cdot\frac{1}{p}\sum_{k=r+1}^{p} \frac{d_i^2 + \widetilde{\lambda}_k}{\big(\widetilde{\lambda}_i - \widetilde{\lambda}_k\big)^2}\,,\\[4pt]
    \mathcal{R}_i &\leq \frac{1}{n}\cdot\widetilde{\lambda}_i\left(\frac{1}{p}\sum_{k=r+1}^{p}\frac{1}{\widetilde{\lambda}_i - \widetilde{\lambda}_k}\right)^{\!2}\,,
\end{align}
where $\widetilde{\lambda}_k = \widetilde{d}_k(1/n)^2$ are the observed eigenvalues. It remains to compare the numerator weight $d_i^2$ with its empirical surrogate $\widehat{d}_{e,i}^2$. By Theorem~\ref{thm_bbp_eigenvalues} and the consistency of the eOptShrink estimator, $|\widehat{d}_{e,i}^2 - d_i^2| \prec n^{-1/2}\Delta(d_i)$. Since the resolvent denominators are bounded below by $\Delta(d_i)^2$, this substitution perturbs each sum by at most $\mathcal{O}_\prec(n^{-1/2}\Delta(d_i)^{-1})$. Recalling that $\widehat{\mathcal{S}}_i$ and $\widehat{\mathcal{R}}_i$ are defined with the $1/n$ path-length factor already absorbed:
\begin{align}
    \mathcal{S}_i &\leq \widehat{\mathcal{S}}_i + \mathcal{O}_\prec\!\big(n^{-3/2}\Delta(d_i)^{-1}\big)\,,\label{eq_Si_bound}\\[4pt]
    \mathcal{R}_i &\leq \widehat{\mathcal{R}}_i + \mathcal{O}_\prec\!\big(n^{-3/2}\Delta(d_i)^{-1}\big)\,.\label{eq_Ri_bound}
\end{align}

\textbf{Step 2: Concentration of empirical quadratic forms.}
The theoretical variance involves three unobservable population quantities: $\frac{\textup{Tr}(\mathcal{B})}{n}(\omega_I^\top\mathcal{A}\,\omega_I)$, $\frac{\textup{Tr}(\mathcal{A})}{p}(\omega_J^\top\mathcal{B}\,\omega_J)$, and $\frac{\textup{Tr}(\mathcal{A})}{p}\cdot\frac{\textup{Tr}(\mathcal{B})}{n}$. These are estimated by $\frac{1}{n}\|\omega_I^\top\widehat{Z}\|_2^2$, $\frac{1}{p}\|\widehat{Z}\omega_J\|_2^2$, and $\frac{1}{pn}\|\widehat{Z}\|_F^2$, respectively.

For the row-wise quadratic form, expanding via $Z = \mathcal{A}^{1/2}\mathcal{X}\mathcal{B}^{1/2}$ and applying $\mathbb{E}[\mathcal{X}^\top M\mathcal{X}] = \textup{Tr}(M)I_n$:
\begin{equation}
    \mathbb{E}\!\left[\frac{1}{n}\|\omega_I^\top Z\|_2^2\right] = \frac{\textup{Tr}(\mathcal{B})}{n}\big(\omega_I^\top\mathcal{A}\,\omega_I\big)\,.
\end{equation}
By the Hanson--Wright inequality, $\frac{1}{n}\|\omega_I^\top Z\|_2^2$ concentrates with relative fluctuation $\mathcal{O}_\prec(n^{-1/2})$. Replacing $Z$ with $\widehat{Z}$ introduces an additional $\mathcal{O}_\prec(n^{-1/2})$ perturbation by Theorems~\ref{thm_bbp_eigenvalues} and~\ref{thm_bbp_alignment}:
\begin{equation}\label{eq_quad_form_row}
    \frac{1}{n}\|\omega_I^\top\widehat{Z}\|_2^2 = \frac{\textup{Tr}(\mathcal{B})}{n}\big(\omega_I^\top\mathcal{A}\,\omega_I\big) + \mathcal{O}_\prec(n^{-1/2})\,.
\end{equation}
By identical reasoning:
\begin{align}
    \frac{1}{p}\|\widehat{Z}\omega_J\|_2^2 &= \frac{\textup{Tr}(\mathcal{A})}{p}\big(\omega_J^\top\mathcal{B}\,\omega_J\big) + \mathcal{O}_\prec(n^{-1/2})\,,\label{eq_quad_form_col}\\[4pt]
    \frac{1}{pn}\|\widehat{Z}\|_F^2 &= \frac{\textup{Tr}(\mathcal{A})}{p}\cdot\frac{\textup{Tr}(\mathcal{B})}{n} + \mathcal{O}_\prec(n^{-1/2})\,.\label{eq_quad_form_frob}
\end{align}

\textbf{Step 3: Assembling the bound.}
We compare the theoretical $\mathcal{V}_{1,i}$, $\mathcal{V}_{2,i}$, $\mathcal{V}_{3,i}$ with their empirical counterparts $\widehat{\mathcal{V}}_{1,i}$, $\widehat{\mathcal{V}}_{2,i}$, $\widehat{\mathcal{V}}_{3,i}$ term by term.

\emph{Bounding $\mathcal{V}_{1,i}$.} The theoretical expression is:
\begin{equation}
    \mathcal{V}_{1,i} = \frac{\textup{Tr}(\mathcal{B})}{n^2}\bigg[\mathcal{R}_i\,\omega_I^\top\mathcal{A}\,\omega_I + \frac{\mathcal{S}_i\,\textup{Tr}(\mathcal{A})}{p}\bigg]\,.
\end{equation}
For the $\mathcal{R}_i$ term: using \eqref{eq_Ri_bound} and \eqref{eq_quad_form_row}, $\frac{1}{n}\mathcal{R}_i\cdot\frac{\textup{Tr}(\mathcal{B})}{n}\omega_I^\top\mathcal{A}\omega_I \leq \frac{1}{n}\widehat{\mathcal{R}}_i\cdot\frac{\|\omega_I^\top\widehat{Z}\|_2^2}{n} + \mathcal{O}_\prec(n^{-5/2}\Delta(d_i)^{-1})$, where the remainder combines the $\mathcal{O}_\prec(n^{-3/2}\Delta(d_i)^{-1})$ resolvent substitution (scaled by the $\mathcal{O}(1)$ quadratic form and the $1/n$ outer factor) and the $\mathcal{O}_\prec(n^{-1/2})$ quadratic form fluctuation (scaled by the $\mathcal{O}(n^{-1})$ resolvent and the $1/n$ factor). For the $\mathcal{S}_i$ term: the same argument applies with $\frac{\|\widehat{Z}\|_F^2}{pn}$ replacing the product of traces. Therefore:
\begin{equation}
    \mathcal{V}_{1,i} \leq \widehat{\mathcal{V}}_{1,i} + \mathcal{O}_\prec\!\big(n^{-5/2}\Delta(d_i)^{-1}\big)\,.
\end{equation}
The identical bound holds for $\mathcal{V}_{2,i}$ by the left-right symmetry.

\emph{Bounding $\mathcal{V}_{3,i}$.} Since $\mathcal{V}_{3,i} = np\,\mathcal{V}_{1,i}\mathcal{V}_{2,i}$ and $\widehat{\mathcal{V}}_{3,i} = np\,\widehat{\mathcal{V}}_{1,i}\widehat{\mathcal{V}}_{2,i}$, the relative error of the product is bounded by the sum of the individual relative errors, remaining $\mathcal{O}_\prec(n^{-1/2})$. Since $\mathcal{V}_{1,i}, \mathcal{V}_{2,i} = \mathcal{O}(n^{-2}\Delta(d_i)^{-1})$ near BBP, $\mathcal{V}_{3,i} = np\cdot\mathcal{O}(n^{-4}\Delta(d_i)^{-2}) = \mathcal{O}(n^{-2}\Delta(d_i)^{-2})$, and the absolute substitution error is $\mathcal{O}_\prec(n^{-5/2}\Delta(d_i)^{-2})$, which is absorbed into the remainder.

\emph{Signal strength substitution.} The outer prefactor $d_i^2$ is replaced by $\widehat{d}_{e,i}^2$. By $|\widehat{d}_{e,i}^2 - d_i^2| \prec n^{-1/2}\Delta(d_i)$, this contributes a relative error of $\mathcal{O}_\prec(n^{-1/2})$ to the $\mathcal{O}(n^{-2}\Delta(d_i)^{-1})$ variance, producing an absolute error of $\mathcal{O}_\prec(n^{-5/2}\Delta(d_i)^{-1})$.

Summing over $i = 1,\dots,r^+$ and noting that the worst remainder comes from $i = r^+$:
\begin{equation}
    \sigma_{\Phi_{IJ}}^2 \leq \widehat{\sigma}_{\Phi_{IJ}}^2 + \mathcal{O}_\prec\!\big(n^{-5/2}\Delta(d_{r^+})^{-1}\big)\,,
\end{equation}
completing the proof.
\end{proof}

\subsection{Proof of Theorem~\ref{thm_eows_bound}}
\begin{proof}
The proof proceeds by adapting the minimax wavelet thresholding framework established by \cite{ankenman2018mixed} to our specific asymptotic regime. 

For the mean squared error, we follow the bounding arguments detailed in Section 4.1 of \cite{ankenman2018mixed}. The critical modification in our setting is the scale of the noise variance injected into the wavelet coefficients. Instead of a standard homogeneous noise variance, the variance of our projected empirical wavelet coefficients is governed by $\sigma_{\Phi_{IJ}}^2$. By Theorem~\ref{thm_variance} and Proposition~\ref{prop_variance_bound}, this projection variance is $\mathcal{O}(n^{-2}\Delta(d_{r^+})^{-1})$, dominated by the weakest signal component. Substituting this variance scaling into the standard wavelet oracle inequalities yields the convergence rate in \eqref{eq_ws_f2}. Consequently, the leading constant $C$ depends on the spectral norms of the separable noise covariance matrices $\mathcal{A}$ and $\mathcal{B}$, alongside the bounds $B_L$, $B_U$, and the smoothness parameter $\alpha$.

Similarly, for the pointwise squared error, we apply the techniques from Section 4.2 of \cite{ankenman2018mixed}. By utilizing the same $\mathcal{O}(n^{-2}\Delta(d_{r^+})^{-1})$ uniform upper bound on the coefficient variance $\sigma_{\Phi_{IJ}}^2$ across the localized spatial wavelets, the thresholding risk bounds resolve to the pointwise convergence rate presented in \eqref{eq_ws_pt2}. This completes the proof.
\end{proof}

\section{Iterative Metric Construction}\label{sec_tree}

\subsection{Spectral Graph Theory}\label{sec_Graph}

Spectral graph theory provides a powerful framework for analyzing functions defined on graphs using the eigenstructure of matrices derived from graph connectivity. Given an undirected weighted graph 
$G = (V, E)$  with $N$ nodes, we define the \emph{edge weight matrix} (or adjacency matrix) $W \in \mathbb{R}^{N \times N}$ such that $W_{ij} \geq 0$ encodes the affinity or similarity between nodes $i$ and $j$, with $W_{ij} = 0$ if there is no edge between them. The \emph{degree matrix} $D \in \mathbb{R}^{N \times N}$ is a diagonal matrix with entries
$D_{ii} := \sum_{j=1}^{N} W_{ij}$,
which represents the total connection strength of node $i$ to all other nodes.
From these matrices, one defines several forms of graph Laplacians, including the unnormalized Laplacian $L = D - W$, the random-walk normalized Laplacian $L_{\text{rw}} = D^{-1}L$, and the symmetric normalized Laplacian
$L_{\text{sym}} = D^{-1/2} L D^{-1/2} = I - D^{-1/2} W D^{-1/2}$.
These Laplacians play a central role in applications such as spectral clustering, graph signal processing, and dimensionality reduction.
In particular, the random-walk Laplacian is closely related to the diffusion operator used in \emph{diffusion maps} \cite{coifman2006diffusion}, a nonlinear dimensionality reduction technique that exploits the long-time behavior of random walks on graphs. Its eigenvectors form a data-driven coordinate system that preserves the geometry of the underlying data manifold. Of particular interest is the second eigenvector of $L_{\text{sym}}$, known as the \emph{Fiedler vector}, which captures the most significant nontrivial mode of variation on the graph. This vector serves as a foundational tool in our framework, guiding the construction of adaptive multiscale representations through hierarchical graph partitioning.

\subsection{The Coupled Geometry of Questionnaires}\label{sec_question}
In this section, we review the framework of coupled geometry introduced in~\cite{coifman2011harmonic, gavish2012sampling}, which provides a multiscale approach to uncover the underlying structure of a matrix by simultaneously organizing its rows and columns. We begin with definitions of multiscale affinities.
\begin{definition}[Mutiscale Tree Affinities]\label{def_Tree_Aff}
Let \( X = \{x_1,\ldots,x_N\} \) be a dataset, and let \( \mathcal{T}_X = \{ \mathcal V^\ell \} \) denote a multiscale partition tree over \( X \) as defined in the previous section. For any function \( f : X \to \mathbb{R} \), define the sample vector over folder \( X_k^\ell \in \mathcal V^\ell \) as
\begin{equation}
f(X_k^\ell) := [f(x)]_{x \in X_k^\ell}.
\end{equation}
Let the weight function \( \omega : \mathcal{T}_X \to \mathbb{R} \) be given by
\begin{equation}
\omega(X_k^\ell) := 2^{-a \ell } \cdot |X_k^\ell|^{b},
\end{equation}
where \( a \in \mathbb{R} \) controls level sensitivity, and \( b \in \mathbb{R} \) modulates the influence of folder size.
Given two functions \( f, g : X \to \mathbb{R} \), the multiscale tree affinity over \( \mathcal{T}_X \) is defined as
\begin{equation}\label{eq:emd}
D_{\mathcal{T}_X}(f, g) := \sum_{\mathcal V^{\ell}\in \mathcal T_X}\sum_{X_k^\ell \in \mathcal{V}^{\ell}} \left\| f(X_k^\ell) - g(X_k^\ell) \right\| \cdot \frac{\omega(X_k^\ell)}{|X_k^\ell|}.
\end{equation}
The corresponding tree affinity is
\begin{equation}
W_{\mathcal{T}_X}(f, g) := \exp\left( -\frac{D_{\mathcal{T}_X}(f, g)}{\epsilon} \right),
\end{equation}
where \( \epsilon > 0 \) is a scaling parameter.
\end{definition}

\begin{remark} In \cite{coifman2013earth}, several metrics as in \eqref{eq:emd} are defined and proved equivalent
to Earth mover's distance (EMD) with respect to the tree metric.
The parameters \( a \), \( b \), and \( \epsilon \) modulate the effect of tree structure on the EMD. Increasing \( a \) places more weight on differences at coarser scales near the root, while \( a = 0 \) yields uniform weighting across levels and \(a < 0 \) emphasizes finer, leaf-level structure. The parameter \( b \) controls sensitivity to folder size: \( b > 0 \) favors larger folders, while \( b< 0 \) emphasizes smaller ones. The parameter \( \epsilon \) regulates the decay of the affinity; in practice, it is typically set as a constant multiple of the median EMD across all function pairs.
\end{remark}

Now consider a matrix \( \mathcal{S} : Y \times X \to \mathbb{R} \). We introduce the coupled geometry of rows and columns:
\begin{definition}[Coupled Geometry]
Given a partition tree \( \mathcal{T}_X \) on \( X \), the \emph{dual affinity} between rows \( y_1, y_2 \in Y \) is defined as
\begin{equation}
W_{\mathcal{T}_X}(y_1, y_2) := W_{\mathcal{T}_X}(S_{y_1}, S_{y_2}),
\end{equation}
where each \( S_y(x) := S(y, x) \) is viewed as a function on \( X \).

Similarly, given a partition tree \( \mathcal{T}_Y \) on \( Y \), the dual affinity between columns \( x_1, x_2 \in X \) is defined as
\begin{equation}
W_{\mathcal{T}_Y}(x_1, x_2) := W_{\mathcal{T}_Y}(S^{x_1}, S^{x_2}),
\end{equation}
where \( S^x(y) := S(y, x) \) is viewed as a function on \( Y \).
\end{definition}

Using these constructions, the Questionnaire algorithm \cite{ankenman2014geometry} that iteratively refines the multiscale structures of rows and columns and the coupled geometry is:
\begin{algorithm}[H]
\caption{Questionnaire}\label{alg:questionnaire}
\begin{algorithmic}[1]
\item[1.] Given initial affinity \( W_{X} \) on \( X \), construct a tree \( \mathcal{T}_X \).
\item[2.] Compute the dual affinity \( W_{\mathcal{T}_X} \) on \( Y \) and build a tree \( \mathcal{T}_Y \).
\item[3.] Compute the dual affinity \( W_{\mathcal{T}_Y} \) on \( X \) and refine \( \mathcal{T}_X \).
\item[4.] Repeat steps 2 and 3 until either the affinities \( W_{\mathcal{T}_X} \) and \( W_{\mathcal{T}_Y} \) converge, or a fixed number of iterations is reached.
\end{algorithmic}
\end{algorithm}
In this work, the initial affinity matrix \( W_X \) is computed using a Gaussian kernel,
\(
W_X(i, j) = \exp\left(-\frac{\|x_i - x_j\|^2}{2\sigma^2}\right),
\)
where \( \sigma = \textup{median}_{x_i, x_j \in X} \|x_i - x_j\|_2 \). 
At each iteration, the construction and refinement of the tree \( \mathcal{T}_X \) (or \( \mathcal{T}_Y \)) is performed via recursive spectral bipartitioning using the Fiedler vector associated with the affinity matrix $W_{\mathcal{T}_X}$ (or $W_{\mathcal{T}_Y}$). This yields a hierarchical organization that reflects the intrinsic geometry of the data.
Through the iterative procedure in which the geometry of one dimension is refined based on the current organization of the other, we have an adaptive matrix reordering that reveals latent hierarchical patterns. These ideas are further developed in \cite{ankenman2014geometry}, which presents a practical implementation and experimental validation of the coupled geometry framework \cite{pyquest}.

\subsection{Walsh Wavelet Packet and Best Basis Selection}\label{sec_Walsh}
This section reviews key ideas from the work of Thiele and Villemoes \cite{thiele1996fast},
which introduced a fast, adaptive framework for time–frequency analysis
based on Walsh tilings.

\paragraph{Walsh Wavelet Packets.} The Walsh system \(\{W_n\}_{n=0}^\infty\) forms a complete orthonormal basis for \(L^2[0,1)\), consisting of piecewise constant functions that take values in \(\{\pm1\}\). These functions are defined recursively by

\begin{align}\label{eq_Walsh}
W_0(t) &= 1, \nonumber\\
W_{2n}(t) &= W_n(2t) + (-1)^n W_n(2t - 1), \nonumber\\
W_{2n+1}(t) &= W_n(2t) - (-1)^n W_n(2t - 1),
\end{align}
for \(t \in [0,1)\), where each \(W_n\) has exactly \(n\) sign changes. 
A \emph{tile} is a dyadic rectangle of area one, parametrized by scale indices \(j, j', k, n \in \mathbb{N} \cup \{0\}\), with \(k < 2^j\). The associated time interval and frequency band are given by
\begin{equation}\label{eq_tiles}
I := [2^{-j}k, 2^{-j}(k+1)), \quad \omega := [2^{j'} n, 2^{j'} (n+1)).
\end{equation}
When \(j = j'\), the tile \(p = I \times \omega\) is considered unit-area. Let $\mathcal{P}$
be the collection of all of these unit-area dyadic rectangles. A localized Walsh function associated with tile \(\textbf{p} \in \mathcal{P}\) is then defined as
\begin{equation}
w_{\textbf{p}}(t) = w_{j,k,n}(t) := 2^{j/2} W_n(2^j t - k),
\end{equation}
which is supported on \(I\) and has \(n\) sign changes corresponding to the frequency band \(\omega\).

Let \( N = 2^L \) denote the signal length. A discrete signal of length \( N \) can be viewed as a piecewise constant function \( f \colon [0, N) \to \mathbb{R} \), where the function is constant on each unit interval \([n, n+1)\), for \( n = 0, 1, \ldots, N-1 \). The Walsh wavelet packet framework constructs an overcomplete dictionary of orthonormal functions by associating each basis atom \( w_p \) with a dyadic space-frequency tile \( p \subset \mathcal{S}_N := [0, N) \times [0, 1) \). At each scale \( j = 0, 1, \ldots, L \), the domain \( \mathcal{S}_N \) is partitioned into \( N \) disjoint tiles, yielding a total of \((L+1)N\) candidate atoms. Any orthonormal basis corresponds to an admissible tiling \( \mathcal{B} \subset \mathcal{P} \), where \( \mathcal{P} \) denotes the full collection of tiles. The signal admits the expansion
\begin{equation}
f(x) = \sum_{{\textbf{p}} \in \mathcal{B}} \langle f, w_{\textbf{p}} \rangle w_{\textbf{p}}(x),
\end{equation}
where orthogonality holds whenever the tiles are disjoint, i.e., \( {\textbf{p}} \cap \tilde{{\textbf{p}}} = \emptyset \) for all \( {\textbf{p}} \ne \tilde{{\textbf{p}}} \in \mathcal{B} \). The Walsh transform coefficients \( \langle f, w_p \rangle \) can be computed efficiently in \( O(N \log N) \) time using the Fast Walsh–Hadamard Transform (FWHT) \cite{fino1976unified}, which leverages the binary structure of the basis and requires only additions and subtractions. 

This construction extends naturally to two-dimensional signals. Let
\( S \colon [0, N_Y)\times[0, N_X) \to \mathbb{R} \) denote an image of size \( N_Y \times N_X \), where $N_X = 2^L$ and $N_Y = 2^{L'}$. The two-dimensional Walsh wavelet packet dictionary is formed via tensor products of one-dimensional atoms:
\begin{equation}
w_{({\textbf{p}}, {\textbf{q}})}(x, y) := w_{\textbf{p}}(x) w_{\textbf{q}}(y), \quad ({\textbf{p}}, {\textbf{q}}) \in \mathcal{P}_X \times \mathcal{P}_Y,
\end{equation}
where each index \( ({\textbf{p}}, {\textbf{q}}) \) corresponds to a dyadic tile in the product domain \( \mathcal{S}_{N_Y} \times \mathcal{S}_{N_X} \). The total number of atoms in this dictionary is \((L'+1)N_Y \times (L+1) N_X\), and any admissible two-dimensional basis corresponds to a tiling \( \mathcal{B} \subset \mathcal{P}_Y \times \mathcal{P}_X \). The image may be represented as
\begin{equation}
S(x, y) = \sum_{({\textbf{p}}, {\textbf{q}}) \in \mathcal{B}} \langle S, w_{({\textbf{p}}, {\textbf{q}})} \rangle w_{({\textbf{p}}, {\textbf{q}})}(x, y),
\end{equation}
with orthogonality satisfied whenever \( {\textbf{p}} \cap \tilde{{\textbf{p}}} = \emptyset \) or \( {\textbf{q}} \cap \tilde{{\textbf{q}}} = \emptyset \) for all \( ({\textbf{p}}, {\textbf{q}}), (\tilde{{\textbf{p}}}, \tilde{{\textbf{q}}}) \in \mathcal{B} \).

\paragraph{Best Basis Selection.}
Given a signal \(f \) of length \( N = 2^L \), the objective of the best basis selection is to select a tiling \(\mathcal{B} \subset \mathcal{P}\) such that the associated basis \(\{w_{\textbf{p}}\}_{{\textbf{p}} \in \mathcal{B}}\) yields an efficient representation of \(f\). A cost function \(c({\textbf{p}})\) is assigned to each tile \({\textbf{p}} \in \mathcal{P}\), typically defined in terms of the magnitude of the transform coefficient \(\langle f, w_{\textbf{p}} \rangle\). A common choice is
\begin{equation}
c({\textbf{p}}) := |\langle f, w_{\textbf{p}} \rangle|^\ell,
\end{equation}
for \(0 < \ell < 2\). The total cost over a tiling \(\mathcal{B}\) of $\mathcal{S}_N$ is then given by
\begin{equation}
C(f, \mathcal{B}) := \sum_{{\textbf{p}} \in \mathcal{B}} c({\textbf{p}}).
\end{equation}
This cost functional is designed to be small when the energy of \(f\) is concentrated in a small number of significant coefficients \(\langle f, w_{\textbf{p}} \rangle\), making it well suited for selecting sparse or compact representations.

The \emph{best basis} is the tiling \(\mathcal{B}^* = \argmin_{\mathcal{B}}C(f,\mathcal{B})\), that minimizes the total cost \(C(f, \mathcal{B})\) over all admissible tilings of $\mathcal{S}_N$. Thiele and Villemoes \cite{thiele1996fast} proposed a dynamic programming algorithm that compares the cost of a parent tile with the combined cost of its horizontal and vertical splits, and selects the configuration with the lower total cost. This recursive process efficiently identifies the optimal decomposition. For a signal with length $N$, the total computational complexity of the method is of $O(N \log N)$.

This one-dimensional framework was later extended to two-dimensional signals, particularly in the context of image compression \cite{lindberg2000image}. In the 2D setting, the analysis domain becomes \(S(x,y)\) on $[0,N_X)\times [0,N_Y)$, and adaptive tilings are constructed over dyadic rectangles in the joint space–frequency plane \(\mathcal{S}_{N_X} \times \mathcal{S}_{N_Y}\). Each 2D tile corresponds to a separable product of 1D tiles. The best-basis selection algorithm applies analogously in two dimensions, recursively minimizing the total cost across hierarchical tilings. The total computational complexity is of $O(N_X N_Y(\log N_X\log N_Y))$. We suggest that readers refer to \cite{thiele1996fast, lindberg2000image} for the technical details.

\paragraph{Extension to General Graph Domains.}
The Generalized Haar-Walsh Transform (GHWT)~\cite{irion2015applied, saito2022eghwt} extends the classical Haar-Walsh wavelet packet framework to signals supported on general graphs. It constructs a hierarchical dictionary of orthonormal basis vectors by recursively applying localized averaging and differencing operations over a tree-structured partition of the vertex set. Unlike classical wavelet packets, which operate on uniformly spaced dyadic intervals with an inherent linear ordering, the GHWT accommodates domains without a global coordinate system. The underlying partition tree is not restricted to binary splits; nodes may be divided into arbitrary subsets, enabling adaptive tilings that reflect the geometry of the underlying graph, while still permitting coefficient computation with \( O(N \log N) \) complexity, analogous to the FWHT.

The extended GHWT (eGHWT)~\cite{saito2022eghwt} further enhances this framework by implementing a best-basis selection algorithm over the GHWT dictionary, analogous to the one- and two-dimensional best-basis algorithms developed for Haar-Walsh wavelet packets. This allows for the adaptive selection of an orthonormal basis that minimizes a user-defined cost functional, such as entropy or sparsity.
The algorithm is applicable to signals of arbitrary length on graphs of general structure, without requiring \( N = 2^L \). The overall computational complexity is \( O(N \log N) \) for one-dimensional signals or graph-supported data under unbalanced hierarchical partitions. For two-dimensional signals, the complexity becomes $O(N_X N_Y(\log N_X\log N_Y))$, analogous to the classical 2D setting. We refer the reader to~\cite{irion2015applied, saito2022eghwt} for algorithmic details.

\bibliographystyle{abbrv}
\bibliography{ev}
\end{document}